\documentclass[12pt]{amsart}

\newtheorem{coro}{Corollary}[section]
\newtheorem{defi}{Definition}[section]
\newtheorem{prop}{Proposition}[section]
\newtheorem{theo}{Theorem}
\newtheorem{lemm}{Lemma}[section]
\newtheorem{ques}{Question}[section]
\newtheorem{rema}{Remark}
\newtheorem{clai}{Claim}

\def\NN{{\mathbb N}}

\def\ZZ{{\mathbb Z}}

\def\La{{\Lambda}}

\def\Om{{\Omega}}

\def\cG{{\mathcal G}} 
\def\cM{{\mathcal M}}

\def\cN{{\mathcal N}}

\def\cO{{\mathcal O}} 
\def\cU{{\mathcal U}}

\def\cP{{\mathcal P}} 
\def\cV{{\mathcal V}}
  
\def\cK{{\mathcal K}} 
\def\cQ{{\mathcal Q}} 
\def\cW{{\mathcal W}}

\def\cR{{\mathcal R}}

\def\rG{\mathrm{G}}
\def\rO{\mathrm{O}}
\def\rR{\mathrm{R}}
\def\fP{\mathfrak{P}}
\def\diam{\mathrm{diam}}
\def\cdiam{\mathrm{cdiam}}
\def\supp{\mathrm{supp}}
\def\MEx{{}_{\mathrm{MEx}}}

\title[Aperiodic classes of $C^1$-generic diffeomorphisms]{\bf Aperiodic chain recurrence classes of $C^1$-generic diffeomorphisms}

\begin{document}

\author{Christian Bonatti and Katsutoshi Shinohara}
\begin{abstract}We consider the space 
of $C^1$-diffeomorphims equipped with the $C^1$-topology
on a three dimensional closed manifold. 
It is known that there are open sets in which 
$C^1$-generic diffeomorphisms display uncountably many chain recurrences classes, while only countably many of them may contain periodic orbits. The classes without periodic orbits, called aperiodic classes, are the main subject of this paper.  
The aim of the paper is to show that aperiodic classes of $C^1$-generic diffeomorphisms can exhibit a variety of
topological properties. More specifically, there are $C^1$-generic diffeomorphisms with 
(1) minimal expansive aperiodic classes,
(2) minimal but non-uniquely ergodic aperiodic classes,
(3) transitive but non-minimal aperiodic classes,
(4) non-transitive, uniquely ergodic aperiodic classes.
\end{abstract}

\maketitle

\footnotesize{

}


\today


\section{Introduction}
\subsection{General setting: Wild dynamical systems}
The notion of chaotic dynamical systems goes back to Poincar\'{e} \cite{Po}: He noticed that certain deterministic behavior, governed by simple equations, presents a very complicated behavior if the stable and the unstable manifold of 
a periodic saddle point intersect transversely. This phenomenon, called a \emph{homoclinic intersection} according to his terminology, not only implies that the system is chaotic, but also every system described by slightly
perturbed equations are 
as such. The behavior is \emph{robustly chaotic}. 

In the middle of the twentieth century, 
Anosov and Smale develop a 
theory of a geometric structure called \emph{(uniform) hyperbolicity}. Hyperbolic systems may be chaotic but they are stable: The systems close to a hyperbolic system are obtained 
by looking the system through a small continuous change of coordinates. More precisely, the theory of 
structural stability (see \cite{Ro,Ro2, Sa, Pa, Ma2, Ha}) shows that the notions of stability called \emph{$C^1$-structural stability} (on the all space) or \emph{$C^1$-$\Om$-stability} (restricted to the non-wandering set) are characterized by geometric structures having loquacious names: \emph{the Axiom A + the strong transversality property} and \emph{the Axiom A + the no cycle condition} respectively. These 
geometric structures allow a very precise 
description of the dynamics from topological and ergodic
viewpoints.

However, these structures are so rigid that they can only describe specific regions in the space of dynamical systems. It has been noticed  in the late sixties by Newhouse \cite{Ne} (for the $C^2$-topology on surface diffeomorphisms) and 
Abraham and Smale \cite{AS} (for the $C^1$-topology, in dimension $\geq 3$)
that there are open sets of non-hyperbolic dynamical systems. 
These systems are robustly unstable: Perturbations of the systems change the qualitative behavior of dynamics. 

How can we describe these robustly unstable systems? 
In the non-conservative setting, the global dynamics 
split into several somehow independent pieces
(we refer  \cite{Bo} for a global overview of the space of $C^1$ dynamical systems). The first global 
description as such is given 
(in the very general setting of homeomorphisms of compact metric spaces) by Conley \cite{Co} using the notion of attracting/repelling 
sets: Two recurrent points are not in the same \emph{chain recurrence class} if and only if there is an attracting/repelling
set containing one and only one of these two points. 

Hyperbolic dynamical systems have only finitely many  chain recurrence classes, and these chain recurrence classes are the \emph{homoclinic classes} (closure of the transverse homoclinic intersections) and are the maximal transitive sets. 
For robustly non-hyperbolic generic systems, 
there are two kinds of typical behavior.
\begin{itemize}
\item  Either they have finitely many classes  whose number remains locally constant in the space of all systems (see \cite{Sh,Ma1} for examples). 
They are called \emph{tame dynamical systems}.
\item Or, locally generically (i.e., on a residual subset of a non-empty open set of all the systems) they have infinitely many classes (see \cite{N, BD2,BD1}). They are called 
\emph{wild dynamical systems}.  
\end{itemize}

Indeed, in the space of $C^1$-diffeomoprhisms, 
by \cite{BC} we know that every $C^1$-generic diffeomorphism 
can be classified one of above kinds of dynamical systems. 
This paper focuses on $C^1$-generic wild diffeomorphisms on $3$-manifolds. Recall that for $C^1$-generic dynamical 
systems, Kupka-Smale theorem \cite{Ku,Sm} implies that the periodic orbits are all hyperbolic and 
Pugh's closing lemma \cite{Pu}, the connecting lemma by \cite{Ha} and \cite{BC} imply that the set of hyperbolic
periodic orbits is dense in the chain recurrent set. 
According to \cite{BC}, for $C^1$-generic diffeomorphisms, 
a chain recurrence class containing a periodic orbit coincides 
with its homoclinic class and is a maximal transitive set. 
This property is an important starting point 
for understanding the dynamics inside these classes. 

While the cited results imply that for $C^1$-generic 
tame systems we can always find a hyperbolic periodic point,
$C^1$-generic wild diffeomorphisms may have \emph{aperiodic classes}, i.e., chain recurrence classes with no periodic orbit. 
As far as we know, very few are known about the dynamics inside aperiodic classes. Up to now, the unique known example have the structure called adding machines or odometers. 
However, there are no theoretical result which asserts 
it should be always the case. 

The aim of this paper is to provide a great variety of topological behaviors in aperiodic classes of $C^1$-generic diffeomorphisms on $3$-manifolds. For instance, some of these classes are not transitive, breaking for the first time the equivalence between the notion of chain recurrence classes and maximal transitive sets for $C^1$-generic diffeomorphisms. The following 
is the results we give in this paper:

\begin{theo}\label{t.minimal-non} Let $M$ be 
a closed (compact and without boundary) $3$-manifold. 
Then there is a non-empty $C^1$-open set $\cO$ of 
$\mathrm{Diff}^1(M)$ such that
there is a residual subset $\cG\subset \cO$ 
in which every $f\in\cG$ has
\begin{enumerate}
 \item \label{i.2} an uncountable set of chain recurrence classes which are all minimal and expansive;
\item \label{i.1} an uncountable set of chain recurrence classes which are all minimal but support infinitely many ergodic measures; 
 \item \label{i.4} an uncountable set of chain recurrence classes which are all transitive but contain at least two minimal sets;
 \item \label{i.3} an uncountable set of chain recurrence classes which are all uniquely ergodic but not transitive.
\end{enumerate}
\end{theo}

\subsection{Presentation of the main results}
Let us see the strategy of the construction.
In the series of papers \cite{BS1,BS2,BS3}, we develop a 
technique for expelling a non-trivial hyperbolic set from a partially hyperbolic chain recurrence class satisfying the condition 
called property $(\ell)$, 
by an arbitrarily $C^1$-small perturbation. 
This means that, before the perturbation the set was contained in the class, but after the 
perturbation its continuation is separated from the class by a filtrating set. Conley theory tells us that 
the separation is $C^0$-robust, hence $C^1$-robust. 
In \cite{BS3} this expulsion process is carefully examined and 
we prove that a small perturbation supported in the filtrating set allow us to recover, in the new class, exactly 
the same properties which allow us 
to perform the expulsion, that 
is, property $(\ell)$. 

In the language introduced in \cite{Bo} this means that this property $(\ell)$ is \emph{a viral property}. 
Namely, once we have a chain recurrence class which satisfies
the property $(\ell)$, then by adding an arbitrarily small perturbation we can produce a new chain recurrence class, 
and the same is true for the newly created class. This enables 
us to repeat the creation of distinct chain recurrence classes.
As a result, we can prove that 
nearby $C^1$-generic systems possess nested 
infinite sequence of filtrating sets which has branches at 
each depth. Each infinite sequence corresponds to an aperiodic 
class. As a result, we have 
an uncountable family of aperiodic 
classes.  

This structure appears among 
all known examples of aperiodic classes for $C^1$-generic diffeomorphisms, even for 
those built before the notion of viral property being invented. 
However, up to now, all the known aperiodic classes are \emph{adding machines} which are  
minimal and uniquely ergodic. They are forced to have such  
simple behavior since the filtrating sets at each level 
has simple combinatorial structure 
(periodic attracting/repelling balls) which prohibits
the complexity of the limit dynamics.  

The aim of the series of papers \cite{BS1, BS2, BS3} is 
to establish techniques to expel chain recurrence classes, 
keeping the complexity of each depth and 
furthermore controlling the combinatorics of nested layers.
In our construction, each level is no longer a sequence of 
attracting/repelling balls but a collection of disjoint union 
of finitely many cylinders which behave in a Markovian way
under the iteration of maps.
The structure is called \emph{partially hyperbolic filtrating Markov partitions}. Their behavior can be modeled by means 
of subshifts of finite type. The techniques which we develop 
enables us to continue the expulsion process holding some 
control over combinatorics of each level. By choosing 
convenient sequence of combinatorics, 
we can prove the creation of 
aperiodic classes having prescribed properties.

As we will see, the dynamics inside aperiodic 
classes of $C^1$-generic diffeomorphism can be very rich.  
Our aim is not to produce an exhaustive catalog of the possible dynamics but just to illustrate the great diversity which co-exist in a single generic diffeomorphism. We know
a few more possible properties that the technology presented here may produce, but we would like to leave them as problems 
to be discussed in the future.

Our main result is Theorem~\ref{t.minimal-non}. 
Let us restate it clarifying the prerequisite of the 
open set $\mathcal{O}$. 
Let $f$ be a diffeomorphism on a closed $3$-manifold. 
In \cite{BS2} we introduced the notion of partially 
hyperbolic filtrating Markov partitions, which is a package of information involving the partial hyperbolicity and 
the recurrence of the points. It is a disjoint union of finitely 
many $C^1$-cylinders in $M$ such that $f$ maps them 
Markovian ways. At the same time, the union of the cylinders 
is a filtrating set. Thus, the information of the chain recurrence 
is localized there. For more precise information, see Section~\ref{s.preli}.

Let $p$ be a hyperbolic periodic saddle point of $f$.  
In \cite{BS3}, we discuss a property called 
property $(\ell)$ about a chain recurrence 
class contained in a partially hyperbolic 
filtrating Markov partition, say $\cR$, 
which is related to the partial hyperbolicity around 
the class and the topological conditions about the 
dynamics in the neighborhood.  
For the precise definition of the property $(\ell)$, 
see Section~\ref{s.preliminar}.
Property $(\ell)$ is a $C^1$-robust property 
and in \cite{BS3} we prove that the property $(\ell)$ is a viral property: Arbitrarily small $C^1$-perturbation $g$ of $f$ provides new filtrating Markov partition $\cR_1\subset \cR$ disjoint from the class of $p_g$ (continuation of $p$ for $g$), and a periodic point $q\in\cR_1$ 
such that $g,\cR_1$ and $q$ satisfies $(\ell)$. The aperiodic classes we build are the intersection of a nested sequences of such filtrating Markov partitions.   
However, in order to get a control of the dynamical behavior inside the aperiodic classes, we need to prescribe some topological relation
about how the Markov partition $\cR_1$ are nested inside $\cR$.  
By deliberately investigating this relation, we can check the creation 
of the desired classes. 

In summary, we can restate Theorem~\ref{t.minimal-non} as follows. Below, we say that a filtrating Markov partition 
is \emph{transitive} if given any two rectangles $U$ and $V$ we can find a sequence of rectangles $(W_i)_{i=0,\ldots,n}$ 
such that $W_0=U$, $W_n=V$ and 
$f(W_i) \cap W_{i+1} \neq \emptyset$ holds for
$i=0, \ldots, n-1$. 
\begin{theo}[Restatement of Theorem~\ref{t.minimal-non}] 
\label{t.minimal-non-2}
Let $\cO$ be a $C^1$-open set of diffeomorphisms on a closed $3$-manifold $M$ admitting a transitive, partially hyperbolic filtrating Markov partition $\cU$ containing a hyperbolic
periodic point $p$.  
Suppose that for every $g \in \cO$, $\cU$ is a partially
hyperbolic filtrating Markov partition and
we can define a continuation $p_g$ of $p$ contained in 
$\cU$.  
If the chain recurrence class $[p_g]$ satisfies 
the property $(\ell)$ for every $g$, 
then there is a $C^1$-residual subset $\cG\subset \cO$ such that every $f\in\cG$ has uncountably many 
chain recurrence classes in (1--4) 
of Theorem~\ref{t.minimal-non}.
\end{theo}

\subsection{Questions}
In this subsection, we discuss a few questions related to 
Theorem~\ref{t.minimal-non}.

\subsubsection{Which subshifts can be aperiodic classes?} 
In Theorem~\ref{t.minimal-non-2}, all the properties announced
are realized as aperiodic chain recurrent classes which 
are projective limits of subshifts of finite type. 
Our method consists of expelling hyperbolic subsets repeatedly.
At each step, we expel a part of hyperbolic subset 
in the previous level. By adjusting the combinatorics 
of each level, we obtain the desired condition. 
While we have some control over the choice of the hyperbolic
subset we will expel, we do not have full freedom, 
and we do know if all aperiodic chain recurrent compact subsets of a subshift of finite type can be realized as an aperiodic chain recurrence class of a $C^1$-generic diffeomorphism. 
For instance, the full shift of finitely many symbols 
contains minimal invariant sets 
with positive topological entropy (see for instance \cite{Gr}),
and it is not difficult to realize it as a chain 
recurrence class of a (non-generic) $C^\infty$-diffeomorphism.
We suspect that our method cannot produce such
aperiodic classes for $C^1$-generic diffeomorphisms, 
since at each step  
we need to abandon certain amount 
of complexity which the previous level has. 
Thus, the following question is interesting to investigate.
\begin{ques} 
For $C^1$-generic diffeomorphisms,
what are the topological entropies of aperiodic classes?
Are they always equal to zero? 
\end{ques}

Let us give one more possible direction of future research.
It would be interesting to ask what kind of 
minimal Cantor sets can be realized as aperiodic classes
of $C^1$-generic diffeomorphisms. 
Let us propose a concrete question.
Theorem~\ref{t.minimal-non} produces minimal, expansive
aperiodic classes. We do not know if our example are 
homeomorphic to some known examples. 
A typical example of such dynamical systems
is obtained by considering Denjoy's minimal sets
(see \cite{By} for instance for the investigation of 
Denjoy's minimal sets as subshifts). 
Then the following would be interesting to consider: 
\begin{ques} 
Can one construct aperiodic classes in (1) of 
Theorem~\ref{t.minimal-non} in such a way 
that they are conjugate to Denjoy's minimal sets? 
\end{ques}
See also \cite{GM}, in which some invariants of
minimal Cantor sets are proposed. 

\subsubsection{Adding machines?} 
Theorem~\ref{t.minimal-non} builds aperiodic classes for locally 
$C^1$-generic diffeomorphisms with several different dynamics.
Note that none of these classes are adding machines, 
because they 
violate either minimality, non-expansivity or unique ergodicity. 
However, we do not know if our example is indeed free from 
an adding machine or not. 
Let us discuss the possibility of the existence of 
 adding machine in our setting. 
First, if it exists in the setting of partially hyperbolic filtrating 
Markov partitions, then it must be contained in finitely many periodic contracting center-stable discs: 
Any point in an adding machine 
admit a neighborhood whose orbits remain 
close under iterations. 
However, points in different center stable discs have 
their positive iterates which separate one for the other by a 
uniform distance. 
Thus, the problem is reduced to the 
$C^1$-locally generic existence of 
adding machines in dimension two, 
which is not known until now. Note that 
this problem is tightly related to Smale's conjecture, 
which asserts that Axiom A diffeomorphisms may be $C^1$-dense on surfaces. A positive answer to the conjecture implies the non-existence of locally generic
production of aperiodic classes. In other words, to have 
locally generic existence of adding machines we need to 
have negative answer to Smale's conjecture. 

\begin{ques} 
Do $C^1$-generic diffeomorphisms satisfying the assumption of 
Theorem~\ref{t.minimal-non-2} exhibit 
aperiodic classes which are conjugated to 
adding machines?
\end{ques}

Since the solution of Smale's conjecture seems to be 
out of the reach of our current state of art, it may be difficult to have some 
progress in the $C^1$-topology. Meanwhile,
one may ask a similar question in a higher regularity setting.
One may well expect that the diffeomorphisms in the assumption of
Theorem~\ref{t.minimal-non-2} display homoclinic tangencies in restriction to some periodic center stable discs. 
Thus $C^2$-generic diffeomorphisms in this setting would present adding machine in periodic discs (see for instance \cite[page 33]{BDV}). Thus, under 
an extra assumption, for instance 
assuming the existence of a periodic orbit 
which is volume expansive 
in the center stable direction, it is very likely that one can prove 
the $C^r$-locally generic existence of adding machines. 
So far, we do not have any conjecture about the general case,
and the investigation of such dynamical systems looks very 
interesting.
Even though there are lacks of perturbation techniques 
which are 
available in the $C^1$-topology, the following problem would 
be interesting to pursue:
\begin{ques} 
Let $r \geq 2$.
Do $C^r$-generic diffeomorphisms satisfying the assumption of 
Theorem~\ref{t.minimal-non-2}
exhibit aperiodic classes which are conjugated to 
adding machines?
\end{ques}

Clearly, the following question is also interesting.
\begin{ques} 
Let $r \geq 2$.
Is Theorem~\ref{t.minimal-non-2} true in the $C^r$-setting?
\end{ques}

\subsubsection{The topology of aperiodic classes}
Our construction leads to aperiodic classes which are totally disconnected and contain a Cantor set. The aperiodic classes presented here are intersection of filtrating sets whose connected components have diameters tending to $0$. 
Recall that the set of periodic orbits are dense in the chain recurrent set of any $C^1$-generic diffeomorphism (see \cite{BC}). Thus, an aperiodic class has empty interior.
It is the unique topological property we know about 
these aperiodic classes and it would be interesting 
what kind of topology aperiodic classes can have.  
For instance:
\begin{ques}
Are there locally $C^1$-generic diffeomorphisms having  aperiodic classes containing a non-degenerate continuum
(i.e., a continuum which is not a point)? 
\end{ques}
Note that Ma\~{n}\'{e} showed that a compact metric
space which admits minimal, expansive homeomorphisms 
is homeomorphic to a totally disconnected 
set \cite{Ma-can}, see also \cite{Ar}. Thus this 
question makes sense only for the case (2-4) in 
Theorem~\ref{t.minimal-non-2}.


\subsection{Organization of the paper}
Finally, let us see the organization of this paper. 
In Section~\ref{s.preli} we give some review about 
fundamental notations and the results of
papers \cite{BS1, BS2, BS3}. We keep the review not 
to be comprehensive in order that the readers can 
understand the whole structure of the proof of 
Theorem~\ref{t.minimal-non} without too much meddled
with the technical parts of previous results, which are 
not mandatory for the understanding of the proof of this paper.
In Section~\ref{s.genurgo} we give several topological criterion 
for the confirmation of properties claimed in Theorem~\ref{t.minimal-non}. 
The results presented here are purely topological. 
Thus, we present our results assuming that the ambient space
is just a compact metric space.
In Section~\ref{s.C1} we provide the method to construct 
a nested sequence of filtrating sets satisfying the 
conditions given in Section~\ref{s.genurgo},
by adding an arbitrarily $C^1$-small perturbations
to the systems in the assumption of Theorem~\ref{t.minimal-non-2}. 

\bigskip

{\bf Acknowledgment.}
This work is supported by the JSPS KAKENHI Grant
Numbers 21K03320. 
KS is grateful for the hospitality 
of Institut de Math\'{e}matiques
de Bourgogne of Universit\'{e} de Bourgogne during his visit.


\section{Preliminaries}
\label{s.preli}

In this subsection, we give some fundamental notions of 
dynamical systems and reviews of results of 
\cite{BS1, BS2, BS3} which are used in this paper.

\subsection{Filtrating sets and chain recurrence}
\label{ss.filt}

Let $f$ be a homeomorphism of 
a compact metric space $(X,d)$. A compact set $A\subset X$ is called an \emph{attracting set} if 
$$f(A)\subset \mathring{A},$$
where $\mathring{A}$ is the interior of $A$. 
A \emph{repelling set} is an attracting set for $f^{-1}$.

\begin{rema}\begin{enumerate}\label{r.atre}
\item If $A$ is an attracting set for $f$, then any compact subset $B\subset A$ containing $f(A)$ in its interior (i.e., $f(A)\subset \mathring{B}$) is an attracting set for $f$ and the maximal invariant sets in $A$ and $B$ coincide. 
\item If $A$ is an attracting set for $f$, then $f(A)$ is also 
an attracting set. 
\item If $A$ is an attracting set for $f$, then the complement $X\setminus \mathring A$ is a repelling set.
\end{enumerate}

\end{rema}

\begin{defi} A filtrating set $U$ is an intersection of an attracting set $A$ and a repelling set $R$.
\end{defi}

Let $f, g \in \mathrm{Homeo}(X)$, where
$\mathrm{Homeo}(X)$ denotes the group of homeomorphisms
of $X$. By
$\supp (g, f)$ we denote the closure of the 
set $\{ x \in M \mid f(x) \neq g(x) \}$ and call it 
the \emph{support} of $g$ with respect to $f$.

\begin{lemm}\label{l.filtrating} 
Let $f \in \mathrm{Homeo}(X)$, 
$O\subset X$ be a compact subset 
and $g \in \mathrm{Homeo}(X)$ such that 
$\supp (g, f) \subset O$. Then, for $O$, 
being an attracting set for $f$ or for $g$ are equivalent.  
The same is true for a repelling set and a filtrating set.
\end{lemm}
\begin{proof} By definition, $f(X \setminus O)=g(X \setminus O)$ and therefore $f(O)=g(O)$.  
Thus, $g(O) \subset \mathring{O}$ if and only if 
$f(O) \subset \mathring{O}$, in other words,
$O$ is an attracting set for $g$ if and only if it is so for $f$. 

For the assertion about the repelling set, note that 
in general, $\supp (g, f) \subset O$ does not imply 
$\supp (g^{-1}, f^{-1}) \subset O$. Instead, we know  
$\supp (g^{-1}, f^{-1}) \subset f(O) = g(O)$.
Thus we know that $f(O) = g(O)$ is an attracting set for 
$f^{-1}$ if and only if it is so for $g^{-1}$. Then, we can 
deduce the conclusion by Remark~\ref{r.atre}, (2).

Assume now that $O=A\cap R$ where $A$ and $R$ are an attracting and a repelling set for $f$. 
As $g$ coincides with $f$ outside $A$ and outside $R$, the 
previous arguments imply that $A$ is still 
both attracting and repelling for $g$. 
Thus $O$ is a filtrating set for $g$. 
\end{proof}

One of the main properties of a filtrating set $U$ is that, 
if $x\in U$ and $f(x)\notin U$ then for any $n>0$ one has $f^n(x)\notin U$ and the similar result holds for $f^{-1}$. 
In other words, for any $x\in X$ the set of $n\in\ZZ$ such that $f^n(x)\in U$ is an interval.  

This property remains true for pseudo orbits.
Recall that a (finite or infinite) sequence $(x_n) \subset X$ is 
called an $\varepsilon$-pseudo orbit if for every $i$ we have
$d(f(x_i), x_{i+1}) < \varepsilon$ holds, whenever it makes sense. 
For a filtrating set $U$, 
there is $\varepsilon>0$ such that if $x_i, i\in\ZZ$ is an $\varepsilon$-pseudo orbit satisfying $x_0\in U$ and $x_1\notin U$, then $x_i\notin U$ for $i>0$. 
A similar result holds for $f^{-1}$.   

Recall that a point $x \in X$ is called
\emph{chain recurrent} if for any $\varepsilon >0$ we can find an $\varepsilon$-pseudo orbit
$(x_i)_{i=0,\ldots,n}$ ($n>0$) such that $x_0 = x_n =x$.
For a chain recurrent point $x$, the chain recurrence class 
$[x]$ is the set of points $y \in X$ for which we can find 
$\varepsilon$-pseudo orbits starting from $x$ and ending at $y$, and vice versa for any $\varepsilon >0$. One can 
prove that $[x]$ is a compact $f$-invariant set of $X$.

An important consequence is that $U$ is saturated for the chain recurrence classes, 
that is, chain recurrence class is either disjoint or contained in $U$, if $U$ is a filtrating set.  

%
%
%

\subsection{Partially hyperbolic filtrating Markov partitions and the property $(\ell)$}\label{s.preliminar}

In this subsection, $f$ denotes a $C^1$-diffeomorphism
on a three dimensional closed manifold $M$.
In \cite{BS2}, we introduced the notion of 
\emph{partially hyperbolic filtrating Markov partitions},
which we just call filtrating Markov partitions in the sequel 
for simplicity.
It is a set of information about a filtrating set 
having some partial hyperbolicity adapted to its shape. 
It enables us to conclude several properties about the 
isolation of chain recurrence classes.

Let us review the concept of filtrating Markov partitions.
Since in this paper we do not need the precise definition 
of them, we only describe some rough ideas 
of what they are. Those who are 
interested in precise definition of them, 
see \cite[Definition~2.5]{BS2}.

By a \emph{cylinder} we mean a compact 
subset $C$ of a three dimensional manifold which is 
$C^1$-diffeomorphic to a cylinder $\mathbb{D}^2 \times I$,
where $\mathbb{D}^2 \subset \mathbb{R}^2$ is a round
disk and $I \subset \mathbb{R}$ is an interval.
For a cylinder, by the lid boundary we mean the subset of 
the boundary which corresponds to 
$\mathbb{D}^2 \times (\partial I)$
and the side boundary is to 
$(\partial \mathbb{D}^2) \times I$.

A (partially hyperbolic) filtrating Markov partition (of saddle type)
is a filtrating 
set $\cU = A \cap R$ of a diffeomorphism $f$ 
having additional properties such as:
\begin{itemize}
\item $\cU$ is a disjoint union of finitely many cylinders:
$\cU = \cup C_i$ where $C_i$ is a cylinder.
Each $C_i$ is referred as a rectangle.
\item For each rectangle, the lid boundary is contained in
$\partial R$ and the side boundary is in $\partial A$.
\item For each rectangle $C_i$, the set 
$f(\cU) \cap C_i$ consists of finitely many cylinders
which properly crosses $C_i$ in the axial direction. 
\item There is a partially hyperbolic structure of the 
form $E^u \oplus E^{cs}$ in the neighborhood of $\cU$
satisfying the following:
\begin{itemize}
\item $E^u$ is one dimensional and uniformly expanding.
It is almost parallel to the axial direction of rectangles.
\item $E^{cs}$ is two dimensional and 
$E^u \oplus E^{cs}$ forms a dominated splitting. 
It is almost parallel to the lid boundary of rectangles.
\end{itemize} 
\end{itemize}

For filtrating Markov partitions, we can deduce the following:
\begin{itemize}
\item Being a Markov partition is a $C^1$-robust property 
(see \cite[Lemma~2.9, (1)]{BS2}).
\item For each rectangle $C_i$, the set 
$f^{-1}(\cU) \cap C_i$ consists of finitely many cylinders
which properly crosses $C_i$ in the horizontal direction 
(see \cite[Proposition~2.10]{BS2}).
\item For each $m, n \geq 0$, the set 
$f^{-m}(\cU) \cap f^{n}(\cU)$ is 
a filtrating Markov partition, too. 
We call it the $(m, n)$-refinement 
of $\cU$ and denote it by $\cU_{(m, n)}$ (see \cite[Corollary~2.13]{BS2}).
\item $ \cap_{m \geq 0} f^m(\cU)$ is a disjoint 
union of $C^1$-lamination of curves which coincides with the 
unstable manifold of the locally maximal invariant set of 
$\cU$ (see \cite[Lemma~2.15]{BS2}).
\item $ \cap_{m \geq 0} f^{-m}(\cU)$ is a 
continuous family of $C^1$-discs such that each disc cuts
the rectangle it belongs to (see \cite[Lemma~2.15]{BS2}).
\end{itemize}

A useful property of a filtrating Markov partition is that 
we can determine the isolation of chain recurrence classes 
by local information. Let us explain this. 
For a Markov partition $\cU$, the 
backward invariant set $\cap_{k\geq 0} f^{-k}(\cU)$
is a continuous family 
of two dimensional center stable manifolds.
Suppose that we have a hyperbolic periodic point $p$ in $\cU$ of 
stable index two. 
We denote the connected component of
$\cap_{k\geq 0} f^{-k}(\cU)$ containing $p$
by $W^{cs}_{\mathrm{loc}}(p)$. 
As explained above one can prove 
that $W^{cs}_{\mathrm{loc}}(p)$ is 
$C^1$-diffeomorphic to a two dimensional disc which cuts the 
rectangle containing $p$ horizontally. 
We say that $p$ has a \emph{large stable manifold} 
if $W^s(p)$ contains $W^{cs}_{\mathrm{loc}}(p)$.
For periodic points having a large stable manifold, we 
can determine the information of the size of the chain 
recurrence class of $p$ just by looking its fundamental 
domain, since the image of the fundamental domain under 
$f^{-1}$ goes out from the filtrating set and never comes 
back (see \cite[Proposition~2.23]{BS2}). 

This information enabled us to conclude the abundance of 
the isolation of saddles 
(see \cite[Corollary~1.2]{BS2}), together with the notion of 
\emph{$\varepsilon$-flexible points}.
We will not review the precise definition of it (see \cite[Section~3.1]{BS2} for the definition). It is a 
hyperbolic periodic point of stable index two whose derivative 
cocycle in the center stable direction admits very peculiar 
deformation of size $\varepsilon$ in the $C^1$-distance:
By a continuous deformation of the cocycle of size 
$\varepsilon$ it can be deformed into a stable index one 
periodic point, and it also can be deformed into a stable 
index two periodic point having non-real eigenvalues in the 
center stable direction. 

By combining these two kinds of deformations, we proved 
in \cite{BS1} that, for an $\varepsilon$-flexible point, 
by performing a perturbation whose $C^1$-size is 
less than $\varepsilon$, we can 
deform the point to be a stable index one 
periodic point having a neutral eigenvalue such 
that in the fundamental domain of the center stable manifold 
the position of the strong stable manifold is given by any
prescribed $C^1$-curve subject to a topological limitation
(see \cite[Theorem~1.1]{BS1}). 

Thus, if we have a filtrating Markov partition which contains 
$\varepsilon$-flexible points with large stable manifolds, 
then by performing 
a $C^1\mbox{-}\varepsilon$-small perturbation
we can make their chain recurrence classes to be 
the periodic orbit themselves.
In \cite{BS1, BS2, BS3}, we have shown that the existence 
such periodic points for an arbitrarily small $\varepsilon>0$
in filtrating Markov partitions which satisfy a condition called 
property $(\ell)$
(see \cite[Proposition~5.1]{BS1}, \cite[Lemma~3.8]{BS2},
\cite[Section~4]{BS3}).
As a result, we proved that such a diffeomorphism is 
\emph{wild}, that is, $C^1$-generically it has infinitely 
many chain recurrence classes. 

\bigskip

Let us review the definition of the property $(\ell)$ for a chain recurrence 
class, which guarantees the existence of $\varepsilon$-flexible
points with large stable manifolds for small $\varepsilon>0$. 
In this article 
by $\mathrm{Diff}^1(M)$ we denote the set of $C^1$-diffeomorphisms of $M$ with the $C^1$-topology.

\begin{defi}[Definition~2, \cite{BS3}]
Let $f \in \mathrm{Diff}^1(M)$ and $p$ be a 
hyperbolic periodic point of stable index two. 
We say that $[p]$, the chain recurrence class of $p$, 
satisfies property $(\ell)$ if the following holds:
\begin{itemize}
\item $[p]$ is contained in a filtrating Markov
partition $\cU$ such that $p$ has a large stable manifold in 
$\cU$.
\item $p$ is homoclinically related to a hyperbolic periodic 
point $p_1$ of stable index two such that $Df^{\pi_1}_{p_1}$
has non-real eigenvalues (where $\pi_1$ is the period of $p_1$).
\item There are two hyperbolic sets $K$ and $L$ such that 
\begin{itemize}
\item $K$ has stable index two and contains $p$.
\item $L$ has stable index one. 
\item $K$ and $L$ form a robust heterodimensional cycle
(see Proposition~5.1 of \cite{BS1}). 
\end{itemize}
In the following, we refer this property as \emph{$p$ is 
in a robust heterodimensional cycle}.
\end{itemize}
\end{defi}

Note that, in the above definition, the fact that $\cU$ is filtrating 
forces $p_1$ and the heteroclinic points between $p$ and $p_1$
to be contained in the interior of $\cU$. A similar property 
holds for the robust heterodimensional cycle, too.

In this article, 
we also consider a local version of property $(\ell)$,
see also \cite[Definition 4.1]{BS3}.
Let $p$ be a hyperbolic periodic point of stable index two 
contained in a filtrating Markov partition $\mathcal{U}$. 
By a \emph{sub Markov partition of $\cU$} 
or a \emph{non-filtrating Markov partition}
we mean the union of the collection of rectangles of $\cU$.
Note that in general a sub Markov partition is not a filtrating set.
Let $\cW$ be a sub Markov partition 
of $\mathcal{U}$ which contains 
the orbit of $p$. By a \emph{relative homoclinic class} 
$H(p, \mathcal{W})$ we mean the closure of the 
transverse homoclinic intersections 
of $W^u(p)$ and $W^s(p)$ whose orbit 
is in $\mathcal{W}$. Notice that, due to the 
filtrating property of $\cU$, $H(p, \mathcal{W})$ is 
a compact $f$-invariant set which is contained in the 
interior of $\mathcal{W}$. 

\begin{defi}
Let $p$, $f$, $\mathcal{W}$ be as above. 
We say that $[p]$ 
satisfies property $(\ell_{\mathcal{W}})$
if the following holds:
\begin{itemize}
\item $p$ has a large stable manifold in $\cU$.
\item There is a hyperbolic saddle $p_1$ of stable index two
such that $p$ and $p_1$ are homoclinically related in $\cW$.
\item $p$ has a robust heterodimensional cycle in $\cW$,
that is, there are hyperbolic sets $K$, $L$ of different indices
contained in $\cW$
such that $K$ contains $p$ and $K$, $L$ form a robust 
heterodimensional cycle in $\cW$.
\end{itemize}
\end{defi}

One can see that having property $(\ell_{\mathcal{W}})$
is a $C^1$-robust property. If $[p]$ has property 
$(\ell_{\mathcal{W}})$ then it implies that it has property 
$(\ell )$ for the Markov partition which contains 
$\mathcal{W}$. In this sense, $(\ell_{\mathcal{W}})$
is a condition which is stronger than $(\ell{})$. The advantage 
of the condition $(\ell_{\mathcal{W}})$ is that it guarantees 
the existence of the flexible point whose orbits are localized 
in $\mathcal{W}$, as we will see in the next subsection.


 \subsection{Tools from \cite{BS3}}
 \label{ss.tool} 
 
 We cite several results from \cite{BS3} which are used 
 to construct aperiodic classes. 
To give the results, we introduce a few definitions. 

Let $\cU$ be a filtrating Markov partition of $f$ and 
$g \in \mathrm{Diff}^1(M)$ is so $C^1$-close to $f$ that 
$\cU$ is still a filtrating Markov partition for $g$.
Then by $\cU_{(m, n;g)}$ we denote the $(m, n)$-refinement
of $\cU$ with respect to $g$. When we do not need to 
indicate with which map we take the refinement, we 
also use the notation $\cU_{(m, n)}$.

Let $\cU$ be a filtrating Markov partition and $S \subset \cU$.
By $\cU(S)$ we denote the sub Markov partition of $\cU$
which consists of the rectangles having non-empty intersection
with $S$. 

 A \emph{circuit} of points of $f$ is a  set of 
 finitely many hyperbolic periodic orbits $\{\cO(q_i)\}$
 of stable-index 2 and transverse homo/heteroclinic orbits 
 $\{\cO(Q_j)\}$ among them, see \cite[Section1.3]{BS3}.
Given a circuit, we can obtain a directed graph by setting 
vertices to be the periodic 
orbits and edges to be the homo/heteroclinic orbits. 
We say that a circuit is transitive
if the directed graph is transitive. In this article, we only treat transitive circuit and 
we always assume so without mentioning it. Note that $S$ is a uniformly 
hyperbolic set.

Suppose that we have circuits of points $K$ of $f$ and $L$ of $g$. 
We say that $K$ and $L$ are \emph{similar} if the following 
holds:
  \begin{itemize}
  \item There is a continuous bijection $h: K \to L $ which conjugates 
  $f$ and $g$: $ f\circ h = h \circ g$ holds on $K$. 
  \end{itemize} 
We say that they are \emph{$\delta$-close} if
the $C^0$-norm of 
$h$ can be chosen smaller than $\delta$,
more precisely, $d(x, h(x)) < \delta$ holds for all $x \in K$.  

 \subsubsection{Expulsion result}
 The first two results will be used to obtain a filtrating set  
 by a small  perturbation of the diffeomorphism 
 in such a way that it contains a prescribed circuit.
 The expulsion process consists of two steps. 
 For expelling a new filtrating set, we need to have that 
 the rectangles containing the circuit be \emph{affine},
 which guarantees the regularity of the shape of the 
 rectangles. Since the precise definition of it is not used 
 in this paper, we do not state it here. The reader interested 
 in can find it in Section~2 of \cite{BS3}.
 
 The second result claims that if a sub Markov partition
 is affine and all the periodic orbits involved are 
 $\varepsilon$-flexible points
 with large stable manifolds, then by adding $2\varepsilon$-perturbation and taking forward refinement we have the 
expulsion of a chain recurrence class, keeping the combinatorial 
information of the Markov partitions. 
Recall that for a filtrating Markov partition we have defined the notion of 
$\alpha$-robustness (see Section~2.7 of \cite{BS3}). We do not 
review the precise definition of it here. 
Roughly speaking, $\alpha$-robustness 
implies that the filtrating Markov partition persists against perturbations of $C^1$-size $\alpha$.
One property which we use about the robustness is that 
taking refinements does not decrease the robustness 
(see Section~2.7 of \cite{BS3}).
 
 Let us give the first result (see Theorem~1.4 of \cite{BS3}):
 \begin{theo}\label{t.aff}
 Let $f \in \mathrm{Diff}^1(M)$ having 
 an $\alpha$-robust filtrating Markov partition $\cU$
containing a circuit (of points) $S$. Assume that every periodic orbit of $S$ 
has a large stable manifold in $\cU$. Then for
any neighborhood $U$ of $S$ and any sufficiently small 
$C^1$-neighborhood $\rO$ of $f$
there is $g \in \rO$ such that the following holds:
\begin{itemize}
\item The support $\supp(g, f)$ is contained in $U$.
\item For $g$ we have a continuation of $S$ which we denote $S_g$.
Then, there exists $m_0, n_0 >0$ such that for every $m \geq m_0, n \geq n_0$,
the Markov partition $\cU_{(m, n; g)}(S_g)$
is affine.
\item The conjugate periodic points in $S$ and $S_g$ have 
the same orbits and the same derivatives along them.
\item $\cU_{(m, n;g)}$ is $\alpha$-robust, too.
\end{itemize}
\end{theo}

To give the second result, we prepare a definition.
\begin{defi}\label{d.macha}
Let $\cU$ be a (possibly non-filtrating) Markov partition
and $\cV$ be another (possibly non-filtrating) Markov partition such that 
\begin{itemize}
\item $\cV \subset \mathring{\cU}$,
\item each rectangle of ${\cU}$ contains one 
and only one rectangle of $\cV$. 
\end{itemize}
We say that $\cV$ is matching to $\cU$ if these 
properties hold.
\end{defi}
The information of being matching implies that the 
two Markov partitions have similar combinatorial information. 
This enables us to determine the properties of new Markov 
partitions we produce.

Then we have the following (see Theorem~1.5 of \cite{BS3}): 
\begin{theo}\label{t.isola}
Let $f \in \mathrm{Diff}^1(M)$ having an 
$\alpha$-robust filtrating 
Markov partition $\cU$ containing a circuit $S$ whose 
periodic orbits are all $\varepsilon$-flexible, have 
large stable manifolds in $\cU$ and $\cU(S)$ is affine. 
Suppose $\alpha > 2\varepsilon$. Then there exists $n_0>0$ such that for every $n\geq n_0$
there is a $C^1$-diffeomorphism $h_n$ 
which is $2\varepsilon$-$C^1$-close to $f$
and satisfying the following:
\begin{itemize}
\item $\supp(h_n, f)$ is contained in the interior of $\cU(S)$. 
\item $h_n$ has a transitive filtrating Markov partition $\cV$ containing 
a circuit $S_n$ which is similar to $S$. 
\item For $S_n$ we have the following:
\begin{itemize} 
\item We have $\cU(S_n)  = \cU(S)$ and 
we can require that the points which are conjugate 
belong to the same rectangle and 
the conjugate periodic orbits of 
$S_n$ has the same orbit as in $S$.
\item Every periodic orbit of $S_n$ has a large stable 
manifold in $\cV$ and it is $\varepsilon$-flexible.
\end{itemize}
\item For $\cV$ we have the following:
\begin{itemize}
\item Each rectangle of $\cV$ is a vertical sub rectangle 
of a rectangle of $\cU_{(0, n; h_n)}(S_n)$.
\item $\cV$ is matching to $\cU_{(0, n; h_n)}(S_n) = \cU_{(0, n; h_n)}(S) $.
\item $\cV$ is $(\alpha -2\varepsilon)$-robust, too.
\end{itemize}
\end{itemize} 
\end{theo}

%
%
 
 \subsubsection{Recovering condition ($\ell$) from flexibility}
 
 Our second tool will allow us to get a circuit satisfying convenient dynamical conditions (see Section~4.1.2 of \cite{BS3}).   
 
 \begin{theo}\label{t.relative} 
 Let $f \in \mathrm{Diff}^1(M)$, 
$p$ be a periodic point of $f$, 
$\cU= \cup U_i$ be a $\alpha$-robust filtrating Markov partition and $\cW$ be a sub Markov partition of $\cU$.
Assume that $p$ is $\varepsilon$-flexible,
has a large stable manifold,
the orbit of $p$ is contained in $\cW$, 
the relative homoclinic class $H(p, \cW)$ is non-trivial,
and $\alpha>4\varepsilon$ holds.
 
Then, there is a diffeomorphism $g$ which 
is $4\varepsilon$-$C^1$-close to $f$ 
such that $p$ satisfies the condition 
$(\ell_{\cW})$. Furthermore, we may assume the following:
\begin{itemize}
\item $\supp(g, f)$ is contained in the interior of $\cW$.
\item $\cU$ is still a filtrating Markov partition for $g$.
\item Suppose that $p$ is contained in a circuit $K$. 
 Then for any $\delta >0$ by choosing appropriate $g$ 
 we may assume 
 there is a circuit which is $\delta$-close to $K$ containing 
 $p_g$.
\end{itemize}  
 \end{theo}

\subsubsection{Abundance of flexible points}
The last result shows that given a periodic 
orbit satisfying condition $(\ell_{\cW})$ for some 
sub Markov partition $\cW$, 
up to an arbitrarily small local perturbation one can 
obtain a lot of flexible points with large periods, 
see Section~4.1.1 of \cite{BS3}.
\begin{theo}\label{t.flex}
Let $f \in \mathrm{Diff}^1(M)$ and $\mathcal{U}$ be a
filtrating Markov partition. Let $\cW$ be a sub 
Markov partition of $\cU$ and $p$ a hyperbolic periodic point whose 
orbit is contained in $\cW$. 
Assume that $p$ satisfies the condition $(\ell_{\cW})$. 
Then given $\varepsilon>0$, $N>0$ and 
a $C^1$-neighborhood $\mathrm{O}$ of $f$,
there is a diffeomorphism $g\in \mathrm{O}$ such that the following holds:
\begin{itemize}
\item $\supp (g, f)$ is contained in the interior of $\cW$.
\item There is a periodic point $x \in H(p_g,\cW;g)$ homoclinically 
related to $p_g$ in $\cW$ 
which is $\varepsilon$-flexible, having a large 
stable manifold in $\cU$, $\varepsilon$-dense in $H(p_g, \cW;g)$ and
whose period is larger than $N$. 
\end{itemize}
\end{theo}


\section{Topological principles for the genericity of certain aperiodic classes}\label{s.genurgo}

Throughout this section $(X, d)$ denotes a compact metric space, 
$\mathrm{G}\subset \mathrm{Homeo}(X)$ is a subgroup of homeomorphisms endowed with a topology finer than or equal to the 
$C^0$-topology such that $\rG$ is a Baire space with it. 
By $f$ we denote a homeomorphism in $\rG$. In the next section (Section~\ref{s.C1}), $X$ will be a compact manifold $M$ and  $\rG$ will be the space of $C^1$-diffeomorphism on $M$ endowed with the $C^1$-topology. 

In this section we discuss several properties of 
invariant sets defined as a limit of a nested sequence of filtrating sets.
Our arguments in this section are purely topological and 
we prefer to give them in a topological setting for avoiding the implicit use of properties resulting from the differentiability assumptions.

\subsection{Stability for local perturbation}

Suppose that $f$ has a compact set $\cK \subset X$ on which
$f$ satisfies a certain property.
We are interested in if the property holds for $g$ which 
differs only inside $\cK$. In order to discuss it we introduce
a few definitions. 

\begin{defi}
Let $(X, d)$, $f$, $g$ and $\cK$ be as above. We say that 
the \emph{support of $g$ (with respect to $f$) is strictly contained in $\cK$}
if the following holds:
\begin{itemize}
\item $\supp (g, f) \subset \cK$.
\item for each connected component $K$ of $\cK$,
we have $K \setminus \supp (g, f) \neq \emptyset$.
\end{itemize}
A property is called \emph{stable under local perturbation}
if whenever it holds for $f$ and $\cK$ then it holds for 
$g$ and $\cK$ where the support of $g$ is strictly contained 
in $\cK$.
\end{defi}
By the discussion of Section~\ref{ss.filt}, we know that
being an attracting set, a repelling set or a filtrating set 
is a property which are stable under local perturbations.
Let us see an important property of 
homeomorphisms which 
is preserved by local perturbations. 
\begin{lemm}\label{l.1step}  Let
$\cK\subset X$ be a compact set. 
Let $f, g \in \rG$ such that $\supp (g, f)$ is 
strictly contained in $\cK$.
Then for every connected component $L$ of $\cK$ one has 
$f(L)=g(L).$
\end{lemm}

\begin{proof}
By assumption we know $f(\cK) = g(\cK)$. 
Thus, the sets of connected components of $f(\cK)$ and $g(\cK)$ are the same. Consider a connected component $L$ of 
$\cK$. Then $f(L) \cap f(\cK) \neq \emptyset$. 
Since $f$ is a homeomorphism, $f(L)$ is one of the 
connected components of $f(\cK)$ and the same holds 
for $g(L)$. Consider a point 
$x\in L$ for which we have $f(x) = g(x)$. Note 
that the existence of such $x$ is guaranteed by the assumption.
Then $f(L)$ and $g(L)$ are connected components of 
$f(\cK) = g(\cK)$ which contains the common point 
$f(x) = g(x)$. Thus, they coincide, that is, $f(L) = g(L)$.
\end{proof}

\begin{rema}
A natural way of stating ``the support is contained in $\cK$''
would be ``$f \equiv g$ outside $\cK$.'' However, 
this condition is not enough to obtain the conclusion 
in Lemma~\ref{l.1step} in general. For instance, 
consider a compact metric space $Y$ which consists of
two homeomorphic connected components. Let $\cK = Y$
and consider $f = \mathrm{id}_Y$ and some map 
which exchanges two connected components which we denote by $g$. Then for an
obvious reason $f \equiv g$ outside $\cK$
but in the level of the connected 
components they are different. 
\end{rema}

Meanwhile, under a mild extra condition we have 
the same conclusion for these two assumptions. Let us see it.

\begin{lemm}\label{l.supo}
Suppose that $X$ is connected and 
$\cK$ is compact, has finitely many connected components and
$\cK$, $X \setminus \cK$ are non-empty.
If $f\equiv g$ on $X\setminus \cK$,  
then the support of $g$ is strictly contained in $\cK$
(thus the conclusion of Lemma~\ref{l.1step} holds). 
\end{lemm}

\begin{proof} As $f$ and $g$ coincides on $X\setminus \cK$ we have $g(\cK)=f(\cK)$, and even $f\equiv g$ on $\partial_X \cK$, where $\partial_X \cK$ denotes the 
boundary of $\cK$ in $X$. Thus the conclusion is the 
direct consequence of the following:
\begin{clai} For any connected component $L$ of $\cK$
we have $L\cap \partial_X\cK\neq\emptyset.$
\end{clai}
Let us prove the Claim. Take a small open neighborhood $O$ of the 
connected component $L$ in $X$. 
Since $K$ is compact and has finitely many connected components,
we may assume $\cK \cap O = L$.
We have $O\setminus L \neq \emptyset$, otherwise it contradicts 
the connectedness of $X$. 
Thus $O$ contains a point $x_O \in X\setminus \cK$.  
By shrinking $O$, we obtain an accumulation of 
such points in $K\cap\partial_X\cK$. 
\end{proof}

\begin{rema}\label{r.equi}
In Section~\ref{s.C1}, we consider the case where $X$ is a 
compact connected manifold and $\cK$ has finitely many connected 
components and is not equal to neither $\emptyset$ nor $X$. 
Thus Lemma~\ref{l.supo} is always applicable.
\end{rema}

In this article, we discuss two kinds of robustness
of properties: 
The robustness with respect to the $\rG$-topology,
that is, the topology of the space of dynamical systems we are 
interested in, and the robustness with respect to 
the location of the support of the perturbation. 
For the proof of our main theorem, the second 
kind of robustness is not necessary, but we will discuss 
it. Let us see why. 

The first reason is that it enables us to make some 
of the proofs simpler. In the construction of aperiodic 
classes, we produce filtrating sets 
by small perturbations using the techniques 
in Section~\ref{ss.tool}. 
In the course of the construction, 
we use them successively and we need to consider 
their interference.   
Then using only
$\rG$-robustness brings some non-trivial 
complication of the proof. 
One advantage of stability under the local perturbation 
is it makes some arguments about the interference of the 
perturbations much simpler, see Remark~\ref{r.local-why}. 

%
%

Another reason is that the stability by the local perturbations 
enables us to construct a concrete example of dynamical systems
which exhibits the chain recurrence class we announced. 
We prove a perturbation result which produces a filtrating 
set by a local perturbation. Doing it successively, we obtain 
a sequence of maps $f_n$ which converges to a 
map $f$ having the desired chain recurrence set. 
This construction is more concrete than the genericity 
approach, because for this $f$ we certainly know the behavior
outside the support of the perturbation. See Remark~\ref{r.concrete}.

\subsection{A criterion for chain recurrence}

In this subsection, we discuss a sufficient condition which 
guarantees the chain transitivity for an invariant set obtained 
as the limit of nested filtrating sets.

%
%

For $x\in \cK$, we denote the connected component of $\cK$
containing $x$ by $\cK(x)$.

\begin{defi}
Let $\cU\subset X$ be a compact set and $f\in \rG$.
\begin{itemize}
\item $\cU$ is said to be \emph{regular} 
if for any pair of connected components $U_1, U_2$
of $\cU$,
$f (U_1) \cap U_2\neq \emptyset$
implies
$f (\mathring{U}_1) \cap \mathring{U}_2 \neq \emptyset$.
\item We say that a (finite or infinite) sequence 
$(x_i) \subset \cU$ is a \emph{$\cU$-chain of 
points} or \emph{$\cU$-pseudo orbit} if for any $i$ we have
$f (\cU (x_i)) \cap \cU(x_{i+1})\neq \emptyset$.
\item Also, we say that a finite or infinite sequence 
$(U_i)$, where $U_i$ is a connected component of $\cU$,
is a \emph{$\cU$-chain of connected components} if for any $i$ we have
$f(U_i) \cap U_{i+1} \neq \emptyset$
\item We say that $\cU$ is \emph{$\cU$-chain transitive} if given any two connected components $V_0$, $V_1$ of $\cU$, there 
is a finite $\cU$-chain of connected components $(U_i)$ 
starting from $V_0$ and ending at $V_1$. 
\item The \emph{minimum $\cU$-period} of $\cU$ is the smallest length of a periodic $\cU$-chain of connected components, see Remark~\ref{r.chai}. 
\end{itemize}
\end{defi}
For $\cK \subset X$, 
by $\mathrm{cdiam}(\cK)$ we denote 
the supremum of the diameter of the connected components
of $\cK$.
\begin{rema}\label{r.chai}
For $\cU$-chains, we have the following.
\begin{enumerate}
\item A $\cU$-chain of points $(x_i)$ defines a unique 
$\cU$-chain of connected components $(U_i)$ such that 
$x_i \in U_i$ (set $U_i = \cU(x_i)$). 
\item Given a $\cU$-chain of connected components
$(V_i)$
we can find a $\cU$-chain of points $(y_i)$ such that 
$y_i \in V_i$ by choosing a point from each $V_i$. 
Due to these correspondences, the definition of the 
chain transitivity and minimum $\cU$-period can be given 
in terms of $\cU$-chain of points. 
\item If $\cdiam (\cU), \cdiam (f(\cU)) < \varepsilon$ then 
$\cU$-chain of points $(x_i)$ defines 
a $2\varepsilon$-pseudo orbit.
\end{enumerate}
\end{rema}

Let us discuss basic properties of $\cU$-chain transitivity.

\begin{lemm}
Suppose that $f \in \rG$ has a compact set $\cU$ which is
regular, has finitely many connected components, and 
is $\cU$-chain transitive for $f$.
Then, for $g$ which is $C^0$-close to $f$, 
$\cU$ is $\cU$-chain transitive for $g$ as well.
\end{lemm}

\begin{proof}
The regularity condition guarantees that 
once we have  $f(U_1) \cap U_2 \neq \emptyset$ then
the $g(\mathring{U}_1) \cap \mathring{U}_2 \neq \emptyset$ holds for $g$ which is sufficiently $C^0$-close 
to $f$. Thus, if the number of the connected components is 
finite, all the condition as such holds for nearby homeomorphisms,
which implies the $\cU$-chain transitivity.
\end{proof}


\begin{lemm}
Suppose that $f$ has a compact set $\cU$ which is 
$\cU$-chain transitive.
If the support $\supp(g,f)$ is strictly contained in $\cU$,
then $\cU$ is $\cU$-chain recurrent for $g$, too.
\end{lemm}

\begin{proof}
The hypothesis ensures that the set of $\cU$-chains of connected components 
are the same for $f$ and $g$, 
see Lemma~\ref{l.1step}.
Thus, the $\cU$-chain tranisitivity of $f$ and that of $g$ 
are equivalent.
\end{proof}

\begin{prop}\label{p.chainrecurrence}

Consider a sequence $\{(f_n,\cU_n)\}_{n \geq 1}$ where $f_n\in \rG$
and $\cU_n$ is a filtrating set for $f_n$ such that:
\begin{itemize}
 \item $\cU_{n+1}\subset \mathring{\cU}_n$,
 \item $\supp (f_{n+1}, f_n)$ is strictly contained in $\cU_n$,
 \item $\cdiam (\cU_n), \cdiam (f_n(\cU_n)) \to 0$ as $n\to 0$,
 \item $\cU_n$ is $\cU_n$-chain transitive for $f_n$, 
 \item The minimum period of $\cU_n$ tends to $\infty$ as $n\to\infty$ 
\end{itemize}
Then the sequence $f_n$ converges to 
a homeomorphism $f$ in the $C^0$-topology 
for which $\La=\bigcap_{n\geq 1} \cU_n$ is an 
aperiodic chain recurrence class (especially, $\Lambda$ is $f$-invariant). 
\end{prop}
\begin{proof} 
For each $n$, consider the sequence $(f_{n+k})_{k \geq 0}$. 
Each $f_{n+k}$ is
equal to $f_n$ on the open set $X\setminus \cU_n$. 
This implies that the limit $f$ is well defined on $X\setminus\La$. 
Also, for each connected component $U$ of
$\cU_n$, the image $f(U)$ is well defined and coincides with 
$f_n(U)$. Note that this implies $\cU_n$ is $\cU_n$-chain transitive 
for $f$, too. 

The hypothesis on the diameters ensures that
for $x \in \Lambda$, 
$\bigcap_n(f(\cU_n(x))$ is a singleton. 
We denote it by $y$.  
We extend $f$ over $\La$ by mapping $x \in \Lambda$ to $y$. 
Then one can check 
that the map $f\colon X\to X$ is a homeomorphism. 
The continuity and the injectivety of $f$ on $X$ are easy to see. 
The surjectivity of $f$ on
$X \setminus \Lambda$ is also easy. 
Let us show $\Lambda \subset f(\Lambda)$. 
Take $x\in \Lambda$. By the $\cU_n$-chain 
transitivity of $f$, for every $n$  
there is a $\cU_n$-chain of connected components for $f$
starting from and ending at $\cU_n(x)$. For each $n$, 
we consider the second last connected component of the chain
and denote it by $V_n$. 
Let us choose a sequence $(z_n) \subset X$ such that $z_n \in V_n$
and take its accumulation point $z$. Then, the compactness of $\cU_n$
and the nested property of $(\cU_n)$ imply that 
$z \in \cap_{n \geq 1} \cU_n = \Lambda$. Also, the assumption 
on the diameters ensures that $f(z) =x$. Accordingly, we have 
$\Lambda \subset f(\Lambda)$.

Let us show that $\La$ is a chain recurrence class for $f$. 
First, the condition that 
the support of $f_{n+1}$ is strictly contained in 
$\cU_{n}$ tells us  
that $\cU_n$ is $\cU_n$-chain transitive for $f$. 
It implies that every point in $\Lambda$ is a chain 
recurrent point (in the usual sense) and every pair of points in $\Lambda$
is chain equivalent, see Remark~\ref{r.chai}. As a result, we see that
$\Lambda$ is contained in a chain recurrence class.
Furthermore, since $(\cU_n)$ is a decreasing sequence of filtrating sets, 
the points outside $\Lambda$ cannot be chain equivalent to points 
in $\Lambda$. Thus, $\Lambda$ itself forms an entire chain recurrence class.

If $\La$ contains a periodic point $p$ of $f$, then the 
$\cU_n$-period of the $\cU_n$ for $f$  is bounded by
the period of $p$ 
since its orbit for $f$ is a $\cU_n$-chain of points for every $f_n$. 
This confirms the assertion about the aperiodicity.    
\end{proof}

\begin{defi}\label{d.chaine}  We say that 
a sequence $\{(f_n,\cU_n)\}_{n \geq 1}$ is a \emph{nested sequence
for an aperiodic class}
if it satisfies the hypothesis of Proposition~\ref{p.chainrecurrence}.
\end{defi}

\subsection{Minimality of aperiodic classes}

In this subsection we discuss the minimality of  
the chain recurrence classes obtained as a limit of 
nested filtrating sets. Recall that a compact invariant set 
$K \subset X$ is called \emph{minimal} if every point of $K$ has a
dense orbit in $K$, that is, for every 
$x \in K$ the set $\mathcal{O}(x) := \{ f^i(x)\}_{i \in \ZZ}$ 
is dense in $K$. 
\begin{defi}\label{d.f-num}
 Let $f \in \mathrm{G}$ and $\cU, \cV \subset X$ 
be compact filtrating sets  
of $f$ satisfying $\cV\subset \mathring\cU$.  We say that 
$\cV$ is \emph{$\cU$-minimal} if given any pair of connected components $U_0$ of $\cU$ and $V_0$ of $\cV$ there is $n>0$ such that 
$f^{n}(V_0)\subset \mathring{U}_0$ holds.
We denote by $n(U_0,V_0)>0$ the smallest positive integer $n$ 
satisfying $f^{n}(V_0)\subset \mathring{U}_0$.
\end{defi}

We give an important observation.

\begin{lemm}\label{l.first-time} If $\cV$ is $\cU$-minimal, then for any connected components $U_0$ of $\cU$ and $V_0$ of $\cV$ 
we have the following equality:
\[
n(U_0,V_0)=\min \{i>0 \mid f^i(V_0)\cap U_0\neq \emptyset\}.
\]

\end{lemm}
\begin{proof} For $i >0$ 
suppose $f^i(V_0)\cap U_0\neq \emptyset$ 
but $f^i(V_0)\not\subset \mathring U_0$. 
Then there is $x\in V_0$ such that $f^i(x)\notin \mathring \cU$. As $\cU$ is a filtrating set, if a point $y\in \cU$ satisfies 
$f(y)\notin\mathring \cU$ then $f^k(y)\notin \cU$ for all $k\geq 0$. Applying this property to the orbit of $x$  we have $f^j(x)\notin \cU$ for any $j>i$, which implies $n(U_0,V_0)\leq i$. 
\end{proof}
%
%

The following lemma tells us that the $\cU$-minimality 
is a $C^0$-robust property under a mild assumption.

\begin{lemm}\label{l.minimal-robust}Let $f \in \mathrm{G}$ and
$\cU$, $\cV$ be filtrating sets of $f$. 
We assume that both $\cU$ and $\cV$ have 
finitely many connected components.  
If $\cV$ is $\cU$-minimal, then there is 
a neighborhood $\mathrm{O}$ of $f$ in 
$\mathrm{G}$ such that 
for any $g\in\mathrm{O}$ the compact
set $\cV$ is $\cU$-minimal. 
\end{lemm}

\begin{proof}
For any pair of connected components of $U_0$ of $\cU$ 
and $V_0$ of $\cV$, we have $f^{n(U_0, V_0)}(V_0) \subset \mathring{U_0}$.
This condition is a $C^0$-robust, hence $\rG$-robust condition. 
Since the number of connected components of $\cU$ and $\cV$ are finite, 
having this condition for every pair of $U_0$ and $V_0$ is also $\rG$-robust.
\end{proof}

\begin{lemm}\label{l.minimale}
Let $f \in \mathrm{G}$ and  $\cU$, $\cV$ 
be compact filtrating sets of $f$ such that 
$\cV$ is $\cU$-minimal. 
Let $g$ be a homeomorphism 
whose support is strictly contained in $\cV$. Then $\cU$, $\cV$ are compact filtrating sets of $g$ and  $\cV$ is $\cU$-minimal for $g$, too.
\end{lemm}
\begin{proof} The fact that $\cU,\cV$ are filtrating sets for $g$ is given by lemma~\ref{l.filtrating}. We just need to prove 
the $\cU$-minimality of $\cV$ for $g$. 

Recall that we have defined the number $n(U_0, V_0)$
for any pair of connected components of $U_0$ of $\cU$
and $V_0$ of $\cV$, see Definition~\ref{d.f-num}.
By definition, we have 
\[
f^{n(U_0, V_0)}(V_0) \subset \mathring{U_0}
\]
for every pair of $U_0$ and $V_0$. We prove that the 
same relation holds for $g$, that is, we have
\[
g^{n(U_0, V_0)}(V_0) \subset \mathring{U_0},
\]
which implies the $\cU$-minimality of $\cV$. 

Assume that this does not hold for some $U_0$ and $V_0$.
We take the minimum number of $n(U_0, V_0)$ among such 
$U_0$ and $V_0$ and denote it by $N$. 
Namely, consider
\begin{equation}\label{e.1}
N=\min\{i>0 \mid \exists U_0, \exists V_0  \mbox{ such that } i=n(U_0,V_0) \mbox{ and }  g^i(V_0)\not\subset \mathring{U}_0 \}.
\end{equation}
Note that we have $N>1$ because $g(V_0)=f(V_0)$. 

Consider the set $g(V_0)=f(V_0)$.  
By assumption, there is  $x\in V_0$ such that $g^{N}(x)\notin \mathring U_0$. 
We take $y\in V_0$ such that $g(x)=f(y)$ holds.
Then we have the following:
\begin{clai}
$g(x) \in \mathring{\cV}$.
\end{clai} 
\begin{proof}
If $g(x) \notin \mathring \cV$ then we have 
$g^2(x)= g(g(x)) = f(g(x)) = f(f(y)) = f^2(y)$, where 
in the second equality we used the fact that $f \equiv g$ outside $\cV$.
Furthermore, since $\cV$ is filtrating for both $f$ and $g$,
we know that the point $g^2(x)=f^2(y)$ does not belong to
$\mathring{\cV}$, either. It enables us to repeat the 
argument for the point $g^2(x)=f^2(y)$. 
Thus, arguing inductively, 
we have $g^k(x)=f^k(y)$ for any $k>0$. 
Now this implies $g^{N}(x) = f^{N}(y)\in \mathring U_0$ which contradicts the choice of $x$. 
\end{proof}
Thus $g(x)=f(y)$ belongs to $\mathring \cV$ and therefore to 
some connected component of $\cV$, which we denote by $V_1$. 
Then consider the integer $n(U_0, V_1)$.
Since $f^{N}(y) = f^{N-1}(f(y)) \in \mathring{U_0}$, 
together with Lemma~\ref{l.first-time} and 
the definition of $n(U_0, V_1)$, we know that 
$n(U_0, V_1) = N-1$. To be precise,
 $f^{N-1}(f(y)) \in \mathring{U_0}$
and Lemma~\ref{l.first-time} shows $n(U_0, V_1) \leq N-1$,
and it cannot be smaller because if so, it contradicts
$n(U_0, V_0)=N$.

 On the other hand, note that we have 
$g^{N-1}(f(y)) = g^{N-1}(g(x)) = g^{N}(x) \not\in \mathring{U_0}$. Thus $g^{N-1}(V_1) \not\subset \mathring{U_0}$. 
However, this contradicts the minimality of $N$,
see (\ref{e.1}).
Thus the proof is completed.
\end{proof}

We prepare a result which enables us to 
construct a concrete example of an aperiodic class.
\begin{prop}[A criterion for minimality]
\label{p.minimale} Let 
$\{(f_n, \cU_n)\}$ 
be a nested sequence for an aperiodic class such that
 $\cU_{n+1}$ 
is $\cU_n$-minimal for $f_{n+1}$ for every $n$.
Then the sequence $(f_n)$ converges in the $C^0$-topology to a 
homeomorphism $f\colon X\to X$ for which $\La=\bigcap \cU_n$ is an aperiodic 
chain recurrence class which is minimal and has $\cU_n$ 
as a basis of (filtrating) neighborhoods.  
\end{prop}

\begin{proof}
The fact that $\Lambda$ is an aperiodic class 
is confirmed in Proposition~\ref{p.chainrecurrence}. Thus 
we only need to check the minimality.

As $\cU_{n+1}$ is $\cU_n$-minimal for $f$,  given $n$ and $ x,y\in\La$  there is $k>0$ such that $f^k(\cU_{n+1}(x))\subset \mathring{\cU}_n(y)$.  As the diameter of $\cU_n(y)$ tends 
to zero, this shows that the orbit of $x$ passes arbitrarily close to $y$.  As $x, y$ are any points in $\La$ this shows that the minimality of $\La$.  
\end{proof}

\begin{rema}\label{r.concrete}
For the proof of Theorem~\ref{t.minimal-non}, 
we only need the case $f_n = f$. 
For instance, in the proof of Proposition~\ref{p.principle-expansive},
we use Proposition~\ref{p.minimale} letting $f_n = f$.
Meanwhile, this version has an advantage. It enables 
us to construct examples of homeomorphisms 
satisfying the desired condition in a concrete way. 
See Proposition~\ref{p.exa}.
\end{rema}

\subsection{A criterion for expansiveness}
\label{ss.exlo}

In this subsection, we describe a sufficient condition 
for the expansiveness of the limit invariant set of a 
nested sequence of filtrating set.
Recall that a compact set $K \subset X$ is \emph{expansive}
if there is $\delta >0$ such that for any $x, y \in X$ $(x \neq y)$
there is $n = n(x, y)$ satisfying $d(f^n(x), f^n(y)) > \delta$.

\begin{defi}
\label{l.def-exp}
Let $f \in \rG$ having a compact filtrating set $\cU$
and $\cK \subset \mathring{\cU}$ be a compact set. 
We say that $\cK$ is \emph{$\cU$-expansive for chains} 
(with respect to $f$)
if given any two bi-infinite $\cK$-chains 
$(x_i)$ and $(y_i)$  $(i\in \ZZ)$ of points
we have the following dichotomy:
\begin{itemize}
\item Either for any $i\in\ZZ$ the connected components $\cK(x_i)$ and $\cK(y_i)$ are equal, 
\item or, there is $i\in \ZZ$ such that $\cU(x_i)\neq\cU(y_i)$. \end{itemize}
 \end{defi}

\begin{lemm}\label{l.exlo}
Let $f$, $\cU$ and $\cK$ be as above.
If the support of $g$ is strictly contained in $\cK$, 
then $\cK$ is $\cU$-expansive for chains
with respect to $g$, too. 
\end{lemm}

\begin{proof}
%
Suppose that we have bi-infinite $\cK$-chains of points 
$(x_i)$ and $(y_i)$
with respect to $g$. We need to show 
if there is $i_0$ such that $\cK(x_{i_0}) \neq \cK(y_{i_0})$ 
then there is $i_1$ such that 
$\cU(x_{i_1}) \neq \cU(y_{i_1})$.

Lemma~\ref{l.1step} shows that 
$(\cK(x_{i}))$, $(\cK(y_{i}))$ are 
$\cK$-chains of connected components with respect to $f$, too.
Since for $f$ we have $\cU$-expansiveness of $\cK$, we know
there is $i_1$ such that 
\[
\cU(x_{i_1}) \neq \cU(y_{i_1}).
\]
For $i_1$ above we have $\cU(x_{i_1}) \neq \cU(y_{i_1})$,
which shows the $\cU$-expansiveness of $\cK$ with 
respect to $g$.
\end{proof}

In the same manner, we can prove the following:
\begin{lemm}\label{l.exper}
Let $f$, $\cU$ and $\cK$ be as above.
Assume that
the number of connected components of $\cK$ is finite.
Suppose that $\cK$ is $\cU$-expansive for chains 
with respect to $f$. 
Then, being 
$\cU$-expansive is a $C^0$-robust property.
More precisely, there exists $\varepsilon >0$ such that if 
$g \in \mathrm{Homeo}(X)$ satisfies
$d(f(x),g(x))<\varepsilon$ for every $x \in X$
then $\cK$ is $\cU$-expansive for $\cK$-chains 
with respect to $g$, too. 
\end{lemm}

\begin{proof}
In order to check the $\cU$-expansiveness for chains
 for $\cK$ 
with respect to $g$, the confirmation of the following 
information is enough: 
\begin{center}
If $g(K_i) \cap K_j \neq \emptyset$ then 
$f(K_i) \cap K_j \neq \emptyset$, equivalently, 
if $f(K_i) \cap K_j = \emptyset$ then 
$g(K_i) \cap K_j = \emptyset$.
\end{center}
Then, repeating the argument in the proof of 
Lemma~\ref{l.exlo} we obtain the condition. 
Note that, by the compactness of $K_i$ and 
$K_j$, the condition above for fixed $K_i$ and $K_j$
is a $C^0$-robust property.  Then the 
finitude of the number of the connected components 
enables us to obtain the conclusion.
\end{proof}

Now we can state a criterion to have
expansive aperiodic classes.
\begin{prop}\label{l.expansive} 
Let $\{(f_n,\cU_n)\}_{n \geq 1}$ be a nested sequence 
for an aperiodic class.  
Assume that 
\begin{itemize}
\item $\cU_{n+1}$ is $\cU_n$-expansive 
for chains with respect to $f_{n+1}$ for $n \geq 1$.
\item $\cU_1$ has finitely many connected components. 
\end{itemize}
Then the sequence $(f_n)$ converges to a homeomorphism $f$ such that $\La=\bigcap_{n \geq 1} \cU_n$ is an expansive chain recurrence class. 
\end{prop}

\begin{proof} 
We only need to prove the expansiveness of $\Lambda$.
Take two different points  $x, y \in \La$ and
consider the orbits $\{f^i(x)\}$ and $\{f^i(y)\}$, $i\in\ZZ$. 
There is $n$ such that the diameter of $\cU_n(x)$ is smaller than $d(x,y)$. Thus $\cU_n(x)\neq\cU_n(y)$, and the $f$-orbit of $x$ and $y$ define distinct $\cU_n$-chains of $f_n$. 
As $\cU_n$ is $\cU_{n-1}$-expansive 
with respect to $f_n$ and $f_n \equiv f$ outside $\cU_n$,
by Lemma~\ref{l.exlo} 
we know that the orbits of $x$ and $y$ are distinct 
$\cU_{n-1}$-chains of points for $f$. 

Arguing inductively, we know that they are distinct 
$\cU_1$-pseudo orbits of $f$. 
In other words, there is $i_0$ such that 
$f^{i_0}(x)$ and $f^{i_0}(y)$ belong to different components 
of $\cU_1$. As $\cU_1$ is assumed to have at most 
finitely many connected components,
there is $\delta>0$ such that two distinct components of 
$\cU_1$ are at distance larger than $\delta$. 
Thus we know $d(f^{i_0}(x),f^{i_0}(y))>\delta$ and 
this shows the expansiveness of $\Lambda$ for $f$.
\end{proof}

Let us state a result which enables us to obtain the generic existence 
of minimal, expansive aperiodic classes. 
We prepare some definitions.

\begin{defi}
We say that a $\rG$-robust property $\cP$ for 
$f$ on $\cU$ is type $\fP_{\MEx}$, where $f\in \rG$ and $\cU$ is 
a filtrating set of $f$, if the following holds: 
If $f$ satisfies $\cP$ then
for any $\delta>0$, $N >0$ and 
any $\rG$-neighborhood $\rO$ of $f\in \rG$ 
there is $g\in \rO$ satisfying the following:
\begin{itemize}
 \item The support of $g$ is strictly contained in $\cU$ with
 respect to $f$.
 \item There are disjoint, regular filtrating sets $\cU_0,\cU_1 \subset \mathring{\cU}$ for $g$ which are 
 $\cU$-minimal and $\cU$-expansive for chains.
 \item $\cU_i$ is $\cU_i$-chain transitive for $i=1, 2$.  
 \item The number of connected components of $\cU_1, \cU_2$
 are both finite,
 \item The minimum period of $\cU_i$  is larger than $N$ for $i=1,2$,
 \item $\mathrm{cdiam}(\cU_i), \mathrm{cdiam} (g(\cU_i))<\delta$ for $i=1,2$,
 \item The restrictions $g$ satisfies the property $\cP$ 
 on $\cU_i$ for $i=1,2$.
\end{itemize}
\end{defi}

If a dynamics have a robust property which is type $\fP_{\MEx}$, 
it implies the locally generic coexistence of uncountably many aperiodic
classes which are minimal and expansive. To explain this, 
we prepare some notation which will be used throughout this paper.
Consider the sets of finite or infinite sequences $\{1,2\}^k$, $k\in\NN$ and $\{1,2\}^\NN$.
For $\omega=(\omega_1,\dots, \omega_k)$ and any $j\leq k$,  or $\omega=(\omega_i)_{i\in\NN}$ and any $j\in\NN$ we write
$$[\omega]_j=(\omega_1,\dots,\omega_j),$$
that is, $[\omega]_j$ is 
the restriction of $\omega$ to its first $j$ terms. 

Now we can state a principle for the locally generic 
coexistence of minimal, expansive aperiodic classes.
\begin{prop}\label{p.principle-expansive}
If $f$ has a filtrating set $\cU$ such that 
$f$ satisfies $\mathrm{G}$-robust property
$\cP$ on a filtrating set $\cU$ which is type $\fP_{\MEx}$, then 
there is a neighborhood $\rO_1$ of $f$ in $\rG$ and a $\rG$-residual subset 
$\rR\subset \rO_1$ such that every $g\in\rR$ admits an uncountable family of aperiodic chain recurrence classes which are minimal and expansive.
\end{prop}
\begin{proof}As $\cP$ is $\mathrm{G}$-robust, 
there is a non-empty open neighborhood $\rO_1$ of $f$
such that every $g \in \rO_1$
satisfies $\cP$. 

Then the property $\cP$ being type $\fP_{\MEx}$ allows us to build by induction, a sequence of open subsets $\rO_n\subset \rO_1$ and a sequence $\delta_n>0$ tending to $0$ as $n\to\infty$ with the following property
(this kind of constructions appear in many papers, for instance  
in Section~3.3 of \cite{BD1} or Section~4.6 of \cite{BS3}, 
so we omit the detail): 
\begin{itemize}
 \item $\rO_{n+1}\subset \rO_n$ and $\rO_{n+1}$ 
 is dense in $\rO_n$,
 \item Every $f\in\rO_n$ has $2^n$ disjoint compact 
 regular filtrating sets $\cU_{\omega}$, 
 $\omega\in \{1,2\}^n$ such that
 the maps $f\mapsto \cU_\omega$ are locally constant
 (see Section~4.6 of \cite{BS3}),
 \item For any $f \in\rO_{n+1}$  and $\omega\in\{1,2\}^{n+1}$ one has $\cU_\omega \subset \mathring\cU_{[\omega]_n}$, 
 and $\cU_\omega$ is 
 $\cU_{[\omega]_n}$-minimal 
 and $\cU_{[\omega]_n}$-expansive. 
 \item $\cU_\omega$ is $\cU_\omega$-chain transitive for every 
 $\omega\in\{1,2\}^n$.  
 \item The minimum $\cU_\omega$ 
 period is larger than $|\omega|$ for every 
 $\omega\in\{1,2\}^n$,
 \item $\max\{\mathrm{cdiam}(\cU_{\omega}),
 \cdiam(f(\cU_{\omega})))\}<\delta_n$ for every 
 $\omega\in\{1,2\}^n$. 
\end{itemize}
Now $\rR=\bigcap_{n \geq 1}\rO_n$ is 
a $\rG$-residual subset of $\rO_1$, and
Proposition~\ref{p.minimale} and \ref{l.expansive} imply 
that for any $f\in\rR$ and any $\omega\in\{1,2\}^\NN$ the 
compact set $\La_\omega=\bigcap_n\cU_{[\omega]_n}$ is 
a minimal expansive chain recurrence class of $f$ (apply
Proposition~\ref{p.minimale} and \ref{l.expansive} letting $f_n=f$). 
\end{proof}

Finally, we prepare a technical lemma which enables us to 
choose a convenient neighborhood for constructing nested sequences.

\begin{lemm}\label{l.min-sml-nbd}
Let $f\colon X\to X$ be a homeomorphism of a 
locally connected compact metric space $X$ and $\cU\subset X$ be 
a compact set having finite number of connected components.  
Assume that $K\subset \mathring{\cU}$ is 
a totally disconnected compact set such that given any $x\in K$ and 
any connected component $C$ of $\cU$ 
there is $n>0$ such that $f^{n}(x)\in \mathring{C}$. 

Then, there is a 
compact neighborhood $\cO$ of $K$ having finitely many connected components
with the following property: Let $g$ be a homeomorphism of 
$X$ such that the support $\supp (g,f)$ is strictly contained in $\cO$. Then
given any connected components 
$A$ of $\cO$ and $C$ of $\cU$, there is $n>0$ such that 
$g^n(A)\subset \mathring C$. 
 
\end{lemm}
\begin{proof}
Due to the finiteness of the number of connected 
components of $\cU$, we first construct 
the desired neighborhood for a fixed 
connected component $C$ of $\cU$, see the last comment for 
the general case.

First, let us see that there is $N>0$ such that 
for any $x\in K$ there is $0<n\leq N$ satisfying $f^{n}(x)\in \mathring{C}$.
To see this, for each $x\in K$ we take
$n >0$ such that $f^n(x) \in \mathring{C}$. 
Then the same holds in a neighborhood of $x$. 
The compactness of $K$ allows us to cover $K$ by 
finitely many 
open sets where the number $n$ is constant. 
Thus, the number $N$ for the component $C$ is uniformly bounded on $K$. 

By hypothesis, $K$ and
$X \setminus \mathring{\cU}$
has positive distance. We denote it by $\delta >0$.
Now let us choose the neighborhood of $K$.
For each $x \in X$, by $B(\varepsilon, x)$ we denote
the closed ball of radius $\varepsilon$ centered 
at $x$. 

We choose $\varepsilon_1$ and $\varepsilon_2$ as follows:
First, $\varepsilon_1$ is a positive real number which 
is so small that for every $x \in K, i=0,\ldots, N$, 
we have $f^i(B(\varepsilon_1, x)) < \delta/5$. 
To define $\varepsilon_2$, we fix an auxiliary 
number $\eta >0$ satisfying the following: 
If $(z_i)_{i=0,\ldots,N}$ 
is an $f$-$\eta$-pseudo orbit then for any
$i=0,\ldots,N$ we have
\[
d(f^i(z_0), z_i) < \delta/5.
\]
Such an $\eta$ can be chosen by the uniform 
continuity of $f$.

Now we choose  $\varepsilon_2$.
Since we assume that $X$ is locally connected, 
for every $\varepsilon>0$ and $x\in X$ we can choose 
a compact connected neighborhood of $x$ whose diameter is 
less than $\varepsilon$. We denote such a neighborhood by $V(\varepsilon, x)$.
For each $\varepsilon>0$, consider the covering 
$\cup_{x \in K}V(\varepsilon, x)$. 
Note that we can choose finite sub-covering, which we denote by
$\{V(\varepsilon, x_i)\}$.

Then, $\varepsilon_2$ is a positive real number which is so small that 
if we consider the finite covering
$\cO := \cup_{i}V(\varepsilon_2, x_i)$, 
then every connected component of $\cO$ has
diameter less than $\eta$ and $\varepsilon_1/2$. 
Note that the existence of $\varepsilon_2$ is 
guaranteed by the totally disconnectedness of $K$
(if such an $\varepsilon_2$ does not exist, then $K$ must 
contain a non-degenerate continuum) 
and $\cO$ has finitely many connected components
(the number is bounded by the number of $\{V(\varepsilon_2, x_i)\}$).
We assume that $\varepsilon_2 < \varepsilon_1$. 

Consider 
the neighborhood $\cO$ of $K$ defined as above. 
Take a connected component $A$ of $\cO$ 
and $x_0 \in A \cap K$. 
Consider $0\leq n\leq N$ such that $f^n(x_0) \in \mathring{C}$. 
Let us show $g^n(A) \subset \mathring{C}$. 
Let $g$ be a homeomorphism whose support is 
strictly contained in $\cO$. By Lemma~\ref{l.1step}, we have 
$f(A) = g(A)$ for every connected
component of $\cO$. Hence, the $C^0$-distance between $f$ and $g$
is less than $\eta$. Thus, for any $x \in A$
the sequence $\{g^i(x)\}_{i=0,\ldots,N}$ defines 
an $f$-$\eta$-pseudo orbit.
As a result, we know that $d(f^n(x), g^n(x)) < \delta/5$ for any $x \in A$. 

Since the diameter of $A$ is less than $\varepsilon_1/2$,
we know $A \subset B(\varepsilon_1, x_0)$. Thus 
we know the diameter of $f^n(A)$ is less than $\delta/5$ by the choice of $\varepsilon_1$.
By the above observation, we know $g^n(A)$ is 
in the $\delta /5$-neighborhood of $f^n(A)$. 
Thus $g^n(A)$ is contained in the $2\delta/5$-neighborhood of $f^n(x_0)$. As a result, we have 
$g^n(A) \subset \mathring{C}$.

For obtaining $\cO$, we construct $A$ for each $C$. For each $C$ we 
have $\varepsilon_2 =\varepsilon_2(C)$. Then we choose the 
minimum $\varepsilon_2>0$ and take $\cO$. It gives us the desired 
neighborhood.   
\end{proof}

\subsection{Infinitely many ergodic measures for a limit set}
\label{ss.infierg}

In this subsection, we consider a sufficient condition 
which guarantees the existence of infinitely many 
ergodic measures on the limit of nested sequence of filtrating sets. 

Let $(X,d)$ be a compact metric space and $\cP(X)$ be the space of probability measures on $X$ endowed with the distance 
\[
\mathfrak{d}(\mu_1,\mu_2)=\sum_{n=1}^\infty \frac1{2^n}\left|\int \phi_n d\mu_1-\int\phi_n d\mu_2 \right|,
\]
where $\mu_1,\mu_2\in \cP(X)$ and $(\phi_n)$ is a
 sequence of continuous functions bounded by $1$ such that 
 $(\phi_n)$ 
 is dense in the unit ball of the set of continuous functions 
 equipped with the supremum norm. 

We say that a finite set of measures is linearly independent if they are independent as vectors in the space of continuous functionals
on the space of continuous functions on $X$. 

\begin{defi}
Let $k\geq 2$ and $\cM=\{\mu_i\}_{i=1,\ldots,k}$ be a finite set of linearly independent probability measures. We define the \emph{independence radius} of $\cM$, denoted by  
$\rho(\cM)$, as the supremum of the set of real numbers for which we have the following property:
\begin{center}
Any $k$-ple of probability measures 
$\cN=\{\nu_i\}_{i=1,\ldots,k}$ satisfying
$\mathfrak{d}(\mu_i, \nu_i)<r$ for every $i=1,\ldots,k$ 
is linearly independent. 
\end{center}
\end{defi}

\begin{rema}
\label{r.m-ineq}
 For any finite set $\cM$ of $k$ linearly independent probability measures, the independence radius is strictly positive. It varies continuously  with $\cM$.  More precisely, for $\tilde \cM=\{\tilde\mu_i\}_{i=1,\ldots,k}$
satisfying $\mathfrak{d}(\tilde \mu_i,\mu_i)<\varepsilon$
for every $i=1,\dots,k$, we have 
\[
\rho(\cM)+\varepsilon\geq  \rho(\tilde \cM)\geq \rho(\cM)-\varepsilon.
\]
\end{rema}
\begin{proof}[Proof of the inequality.] It is mere an application of 
triangular inequality: 
Consider  $\cN=\{\nu_i\}_{i=1,\ldots, k}$.  
Then for every $i$ we have
 $$\mathfrak{d}(\tilde \mu_i,\nu_i)+\mathfrak{d}(\tilde \mu_i,\mu_i)\geq \mathfrak{d}( \mu_i,\nu_i)\geq \mathfrak{d}(\tilde \mu_i,\nu_i)-\mathfrak{d}(\tilde \mu_i,\mu_i)$$
and thus
 $$\mathfrak{d}(\tilde \mu_i,\nu_i)+\varepsilon > \mathfrak{d}( \mu_i,\nu_i) > \mathfrak{d}(\tilde \mu_i,\nu_i)-\varepsilon.$$
 We now consider the infimum of these quantities 
 when $\cN$ is not independent and we obtain the announced inequality. 
\end{proof}

\begin{lemm} \label{l.inde}
Consider a triangular sequence of probability measures 
$\cM_n=\{\mu^n_1,\dots, \mu^n_n\}\subset \cP(X)$,
$n \geq 1$, such that $\mathcal{M}_n$ is linearly independent
for every $n$.
Suppose that for any $n$, 
$i\in\{1,\dots,n\}$ and $k \geq 1$ one has
\[
\mathfrak{d}(\mu^n_i,\mu^{n+k}_i)<\frac12\rho(\cM_n).
\]
For $i \geq 1$, let $\mu_i$ be an accumulation
point of the sequence $\{\mu^m_i\}_{m\geq i}$. 
Then for any finite $N>0$, the set of 
measures $\{\mu_i\}_{i=1,\ldots, N}$ is
linearly independent. 
\end{lemm}
\begin{proof} Fix any $N>0$.  
Then $\mathfrak{d}(\mu^{N}_i, \mu_i)\leq \frac12\rho(\cM_N)$
for $i\in\{1,\dots,N\}$.  Thus $\{\mu_1,\dots,\mu_N\}$ is independent by definition of $\rho(\cM_N)$
and Remark~\ref{r.m-ineq}.
\end{proof}

Hereafter, $\lambda_0$ denotes a fixed positive real number 
in $(0, 1)$ satisfying 
\[
\prod_{i=1}^{+\infty} (1+\lambda_0^i) \cdot \left(\sum_{i=1}^{+\infty}\lambda_0^i  \right) \leq \frac12. 
\]
Let us give a criterion to have infinitely many distinct ergodic measures 
on an aperiodic class.
\begin{prop}\label{p.nuergodique} Let $\{(f_n, \cU_n)\}_{n \geq 1}$ be a nested sequence for 
an aperiodic class (see Proposition~\ref{p.minimale}). 
Recall that the sequence homeomorphisms $f_n$ converges 
to a homeomorphism $f$ for which $\La=\bigcap_n \cU_n$ is 
a chain recurrence class. 

Assume that for every $n\geq 1$ the filtrating set 
$\cU_n$ contains a family of mutually disjoint periodic orbits 
$\Gamma_n := \{\gamma^n_1,\dots,\gamma^n_n\}$ of $f_n$. 
We denote by $\mu^n_i$ the ergodic probability Dirac measure associated to $\gamma^n_i$. 
Note that the set of measures 
$\cM_n:=\{\mu_i\}_{i=1,\ldots, n}$ are linearly independent.

Suppose that 
for any $f_{n+k}$, $k \geq 1$, 
$\Gamma_n$ is still a family of periodic orbits
and that for any $n \geq 1$ and 
$i\in\{1,\dots,n\}$ we have
\[
\mathfrak{d}(\mu^n_i,\mu^{n+1}_i)<\lambda_0^{n+1}\cdot \rho(\cM_n).
\]
Then there are infinitely many $f$-invariant 
ergodic probability measures on $\La$.
\end{prop}

We prepare a lemma. 
\begin{lemm}\label{l.nuergodique} Under the assumption of Proposition~\ref{p.nuergodique}, for any $n$, $i\in\{1,\dots,n\}$ and 
$k>0$ one has 
\[
\mathfrak{d}(\mu^n_i,\mu^{n+k}_i)<\frac12\cdot \rho(\cM_n).
\] 
\end{lemm}
\begin{proof} The assumption $\mathfrak{d}(\mu^n_i,\mu^{n+1}_i)<\lambda_0^{n+1}\cdot \rho(\cM_n)$ implies $\rho(\cM_{n+1})<(1+\lambda_0^{n+1})\rho(\cM_n)$, and thus, by induction
$$\rho(\cM_{n+k})<\rho(\cM_n)\cdot\prod_{j=1}^{k}(1+\lambda_0^{n+j}).$$

Then, for every $n$, $k>0$ and $i\in\{1,\dots,n\}$ the distance $\mathfrak{d}(\mu^n_i, \mu^{n+k}_i)$ is bounded by 
\begin{align*}
 \mathfrak{d}(\mu^n_i, \mu^{n+k}_i)
&\leq 
\sum_{j=0}^{k-1}\mathfrak{d}(\mu^{n+j}_i,\mu^{n+j+1}_{i})
\leq
\sum_{j=0}^{k-1}\lambda_0^{n+j+1}\cdot \rho(\cM_{n+j})\\ 
&\leq \left(\sum_{j=1}^{+\infty}\lambda_0^{n+j}\right)\cdot\rho(\cM_n)\cdot\left(\prod_{j=1}^{+\infty}(1+\lambda_0^{n+j})\right)
\leq \frac{\lambda_0^n}2 \rho(\cM_n).
\end{align*}
Note that the last inequality follows by the choice of $\lambda_0$.  
\end{proof}

\begin{proof}[Proof of Proposition~\ref{p.nuergodique}]
Fix some $N>0$. Consider $\{\mu_1,\dots, \mu_N\}$ 
where $\mu_i$ is an accumulation point of the 
sequence $\{\mu^m_i\}_{m \geq i}$. 
Note that $\{\mu_i\}$ are 
probability measures supported on 
$\La$ and they are $f$-invariant as they are 
limits of $f_n$-invariant probability measures. 

According to Lemma~\ref{l.nuergodique}, one has 
$$\mathfrak{d}(\mu^N_i,\mu_i)\leq \frac12\cdot \rho(\cM_N).$$
Then Lemma~\ref{l.inde} implies that $\mu_1,\dots \mu_N$ are linearly independent. For each of them, consider its decomposition 
into ergodic measures for $f$.  The measures $\{\mu_i\}$ belong to the vector space generated by these ergodic measures. 
This implies that the vector space has dimension at least $N$. 
Thus there are at least $N$ ergodic measures supported on $\La$. 

As $N$ is any positive number, one deduces that $\La$ supports infinitely many ergodic measures for $f$. 
\end{proof}

\subsection{A criterion and a principle for  minimal non ergodic classes}

In this subsection give a criterion to have minimal chain 
recurrence classes with 
infinitely many ergodic measures.

\begin{prop}\label{p.mini-non-erg}
Let $\{(f_n, \cU_n, \Gamma_n)\}_{n \geq 1}$ be a sequence of 
triples such that
\begin{itemize}
 \item $\{(f_n, \cU_n)\}$ is a nested sequence for an 
 aperiodic class $\Lambda = \cap \, \cU_n$.
 \item $\cU_{n+1}$ is $\cU_{n}$-minimal for every $n \geq 1$.
 \item $\Gamma_n =\{ \gamma^n_i \}_{i=1,\ldots, n}$
 is a set of $n$ distinct periodic orbits in $\cU_n$.
 \item Each member of $\Gamma_n$ is outside $\cU_{n+1}$ and 
 they are still the periodic orbit of the same orbit for $f_{n+1}$ (thus for $f_{n+k}$ for any $k \geq 1$).
 \item For any $i\in\{1,\dots, n\}$ one has  
\[\mathfrak{d}(\mu^n_i,\mu^{n+1}_i)< \lambda_0^{n+1} \rho (\{\mu_i^n \mid i=1,\ldots, n\}),\]
where $\mu_i^n$ is the Dirac probability measure associated to
$\gamma_i^n$. 
 \end{itemize}
Then $(f_n)$ converges to a homeomorphism such that 
$\Lambda$ is minimal and supports infinitely many 
ergodic measures.
\end{prop}
The proof is a direct consequence of Proposition~\ref{p.minimale}
and Proposition~\ref{p.nuergodique}. 

\bigskip

Now we state a version of Proposition~\ref{p.mini-non-erg}
which assures the locally generic coexistence of minimal aperiodic 
classes with infinitely many distinct ergodic measures. 
To state it, we prepare a definition.

Let $f \in \rG$ and $\gamma$ be a periodic orbit of $f$. 
We say that $\gamma$ has a \emph{continuation} in $\rG$ if 
there are a neighborhood $\rO$ of $f$ in $\rG$ and a map
$T: g \mapsto \gamma_g$ defined over $\rO$
such that $T(f) = \gamma$, $T$ is continuous with respect to the Hausdorff distance and $\gamma_g$ is a periodic orbit for $g$.
\begin{defi}
\label{d.mini-infi-prop}
Consider a triple $(f, \cU, \Gamma_n)$,
where $f\in \rG$, $\cU$ is a filtrating set of $f$ 
and $\Gamma_n = \{\gamma_1,\dots, \gamma_n \}$ is 
a set of $n$ distinct periodic orbits of $f$ in $\cU$ such that
each $\gamma_i$ has a continuation in $\rG$.

A family of $\rG$-robust 
properties $(\cQ^n)_{n \geq 1}$ for $(f, \cU, \Gamma_n)$
is said to be \emph{type $\fP_{M, \infty}$} if the following holds:
If $(f, \cU, \{\gamma_1,\dots,\gamma_n\})$, 
satisfies $\cQ^n$, then for any $\delta >0, N>0$ and 
any neighborhood $\rO$ of $f$ in which the continuations of orbits of
$\Gamma_n$ are defined, 
there is $g\in \rO$ such that the following holds:
\begin{itemize}
 \item The support $\supp (g, f)$ is strictly contained in $\cU$.
 \item There are disjoint regular filtrating sets 
 $\cU_1,\cU_2 \subset \mathring{\cU}$ for $g$ 
 which are both $\cU$-minimal for chains.
 \item For $j=1, 2$, $\cU_j$ has finitely many connected 
 components and $\cU_j$-chain transitive.
 \item $\cdiam(\cU_j)$, $\cdiam(g(\cU_j)) < \delta$ for $j\in\{1,2\}$.
 \item For $j=1, 2$, the minimum period of $\cU_j$ is larger than $N$.
 \item There are periodic orbits $\gamma_i^j\subset \cU_j$,$j\in\{1,2\}$,  $i\in\{1,\dots, n+1\}$
 such that $(g, \cU_j, \{\gamma_1^j,\dots,\gamma_{n+1}^j\})$ satisfies the property $\cQ^{n+1}$ for $j\in\{1,2\}$.
 \item For $j\in\{1,2\}$ and $i\in\{1,\dots, n\}$ one has  
\[\mathfrak{d}(\mu^g_i,\mu^{j}_i) < 
\lambda_0^{n+1} \rho (\{\mu^g_i \mid i=1,\ldots, n\})\]
where $\mu^g_i,\mu_i^j$ are the Dirac probabilities associated to 
$\gamma_{i, g}$ (the continuation of $\gamma_i \in \Gamma_n$) 
and $\gamma_i^j$, respectively. 
 \end{itemize}
\end{defi}

In the next section, we give an example of a family of local 
properties $(\cQ^n)$ which is type $\fP_{M, \infty}$, see Proposition~\ref{p.C1minnergo}. 
Let us see the consequence of such a property.

\begin{prop}\label{p.principle-genurgo}
Let $(\cQ^n)_{n \geq 1}$ a family of $\rG$-robust
properties which is type 
$\fP_{M, \infty}$.
Assume that $(f, \cU, \gamma_1)$ satisfies property $\cQ^1$ 
and it persists over a $\rG$-open 
neighborhood $\rO$ of $f$. Then there is a 
$\rG$-residual subset $\rR$ of $\rO$ such that 
every $g \in \rR$ has uncountably many aperiodic 
classes which are minimal and supports infinitely many ergodic 
probability measures.
\end{prop}

The proof is similar to 
the proof of Proposition~\ref{p.principle-expansive}.
Let us give it.

\begin{proof}
The fact that $(\cQ^n)$ is type $\fP_{M,\infty}$ allows us to build by induction, a decreasing sequence of open subsets $(\rO_n)$ and a sequence $\delta_n>0$ tending to $0$ as $n\to\infty$ with the following property: 
\begin{itemize}
 \item $\rO_n$ is dense in $\rO$ for every $n$. 
 \item Every $f\in\rO_n$ has $2^n$ disjoint regular filtrating sets $\cU_{\omega}$, $\omega\in \{1,2\}^n$.
 \item $\cU_{\omega}$ is $\cU_{\omega}$-transitive 
 for every $\omega\in \{1,2\}^n$.
 \item For any $f \in\rO_{n+1}$  and $\omega\in\{1, 2\}^{n+1}$ one has 
 $\cU_\omega \subset \mathring\cU_{[\omega]_n}$ and $\cU_\omega$ is 
 $\cU_{[\omega]_n}$-minimal for chains. 
 \item $\max\{\mathrm{cdiam}(\cU_{\omega}),
 \cdiam(f(\cU_{\omega}))\} < \delta_{|\omega|}$. 
 \item The minimum period of $\cU_\omega$
 is larger than $|\omega|$ for every $\omega\in\{1, 2\}^{n}$. 
  \item The maps $f\mapsto \cU_\omega$ are locally constant
 and there are $n$ periodic orbits 
 $\Gamma_{n, \omega} = \{\gamma^{\omega}_i\}_{i=1,\ldots,n}$ in 
 $\cU_{\omega}$.
 Note that $\Gamma_{n, \omega}$ may be perturbed to construct a new
 chain recurrence classes, but the size of the perturbation can be chosen
 arbitrarily small and it is only once (to construct $\Gamma_{n+1}$) 
 so $\Gamma_{n, \omega}$ is a well defined object.
\item For every $n$ and $\omega \in \{ 0, 1\}^{n+1}$ we have
\[\mathfrak{d}(\mu^{n, [\omega]_n}_i, \mu^{n+1, \omega}_{i})<\lambda_0^{n+1}\cdot \rho(\{ \mu^{n, [\omega]_n}_i \mid i=1,\ldots, n\}).\]
\end{itemize}
Now $\rR=\bigcap_{n \geq 1}\rO_n$ is 
a residual subset of $\rO$ and 
it has uncountably many aperiodic minimal chain recurrence classes. 
By proposition~\ref{p.nuergodique} we know that 
they support infinitely many ergodic measures.
\end{proof}

\subsection{A criterion for unique ergodicity}
\label{ss.unique}

Suppose that we have a compact set $\cU$ of $X$. 
The set of Dirac measures of $\cU$-pseudo orbits, 
denoted by $\cP_{k, \mathrm{pseudo}}(\cU)$ is the 
collection of measures
\[
\left\{ \frac{1}{k} \sum_{i=1}^k \delta_{x_i}\, \middle\vert \,(x_i)_{i=1,\ldots,k} \mbox{ is a } \cU\mbox{-chain of points of length $k$.} \right\} 
\] 
where $\delta_x$ denotes the Dirac probability measure supported at $x$.
By $\cP_{\infty, \mathrm{pseudo}}(\cU)$, we denote 
the set of accumulation points of sequences $(\nu_i)$
where $\nu_i \in \cP_{i, \mathrm{pseudo}}(\cU)$.

\begin{lemm}\label{l.urgo} Let $\{(f_n,\cU_n)\}_{n\geq 1}$
be a nested sequence for an aperiodic class.  Assume that there is a sequence 
 of positive real numbers $(\varepsilon_n)$ tending to $0$ as $n\to\infty$ such that, 
 for every $n$ one has
\begin{center}
if $\mu_1$ and $\mu_2$ are probability measures which are accumulated by the Dirac probabilities along 
$\cU_n$-pseudo orbits (that is, 
if $\mu_1, \mu_2 \in \cP_{\infty, \mathrm{pseudo}}(\cU_n)$)
then $\mathfrak{d}(\mu_1,\mu_2)<\varepsilon_n$.
\end{center}
Then, the sequence $(f_n)$ converges to a homeomorphism $f$ such that $\La=\bigcap \cU_n$ is an uniquely ergodic, aperiodic chain recurrence class. 
\end{lemm}
\begin{proof}The proof of the fact that the sequence $(f_n)$ converges to a homeomorphisms for which $\La$ is an 
aperiodic class is 
due to Proposition~\ref{p.minimale}. 
Let us prove the unique ergodicity. 

The $f$-orbits in $\La$ are $\cU_n$-pseudo orbits for $f_n$, for every $n$.  Thus the $\mathfrak{d}$-distance between any two probability measures accumulated by Birkhoff averages of the Dirac measure along $f$-orbits 
in $\La$ is less than $\varepsilon_n$ for every $n$. 
In other words, they are equal. 
Therefore, there is a unique measure accumulated by these Birkhoff averages along $f$ orbits in $\La$. It implies the unique ergodicity of $\La$.
\end{proof}

\begin{rema}\label{r.robust} 
The hypothesis for each fixed $n$ is robust in the $C^0$-topology,
if the number of connected components is finite. 
Also, it is stable for local perturbations, since local perturbations 
do not change the set of $\cU_n$-chains. 
\end{rema}

\subsection{A criterion for transitivity}

In this subsection, we discuss conditions which guarantees 
transitivity of limit aperiodic classes. We say that a compact 
$f$-invariant set 
$K \subset X$ is \emph{transitive} if it has a point 
having a dense orbit in $K$. 

\begin{defi}\label{d.trachai}
Let $f$ be a homeomorphisms and $\cU_0$  and $\cU_1$ be two compact filtrating sets such that $\cU_1\subset \mathring\cU_0$ 
and $\cU_i$ is $\cU_i$-chain transitive for $i=0,1$. 
Let $C$ be a connected component of $\cU_1$.
We say that $(\cU_1, C)$ is \emph{$\cU_0$-transitive} 
(or just $\cU_1$ is $\cU_0$-transitive)
if the following conditions hold:
\begin{itemize} 
\item For any component $A$ of $\cU_0$ there is $i^+(A), i^-(A) \geq 0$ 
such that 
$f^{i^+(A)}(C)\subset \mathring A$ and $f^{-i^-(A)}(C)\subset \mathring A$. 
\item For every connected component $C_0$ of $\cU_0$, there is a connected
component $C_1$ of $\cU_1$ such that $C_1 \subset \mathring{C}_0$.
\end{itemize}
In the following, $i^+(A), i^-(A)$ denote the smallest non-negative integers 
satisfying the condition above.
\end{defi}

\begin{lemm}\label{l.stable-transitive} 
If the number of the connected components of $\cU_0$ is finite,
then $\cU_0$-transitivity for some $C$ is a $C^0$-robust property: Any homeomorphism $g$ which is enough $C^0$-close 
to $f$ has $\cU_0$ and $\cU_1$ as filtrating sets 
and $\cU_1$ is $\cU_0$-transitive. 
\end{lemm}
\begin{proof}
Recall that being a filtrating set is a $C^0$-robust
property. Thus we know that $\cU_1$, $\cU_0$ are filtrating 
sets for homeomorphisms sufficiently close. To see the 
$C^0$-robustness of the transitivity, there are only finitely 
many conditions which we need to require and they 
are obviously $C^0$-robust. Thus we have the conclusion.
\end{proof}

\begin{rema}\label{r.travel}
If the number of connected components of $\cU_0$ is finite, 
then there is $N^+, N^- >0$ 
such that $\{f^i(C)\}_{i=0,\ldots,N^+}$,
$\{f^i(C)\}_{i=-N^-,\ldots,0}$
passes through all the connected components of $\cU_0$.
\end{rema}
\begin{proof}
Set $N^{\pm}$ to be the maximum of $i^{\pm}(A)$ respectively. 
\end{proof}

While $\cU_0$-transitivity is a robust property, 
it is not stable by perturbations supported on $\cU_1$. 
Below we give a sufficient condition for stable $\cU$-transitivity.

\begin{lemm}
\label{l.stab-tran}
For $\cU_1 \subset \mathring{\cU}_0$ and $f$ as above,
assume that $(\cU_1, C)$ is $\cU_0$-transitive.
Furthermore, assume that the following condition holds:
\begin{itemize}
\item For any $\cU_1$-chain of connected components
$(U_k)_{k=0,\ldots,i_0}$ satisfying $U_0 = C$
where $0\leq i_0 \leq i^+(A)$, we have $f^{i^+(A)-i_0}(U_{i_0}) \subset \mathring A$.
\item The same holds for $f^{-1}$ and $i^{-}(A)$: 
For any $\cU_1$-chain of connected components
$(U_k)_{k=-i_0,\ldots, 0}$ satisfying $U_0 = C$
where $0\leq i_0 \leq i^-(A)$, we have $f^{-(i^-(A)-i_0)}(U_{-i_0}) \subset \mathring{A}$.  
\end{itemize}
Then $(\cU_1, C)$ is 
stably $\cU_0$-transitive for local perturbations
whose supports are strictly contained in $\cU_1$.
\end{lemm}

\begin{proof}
Consider $g$ whose support is strictly contained in $\cU_1$.
First consider the condition about the forward iteration.
Note that the condition about the support guarantees that the set of 
$\cU_1$-chains of connected components 
are the same for $f$ and $g$.
Given $x \in C$ and a connected component $A$ of $\cU_0$, consider the orbit 
$\{g^i(x)\}_{i=0,\ldots, i^+(A)}$. If the whole orbit is
contained in $\cU_1$, it defines a $\cU_1$-chain of connected 
components of $g$, which appears
in the set of chains connected components of $f$. 
Since the sets of the chains of connected components are the same for $f$ and $g$,
we know $g^i(x) \in \mathring{A}$ (just consider the case $i_0 = 0$). 
If not, there is $g^k(x)$ such that 
$g^k(x) \not\in \cU_1$. We choose $k$ 
to be the minimum positive integer among such numbers. 
By definition we know there is 
$y \in \cU_1(g^{k-1}(x))$ such that $f(y) = g^k(x)$.
Since $\cU_1$ is a filtrating set, we know that 
$g^{k+i}(x) \not \in \mathring{\cU}_1$ for $i>0$. 
Then, since outside $\cU_1$ we have $g \equiv f$, 
we conclude $g^{i^+(A)-k}(g^k(x)) = f^{i^+(A)-k}(f(y))$.
By definition, we know that this point belongs to $\mathring{A}$. 
Thus, we are done.

Now, let us consider the case for the backward iterations. 
One point which differs from the forward case is the domain of the 
support. For $g^{-1}$, we only have 
$\supp(g^{-1}, f^{-1}) \subset f(\cU_1) = g(\cU_1)$. 
Take $x \in C$. If $\{g^{-i}(x)\}_{i=0, \ldots, i^-(A)}$ is contained in 
$\cU_1$, then we can conclude $g^{-i^-(A)}(x) \in \mathring{A}$
as in the forward case. If not, then we choose minimum $k\geq 1$ such that 
$g^{-k}(x) \not\in \cU_1$.
Now, note that $g^{-k+1}(x) \in \cU_1 \setminus (\cU_1 \cap g(\cU_1))=
\cU_1 \setminus (\cU_1 \cap f(\cU_1))$.
This means $\{g^{-k-l}(x)\}_{l\geq 0}$ is outside $f(\cU_1)$,
due to the filtrating property of $f(\cU_1)$ with respect to $g$.
Thus we have $g^{-k-l}(x) = f^{-l}(g^{-k}(x))$ for $l\geq 0$ 
and the conclusion follows by considering $l = i^{-}(A) -k$. 
\end{proof}

\begin{rema}\label{r.statran}
The assumption of Lemma~\ref{l.stab-tran} implies the 
first condition of $\cU$-transitivity in Definition~\ref{d.trachai}
(considering the case $i_0=0$).
Thus, the second condition of Definition~\ref{d.trachai}
and the assumption of Lemma~\ref{l.stab-tran} imply
the stable $\cU$-transitivity for local perturbations.
\end{rema}


\subsection{A criterion for non-minimality}
The next lemma is a criterion for building transitive aperiodic classes admitting several minimal sets: 
\begin{lemm}\label{l.transitive} Let $\{(f_n,\cU_n)\}_{n \geq 1}$
be a nested sequence for an aperiodic class.
Assume that for every $n$ one has
\begin{itemize}
 \item $(\cU_{n+1},C_{n+1})$ is stably $\cU_n$-transitive for $f_{n+1}$
 where $C_{n+1}$ is a connected component of $\cU_{n+1}$.
\item The number of connected components of $\cU_n$ is finite.
\end{itemize}

Then the sequence $(f_n)$ converges to a homeomorphism $f$ for 
which $\La=\bigcap\cU_n$ is a transitive 
chain recurrence class. 

If, furthermore, there are two disjoint compact sets $A, B$ such that for every $n$, there are two closed $\cU_n$-pseudo 
periodic orbits for $f_n$ contained in 
$\mathring{A}$ and $\mathring{B}$ respectively, 
then $\La$ contains at least two minimal sets. 
 
\end{lemm}
\begin{proof}According to Proposition~\ref{p.chainrecurrence} the sequence 
$(f_n)$ converges to $f$ and $\La$ is an aperiodic 
chain recurrence class of $f$. Let us prove the transitivity.
Consider $(\cU_n, C_n)$ for $f_n$ and 
we fix $N^+>0$ by Remark~\ref{r.travel}. 
Note that by definition of $\cU_n$-transitivity, $C_n$ contains 
at least one connected component of $\cU_{n+1}$.
Inductively, we know that $C_n \cap \Lambda \neq \emptyset$.
Take $x \in C_n \cap \Lambda$.
By definition and the assumption about the support,
 $\{x,f(x),\dots, f^{N^+}(x)\}$ meets every component of $\cU_{n-1}$.  
If $n$ is large enough these components are arbitrarily small, say, smaller that 
$\varepsilon>0$ so that the positive orbit 
of $x$ is $\varepsilon$-dense in $\Lambda$.
The same holds for its backward orbits. 
 
Thus given any two open sets $U$, $V$ of $\La$
there is $x \in X$ whose positive orbit meets both $U$ and $V$ and so does its negative orbit. This means that the positive iteration of $U$ meets $V$
and negative iteration of $U$ meets $V$, too.
This implies the topological transitivity of $f$, see Remark~\ref{r.tergo}.  
 
Let us show that $A$ and $B$ contain at least one 
 minimal set. 
 By $(x^k_i)$ we denote the pseudo periodic orbit in $\mathring{A}$
 in $\cU_k$.
 Consider the set of accumulation points of these orbits
 (that is, the accumulation of points of sequences $(y_n)$
 where $y_n$ is one of $(x^n_i)$).
 We denote it by $\Lambda_A$. One can check that it 
 is compact, $f$-invariant and contained in $A$. Thus, 
 there is a minimal set which is strictly contained in $A$.
 The same holds for $B$, which concludes the proof.
\end{proof}
\begin{rema}\label{r.tergo}
The condition about $U$ and $V$ in the proof implies a property 
stronger than the transitivity: It implies the existence of a point whose 
$\alpha, \omega$-limit sets coincide with $\Lambda$.
\end{rema}
\begin{rema}\label{r.nonmini}
If $\cU_n$ is regular, then having a pseudo periodic orbit inside a 
fixed compact set is a $C^0$-robust property.
\end{rema}

%


We can now express our principle for the local genericity of transitive non minimal aperiodic classes. 

Let $\cV$ be a possibly non-filtrating set and
$A$ a compact set such that every connected component of $\cV$ is either
being contained in $\mathring{A}$ or disjoint of $A$.
Then, by the \emph{restriction of $\cV$ to 
$A$}, denoted by $\cV|_A$, 
we mean the compact set $\cV \cap A$.

\begin{defi}\label{d.tra-nonmi}
Let $A, B$ be compact disjoint subsets of $X$, 
$f \in \rG$ and $\cU$ a filtrating set of $f$.

A property $\cP$ on $(f, \cU)$ is called \emph{type 
$\fP_{\mathrm{TnM}, A, B}$} if the following holds: 
For any $N >0$, $\delta >0$, $\rG$-neighborhood $\rO$ of $f$ there is 
$g \in \rO$
such that:
\begin{itemize}
 \item $\supp(g,f)$ is strictly contained in $\cU$.
 \item There are disjoint regular filtrating sets 
 $\cU_1,\cU_2 \subset \mathring{\cU}$ for $g$. 
 \item The number of connected components of 
 $\cU_1,\cU_2$ are finite.
 \item $\cU_i$ is $\cU_i$-chain transitive for $i=1, 2$. 
 \item $\max \{\cdiam(\cU_i), \cdiam g(\cU_i)\}<\delta$ for $i=1, 2$.
 \item Minimum periods of $\cU_i$ are larger than $N$ for $i=1, 2$.
 \item $(\cU_i, C_i)$ is $\cU$-transitive for $i=1, 2$ (where 
 $C_i$ is some connected component of $\cU_i$). 
 \item $\cU_i|_A$ (resp. $\cU_i|_B$) can be defined and
  $A$ (resp. $B$) contains a $\cU_i|_A$-pseudo periodic orbit
  for $i=1, 2$.
 \item $(g, \cU_i)$ satisfies property $\cP$ for $i=1, 2$.
\end{itemize}
\end{defi}

Then the following is a direct corollary of 
Lemma~\ref{l.transitive}.

\begin{prop}\label{p.trnomin}
Suppose that $(f, \cU)$ satisfies a $\rG$-robust 
property $\cP$ which is type $\fP_{\mathrm{TnM},A, B}$ for some two 
disjoint compact sets $A$ and $B$.
Then there is a neighborhood $\rO_0$ of $f$ in $\rG$ and a residual subset 
$\rR\subset \rO_0$ such that every $g\in\cR$ admits an uncountable family of aperiodic transitive classes containing at least two minimal sets.  
\end{prop}

The proof is almost the same as Proposition~\ref{p.principle-expansive}
or Proposition~\ref{p.principle-genurgo}: We need to construct a branching 
family of nested sequences for aperiodic classes satisfying the conditions 
of Lemma~\ref{l.transitive}. So we omit it.

\subsection{A criterion for non-transitive class}

Let us discuss a condition which guarantees the non-transitivity
of limit aperiodic classes.

 \begin{lemm}\label{l.nontr} 
 Let $\{(f_n,\cU_n)\}_{n \geq 1}$ 
 be a nested sequence for an aperiodic class.  
Assume that there are two sequences $(A_n)$ and $(B_n)$ of connected components of $\cU_n$ such that for every $n$ one has
\begin{itemize}
 \item $A_{n+1}=A_n \cap \cU_{n+1}$, $B_{n+1}=B_n\cap \cU_{n+1}$,
 \item any $\cU_n$-chain of connected components
 of length $n+1$ starting from or ending at $A_n$
 do not contain $B_n$. 
 \end{itemize}
Then the sequence $(f_n)$ converges to a homeomorphism $f$ 
such that $\La=\bigcap \cU_n$ is a non-transitive,
aperiodic chain recurrence class. 
\end{lemm}
\begin{proof}
First, by Proposition~\ref{p.chainrecurrence} 
we know that $\La := \bigcap_n \cU_n$ is an aperiodic chain 
recurrence class. We will prove that $\La$ contains 
two points $a, b$ such that they are isolated in $\Lambda$
and have distinct orbits. This implies the non-existence 
of the point $p \in \La$ such that $\overline{\cO(p)} = \La$.
Indeed, if $p \in \La$ satisfies $a\in \overline{\cO(p)}$ then 
the assumption that $a$ is isolated implies 
$\cO(p) = \cO(a)$. Similarly, $b\in \overline{\cO(p)}$ implies
$\cO(p) = \cO(b)$. However, these cannot hold simultaneously if 
the orbits of $a$ and $b$ are different.
 
Now, let us prove above assertion.
The nestedness of $(A_n)$ and $(B_n)$, together
with $\cdiam (\cU_n) \to 0$ imply that
$\bigcap_n A_n$, $\bigcap_n B_n$ exist and they 
are singletons. 
We denote them by $a$ and $b$ respectively.
The definitions of $(A_n)$ and $(B_n)$
imply that $a, b$ are isolated in $\Lambda$.

Recall that $f_n$ and $f_{n+k}$
coincides outside $\cU_n$ for every $k\geq 0$.
This implies that for every $n$ and $k$ 
the sets of $\cU_n$-chains for $f_n$ and $f_{n+k}$
are the same. In particular, they are same for $f_n$ and $f$.
Suppose $\cO (a) \cap \cO (b) \neq \emptyset$ for $f$.
Then, it implies there is $K \in \ZZ$ such that $f^K(a) =b$,
but this contradicts the fact that
in $\cU_K$ we do not have $\cU_K$-chain of 
length $K$ between  
$A_K$ and $B_K$.
\end{proof}

\begin{rema}\label{r.nontra}
The non-existence of a certain $\cU_k$-chain 
of finite length
is a $C^0$-robust property,
if the number of the 
connected components of $\cU_k$ is finite. 
The condition is stable under a local perturbation whose 
support is strictly contained in $\cU_k$, because such a perturbation does not change
the set of the $\cU_k$-chains.
\end{rema}

\begin{rema}
By carefully repeating the above proof one can see that $f$ has 
no point whose $\omega$-limit set (resp. $\alpha$-limit set)
coincides with $\Lambda$.
\end{rema}

\subsection{A principle for the local genericity of non-transitive uniquely ergodic classes}
Now we can give a principle which ensures the 
generic existence of non-transitive uniquely ergodic 
aperiodic classes.

\begin{defi}
Let $A$, $B$ be compact disjoint sets of $X$.
Let $\cP$ 
be a $\rG$-robust property on $(f, \cU)$ 
where $f\in \rG$, $\cU$ is a compact filtrating set of $f$
such that $\cU$ has unique connected components 
contained in $\mathring{A}$ and $\mathring{B}$ respectively. 
We say that a property $\cP$ is type $\fP_{\mathrm{nT}, A, B}$ 
if the following holds:
For any $\delta_1>0, \delta_2>0, N>0, M>0$  
and 
any neighborhood $\rO$ of $f\in \rG$ there is $g\in \rO$ such that:
\begin{itemize}
 \item $\supp (g, f)$ is strictly contained in $\cU$.
 \item There are disjoint regular filtrating sets $\cU_1,\cU_2 \subset \mathring{\cU}$ for $g$ which are $\cU_1$, $\cU_2$-chain transitive
 respectively. 
 Their minimal periods are larger than $N$. 
 Furthermore, the number of connected components of 
 $\cU_i$ is finite for $i=1,2$.
 \item  $\max \{\cdiam(\cU_i), \cdiam g(\cU_i)\}<\delta_1$ for $i=1,2$. 
 \item The ergodic diameter of $\cU_i$ is less than $\delta_2$ for $i=1, 2$ (see Section~\ref{ss.unique}): 
 If $\mu_{i}$ and $\mu'_{i}$  are probability measures which are accumulated by the Dirac probabilities along $\cU_i$-pseudo orbits, then 
 \[
 \mathfrak{d}(\mu_{i},\mu'_{i})<\delta_2.
 \]
 \item For $i=1, 2$, $\cU_i$ has two connected components $A_i$, $B_i$
 which are unique components contained in the connected component 
 $A_0$, $B_0$ of $\cU$ contained in $A$, $B$ respectively.
 \item Any $\cU_i$-chain of connected components
 of length $M$
 starting from or ending at $A_i$ do not contain $B_i$ for $i=1, 2$.
 \item  $g$ satisfies property $\cP$ on $\cU_i$ for $i=1, 2$.
\end{itemize}
\end{defi}

We have the following principle.

\begin{prop}\label{p.principle-nontrurgo}
Suppose that $(f, \cU)$ satisfies a $\rG$-robust property $\cP$ 
which is of type $\fP_{\mathrm{nT}, A, B}$ 
for some compact disjoint sets $A$ and $B$.
Then, there is a neighborhood $\rO_0$ of $f$ in $\rG$ and a residual subset 
$\rR\subset \rO_0$ such that every $g\in\rR$ admits an uncountable family of aperiodic chain-recurrent classes which are not transitive but uniquely ergodic. 
\end{prop}
\begin{proof}As $\fP_{\mathrm{nT}, A, B}$ is $\rG$-robust, there is a non-empty open neighborhood $\rO_1$  of $f$
such that every $f_1 \in \rO_1$ satisfies $\fP_{\mathrm{nT}, A, B}$ 
on $\cU$.

Then the assumption of Proposition~\ref{p.principle-nontrurgo} 
allows us to build a sequence of open subsets $\rO_n\subset \rO_1$ and  sequences $\delta_n>0$, $\varepsilon_n>0$ tending to $0$ as $n\to\infty$ with the following property (see also Remark~\ref{r.robust} and Remark~\ref{r.nontra}): 
\begin{itemize}
 \item $\rO_{n+1}\subset \rO_n$ and $\rO_n$ is dense in $\rO_n$.
 \item Any $f\in\cO_n$ has $2^n$ disjoint
 regular filtrating sets $\cU_{\omega}$ ($\omega\in \{1,2\}^n$)
 having finitely many connected components 
 such that the map $f\mapsto \cU_\omega$ is locally constant. 
  \item $\max\{\cdiam(\cU_\omega),\diam (g(\cU_\omega))\}<\delta_n$ for every $\omega\in\{1,2\}^n$. 
 \item Each $\cU_{\omega}$ is $\cU_{\omega}$-chain 
 transitive and has its minimum period larger than $|\omega|$. 
 \item For any $f\in\rO_{n+1}$ and $\omega\in\{1, 2\}^{n+1}$ one has 
 $\cU_\omega\subset \mathring\cU_{[\omega]_n}$.
 \item For any $f\in\rO_{n+1}$  and $\omega\in\{1,2\}^{n+1}$ one has 
 two connected components $A_\omega$ and $B_\omega$ of $\cU_\omega$ 
 such that 
 $$A_\omega= \cU_\omega\cap A_{[\omega]_n}\mbox{ and }A_\omega= \cU_\omega\cap A_{[\omega]_n}.$$
 \item For any $f\in\rO_{n}$ and $\omega\in\{1, 2\}^{n}$, 
every $\cU_{\omega}$-chain of connected components of length $|\omega|$
starting from or ending at $A_{\omega}$ does not contain $B_\omega$. 
 \item If $\mu_{1,\omega}$ and $\mu_{2,\omega}$ are probability measures which are accumulated by the Dirac probabilities along $\cU_\omega$-pseudo orbits ($\omega\in\{1,2\}^n$), then 
 $$\mathfrak{d}(\mu_{1,\omega},\mu_{2,\omega})<\varepsilon_n.$$
\end{itemize}

Now $\rR=\bigcap\rO_n$ is a residual subset of $\rO_0$ and
 Lemma~\ref{l.urgo} implies that for any $f\in\rR$ and any 
 $\omega\in\{1,2\}^\NN$ the compact set $\La_\omega=\bigcap_n\cU_{[\omega]_n}$ is a uniquely ergodic chain recurrence class of $f$,  and Lemma~\ref{l.nontr} implies that 
it is not transitive. Thus the proof is completed.
\end{proof}


\section{Principles for the local $C^1$-genericity of aperiodic classes with prescribed dynamics}
\label{s.C1}

In this section, we discuss the existence of 
properties proposed in the previous section and 
complete the proof of Theorem~\ref{t.minimal-non}.
In this section $M$ denotes a 
smooth closed connected three dimensional manifold
and $\rG = \mathrm{Diff}^1(M)$
with the $C^1$-topology.

\subsection{Basic properties}
In this subsection we investigate several basic properties
of filtrating Markov partitions. 

In the course of the construction we use 
properties of finiteness of the number of connected 
components of Markov partitions and regularity. 
Note that by definition of the filtrating Markov partition 
(of saddle type) we have the following properties (see \cite[Definition 2.5]{BS2}):
\begin{itemize}
\item Filtrating Markov partitions have at most finitely many 
connected components.  
\item Filtrating Markov partitions are regular:
Given two connected components $U_1, U_2$
of a Markov partition $\cU$ of a diffeomorphism $f$, 
if $f(U_1) \cap U_2 \neq \emptyset$
then each connected component of the intersection
is a vertical cylinder.
Thus their interiors also have non-empty intersection. 
\end{itemize}

%

We prepare some fundamental results 
about filtrating Markov partitions. 
The following lemma is given in \cite[Section~2.2]{BS3}. So we omit the proof.
For the definition of $\cU_{(m, n)}$ and
$\cU_{(m, n)}(K)$, see Section~\ref{s.preliminar} 
and~\ref{ss.tool}.

\begin{lemm}\label{l.mrk-sml}
Let $\cU \subset M$ be a filtrating Markov partition
of $f \in \mathrm{Diff}^1(M)$
and $K$ be a circuit of points 
such that every periodic point
has a large stable manifold. 
Then the 
Markov partition $\cU_{(m, n)}(K)$ converges to $K$
as $m, n \to \infty$. More precisely, 
for any neighborhood $O$ of $K$ there exists 
$m_0, n_0>0$ such that for every $m\geq m_0, n\geq n_0$
we have $\mathcal{U}_{(m, n)}(K) \subset O$.
\end{lemm}

As a consequence, we have 
the following (see also \cite[Section~4.4]{BS3}).
\begin{lemm}\label{l.mini-peri}
Let $\cU$ be a filtrating Markov partition
and $K$ be a circuit of points such that every 
periodic point has its period larger than $\pi$ and 
a large stable manifold. 
Then, for sufficiently large $m$ and $n$ the 
Markov partition $\cU_{(m,n)}(K)$ has minimum 
period larger than $\pi$ as well.
\end{lemm}

\begin{proof}
When $m$ and $n$ are large, $\cU_{(m, n)}(K)$ is 
close to $K$. Consider a pseudo periodic orbit in $\cU_{(m, n)}(K)$.
If it follows the periodic orbit, then the period must 
be larger than $\pi$. If not, then it must follows one 
of a homo/heteroclinic orbit, but when $m$ and $n$ are large,
then
the number of iterations needed to follow the homo/heteroclinic orbits
to arrive at near the initial point must be large, in 
particular larger than $\pi$. 
Thus in both cases we have the conclusion.
\end{proof}

Finally, let us observe the following. 

\begin{rema}\label{r.peri-sub}
Let $\cU$ be a filtrating Markov partition
and $\cV$ be a filtrating Markov partition 
matching to $\cU$ (see Definition~\ref{d.macha}).
If the minimum period of $\cU$ is larger than 
$\pi$, then so is for $\cV$.
In other words, the property ``having minimum period
larger than $\pi$'' inherits to matching 
Markov partitions.
\end{rema}

\begin{proof}
Consider a pseudo periodic orbit of $\cV$. 
Then the fact that $\cV$ is matching to $\cU$
enables us to take the lift of the pseudo periodic 
orbit to $\cU$ with the same period. 
Hence, the pseudo orbit of $\cV$ must
have the period larger than $\pi$. 
\end{proof}

The following lemma enables us to find a 
homo/heteroclinic point.

\begin{lemm}\label{l.homoclinic}Let $\cU$ be a 
sub Markov partition of some filtrating Markov partition 
and let $p$ be a periodic point whose orbit is contained in $\cU$ having a large stable manifold. 
If we have a finite sequence of $\cU$-chains of points $(x_i)_{i=0,\ldots,n}$
such that $x_0$ and  $x_n$ belong to the orbit of $p$,
it is $\cU$-shadowed by 
an orbit of the point of transverse intersection 
between $W^s(p)$ and $W^u(p)$, that is, there is a point $y$ of the transverse intersection 
between $W^s(p)$ and $W^u(p)$ 
such that $x_i$ and $f^i(y)$ belong
to the same component of $\cU$ for every 
$i=0,\ldots,n$.  
\end{lemm}

\begin{proof}
By the existence of the chain $(x_i)$, we know that 
there is a segment of the local unstable manifold 
$\sigma \subset W^u_{\mathrm{loc}}(p)$ such that 
$f^i(\sigma)$ belongs to the rectangle $\cU(x_i)$
and $f^n(\sigma)$ is a segment which properly 
crosses $\cU(p)$ in the vertical 
direction (see also Lemma~2.15 of \cite{BS2}). 
Thus the largeness of the 
stable manifold of $p$ assures there exists a point
$y \in W^s(p) \cap \sigma \neq \emptyset$. 
This gives the desired point.
\end{proof}

\begin{rema}\label{r.nontri}
It may be that the point obtained is $p$ itself. 
However, if we know that one of the rectangles 
$\cU(x_i)$ does not contain the orbit of $p$, 
then we can conclude that it is a homoclinic point.
In that case, $y$, $f^n(y)$
belong to $W^u_{\mathrm{loc}}(p)$, 
$W^s_{\mathrm{loc}}(p)$ respectively.
\end{rema}

As a simple consequence of Lemma~\ref{l.homoclinic}
we obtain the following:
\begin{coro}\label{c.homo-ubiq}
Let $\cU$ be a transitive filtrating
Markov partition containing 
a periodic point $p$ with a large stable manifold. 
Then the relative homoclinic class $H(p, \cU)$
has non-empty intersection with every rectangle 
of $\cU$.
\end{coro}

\subsection{Expansiveness of sub Markov partitions}

\subsubsection{Expansiveness and circuits}

Given a possibly non-filtrating
Markov partition $\cU$ we say that $\cU$ is \emph{generating} if 
given any components $C_1,C_2$ of $\cU$,
$f(C_1)\cap C_2$ is either empty or connected.  The following is given in Section~2.3 of \cite{BS3}.
\begin{rema}\label{r.expa-sub}
Given a filtrating Markov partition 
$\cU$ and $(m, n) \neq (0, 0)$, 
its $(m,n)$-refinement is generating. 
\end{rema}

Recall that we defined the notion of $\cU$-expansiveness,
see Definition~\ref{l.def-exp}. 
In this subsection,
we prove the following:
\begin{lemm}\label{l.circuit-expansive} Let $\cU$ be a  generating filtrating
Markov partition and $K\subset \cU$ be a circuit of points such that all of its periodic points have large
stable manifolds.  
For every sufficiently large $m, n>0$, 
the Markov partition $\cU_{(m, n)}(K)$ 
is $\cU$-expansive for chains. 
\end{lemm}

Given a filtrating Markov partition $\cU$ and a $f$-invariant set $K\subset \cU$ we say that $K$ is \emph{$\cU$-expansive} if given any $x\neq y$ in $K$ there is $n\in\ZZ$ such that $f^n(x)$ and $f^n(y)$ belongs to different components of $\cU$. Note that this 
is different from $\cU$-expansiveness 
for chains we discussed in Section~\ref{ss.exlo}.

Let us give some preparatory results for the proof of Lemma~\ref{l.circuit-expansive}.
\begin{lemm} Let $\cU$ be a filtrating Markov partition and $K\subset \cU$ be a circuit of points whose periodic orbits have large stable manifolds. 
Then  every point in $K$ has a large stable manifold (that is, the center stable 
manifold which contains the point is contained in the stable manifold of the point).
\end{lemm}
\begin{proof} The property of having a large stable manifold is invariant by negative iteration (if $x$ has a large stable manifold, $f^{-1}(x)$ has a large stable manifold, too). 
If a periodic point has a large stable manifold, so does every point in its orbit, and also every point 
in its stable manifold. In a circuit every point belongs to the stable manifold of a periodic point which is assumed to have a large stable manifold. 
Thus we have the desired property.
\end{proof}

\begin{lemm}Let $\cU$ be a filtrating Markov partition and $K\subset \cU$ be a circuit of points whose periodic orbits have large stable manifolds. Then any neighborhood of $K$ contains a hyperbolic basic set $\La$ such that
\begin{itemize}
 \item $K\subset \La$, and
 \item every point in $\La$ has a large stable manifold in $\cU$. 
\end{itemize}
\end{lemm}
\begin{proof} Any non-trivial circuit of points 
is contained in a hyperbolic basic set $\La_0$ and $\La_0$ can be chosen 
arbitrarily close to $K$. Recall that in $\La_0$ the local stable manifolds vary continuously 
with respect to the point. 
As every point of the circuit has a large stable manifold in $\cU$, every point in $\La_0$ close enough to $K$ has a large stable manifold in $\cU$, too. Thus, by choosing $\La_0$
very close to $K$ we obtain the conclusion.
\end{proof}

\begin{lemm}\label{l.basic-expansive} Let $\cU$ be 
a generating filtrating Markov partition and 
$\Lambda_0 \subset \cU$ be a hyperbolic basic set 
such that every point of $\Lambda_0$ has a large stable manifold. 
Then $\Lambda_0$ is $\cU$-expansive. 
\end{lemm}
\begin{proof}
Suppose that $x,y\in \Lambda_0$ belong to the same rectangle $U_0$ of $\cU$ and 
$f(x),f(y)$ belong to the same rectangle $U_1$. 
Then $x, y$ belong to $U_0\cap f^{-1}(U_1)$ which 
is connected as $\cU$ is generating. 
Arguing by induction, one gets that, if $f^i(x),f^i(y)$ belong to the same rectangles for $i\in\{-m,\dots,n\}$ then $x$ and $y$ belong to the same connected component of 
$\cU_{(m, n)}$. 

As the points in $\Lambda_0$ have large stable manifolds in $\cU$, 
the diameter of the components of $\cU_{(m, n)}(\Lambda_0)$ tends to $0$. 
Thus letting $m, n \to \infty$ we obtain the conclusion.
\end{proof}

Now let us complete the proof of Lemma~\ref{l.circuit-expansive}.
\begin{proof}[Proof of Lemma~\ref{l.circuit-expansive}]
 
Note that $K$ is a hyperbolic set and $\cU_{(m,n)}(K)$ is very close to $K$
by Lemma~\ref{l.mrk-sml} when $m$, $n$ are 
sufficiently large.
Hence, for $m, n$ large enough, the maximal invariant set in $\cU_{(m,n)}(K)$ is 
a hyperbolic basic set $\La_{(m,n)}(K)$ whose points have large stable manifolds in $\cU$. 

By an argument similar to the proof of Lemma~\ref{l.basic-expansive},
given any $\cU_{(m,n)}(K)$-pseudo orbit of points 
$(x_i)_{i\in\ZZ}$, 
there is unique point $y \in \La_{(m,n)}(K)$ such that $f^i(y)$ and $x_i$ belong to the same rectangle of 
$\cU_{(m,n)}(K)$. More precisely, for 
$(x_i)$ consider the rectangle in 
$\cU_{(m+m',n+n')}(K)$ 
which has the same itinerary as $(x_i)$
and let $m', n' \to \infty$: Then we obtain a point 
in a locally maximal invariant set of $\cU_{(m,n)}(K)$
which has the same itinerary as $(x_i)$.
 We say that the point $y$ $\cU_{(m,n)}(K)$-shadows 
$(x_i)_{i\in\ZZ}$.
 
 Assume that $(x^j_i)_{i\in\ZZ}$, $j=1,2$ are two 
 $\cU_{(m,n)}(K)$-pseudo orbits of points such that $x^1_i$ and $x^2_i$ belong to the same rectangle 
 of $\cU$ for $i\in \ZZ$.  
  Let $y^1$ and $y^2$ be the points of $\La_{(m,n)}(K)$ which 
  $\cU_{(m,n)}(K)$-shadow $(x^1_i)_{i\in\ZZ}$ and $(x^2_i)_{i\in\ZZ}$, respectively. 
  This means that $f^i(y^1)$ and $f^i(y^2)$ belong to the same rectangle of $\cU$ 
 for every $i\in\ZZ$. 
 
 According to Lemma~\ref{l.basic-expansive} this implies $y_1=y_2$.  
 As a consequence, $x^1_i$ and $x^2_i$ belong to the same rectangle of $\cU_{(m,n)}(K)$ 
 for every $i\in \ZZ$. 
\end{proof}

\subsubsection{Expansiveness for 
refinements and sub Markov partitions}
In this subsection, we discuss two auxiliary results 
related to expansiveness.

\begin{lemm}\label{l.expu-refi}
Let $\cU$ be a (possibly non-filtrating) 
Markov partition of a 
diffeomorphism $f$, $\cV$ is a generating filtrating 
Markov partition such that $\cV \subset \mathring{\cU}$ holds. Let 
$S$ be a compact 
invariant subset of $\cV$ and assume 
$\cV(S)$ is $\cU$-expansive. 
Then for every $n >0$
the refinement $\cV_{(0, n)}(S)$ is $\cU$-expansive.
\end{lemm}

\begin{proof}
We only prove that $\cV_{(0, 1)}(S)$
is $\cU$-expansive. The general case can be done 
by induction. 
This is a direct consequence
of generating property of $\cV$: Suppose we have 
bi-infinite chains $(V_i)$, $(W_i)$ of connected components of $\cV_{(0, 1)}(S)$.
For each $V_i$, we denote by $V'_i$ the rectangle 
of $\cV(S)$ which contains $V_i$. 
Similarly, we construct $(W'_i)$.

Assume that there is $j$ such that $V_j \neq W_j$.
If $V_j$ and $W_j$ are contained in different 
components 
of $\cV(S)$, then it means that $(V'_i) \neq (W'_i)$
as a chain. Hence by assumption we know there 
is $k$ such that $V'_k$ and $W'_k$ belong to the 
different rectangles of $\cU$, and the same 
holds for $V_k$ and $W_k$.
If $V_j$ and $W_j$ are contained in the same component 
of $\cV(S)$, as $\cV$ is 
generating, we know that $f^{-1}(V_j)$ and 
$f^{-1}(W_j)$ are contained in a different component
of $\cV$. Thus we have the conclusion 
by considering $i =j-1$.
\end{proof}

The following is a direct consequence of the 
definition of the $\cU$-expansiveness, 
so we omit the proof.
\begin{lemm}\label{l.expu-macha}
Let $\cU$ be a filtrating Markov partition, $\cU'$ be 
its sub Markov partition and $\cV$ be another
Markov partition contained in $\mathring{\cU'}$
which is $\cU'$-expansive. 
Suppose $\cW$ is a filtrating Markov partition which is 
matching to $\cV$. Then $\cW$ is $\cU'$-expansive 
(thus $\cU$-expansive), too.
\end{lemm}

\subsection{Local genericity of minimal expansive 
aperiodic classes}

Now we are ready 
to identify a property which is type $\fP_{\MEx}$.

\begin{prop}\label{p.C1-expansive} 
Let $M$ be a smooth closed connected 
three dimensional manifold. 
Let us consider the following $C^1$-robust 
property $\cP_1$ on a filtrating set $\cU$: 
\begin{itemize}
\item $\cU$ is a generating, transitive filtrating 
Markov partition.
\item There is a periodic point $p\in\cU$ such that 
the chain recurrence class 
$[p] \subset \cU$ satisfies property $(\ell)$. 
\end{itemize}
Property $\cP_1$ is type $\fP_{\MEx}$:
In other words, if $(f, \cU, p)$ 
satisfies the condition above
then for any $\delta>0$, $N>0$ and 
any $C^1$-neighborhood 
$\rO$ of $f$ there is $g\in \rO$ such that the 
following holds:
\begin{itemize}
 \item $\supp(g, f)$ is strictly contained in $\cU$.
 \item There are disjoint transitive filtrating Markov partitions $\cU_1,\cU_2 \subset \mathring{\cU}$ for $g$ which are both $\cU$-minimal and $\cU$-expansive.
 \item The minimum period of 
 $\cU_i$ is greater than $N$ for $i=1,2$.
 \item  $\max \{\cdiam(\cU_i), \cdiam (g(\cU_i))\}<\delta$ for $i=1,2$. 
 \item There are periodic orbits 
 $\cO(p_i)\subset \cU_i$ such that 
 $[p_i] \subset \cU_i$ satisfies $(\ell)$ for $i=1, 2$
 (hence $(g, \cU_i, p_i)$ satisfies property $\cP_1$).
\end{itemize}
\end{prop}
\begin{proof}The proof consists of several 
steps.

{\bf Step 1. Preparation.}

First, we fix $\varepsilon>0$ such that every $10\varepsilon$-$C^1$-small perturbation of $f$ belongs to $\rO$. By letting $\varepsilon$ small, 
we may assume that the filtrating Markov partition 
$\cU$ is $10\varepsilon$-robust.  
We also choose $D>0$ such that every diffeomorphisms
which is $10\varepsilon$-$C^1$-close to $f$ is 
$D$-Lipschitz continuous.

As the chain recurrence class $[p]$ satisfies 
property $(\ell)$, Theorem~\ref{t.flex} asserts that, 
by performing an arbitrarily small $C^1$-perturbation of $f$ supported in $\cU$, we can produce 
two $\varepsilon$-flexible points $p_1, p_2$ 
(which are not equal to $p$), with large stable manifolds in $\cU$ and such that the orbits of the points $p_i$, $i=1,2$ meet every rectangle of $\cU$ (for the last property 
apply Corollary~\ref{c.homo-ubiq} and let $\varepsilon$ very small).
Note that we may assume that the period of $p_i$
is larger than $N$.

For $i=1, 2$,
as $p_i$ has a large stable manifold, it admits 
homoclinic orbits, and thus there are hyperbolic circuits $K_i$ which consists of a the unique periodic orbit
 $\cO(p_i)$ and a homoclinic orbit of $p_i$.
Now, we apply Theorem~\ref{t.aff} to $K_1$ and $K_2$. Thus up to an arbitrarily small perturbation 
we have that for sufficiently large $m$ and $n$ 
the refinement
$\cU_{(m, n)}(K_i)$ is affine. By 
abuse of notation, we denote the perturbed 
diffeomorphism by $f$ as well. Note that 
we can assume all of the perturbations has support
which are strictly contained in $\cU$. Recall that
in the setting of filtrating Markov partitions, 
the two notions ``coincides outside $\cU$'' and 
``the support is strictly contained in $\cU$''
are synonymous, see Remark~\ref{r.equi}. Thus
in the following 
we use these two phrases interchangeably.

{\bf Step 2: Expulsion.}

Consider the (possibly non-filtrating) transitive Markov partition $\cU_{(m, n)}(K_i)$. 
According to Lemma~\ref{l.circuit-expansive}, 
for $m, n$ large enough the sub 
Markov partition $\cU_{(m, n)}(K_i)$ 
is $\cU$-expansive 
for $\cU_{(m, n)}(K_i)$-chains. 
Lemma~\ref{l.mini-peri} guarantees that 
for large $m, n$ the minimum period of $\cU_{(m, n)}(K_i)$ are 
larger than $N$. Lemma~\ref{l.mrk-sml} 
shows that each rectangle of $\cU_{(m, n)}(K_i)$ 
has diameter less than $\delta/D$. Note that 
by the choice of $D$ we know that for every 
diffeomorphism $f_0$ which is $10\varepsilon$-close 
to $f$ and every rectangles $U$ of $\cU_{(m, n)}(K_i)$,
the diameter of its image under $f_0$
is less than $(\delta/D) \cdot D = \delta $.   
Also, by Lemma~\ref{l.min-sml-nbd} 
we choose a small neighborhood $\cO_i$ of $K_i$
and choose $m, n$ such 
that $\cU_{(m, n)}(K_i) \subset \mathring{\cO_i}$,
which guarantees the $\cU$-minimality of $\cU_{(m, n)}(K_i)$. 
Note that $\cU_{(m, n)}$ is $10\varepsilon$-robust.

Now Theorem~\ref{t.isola} allows us to obtain a $2\varepsilon$-perturbation $g = g_{\nu}$ of the diffeomorphism with circuits $K'_{i}$ similar to $K_i$,
supported in 
$\cU_{(m, n)}(K_{1}) \cup \cU_{(m, n)}(K_{2}) $ 
such that $g$ admits transitive, filtrating Markov partitions 
$\cU_i\subset \cU_{(m, n)}(K_i) = \cU_{(m, n)}(K'_{i})  $, $i=1,2$ matching to 
$\cU_{(m, n+\nu;g)}(K_{i})$ (where $\nu$ is 
some positive integer).

Remark~\ref{r.peri-sub} guarantees that we have
the largeness 
of minimum period for $\cU_i$. The fact that 
the diameter of $\cU_{(m, n)}(K_i)$ is less than 
$\delta /D$ shows that the same holds for $\cU_i$. 
Lemma~\ref{l.min-sml-nbd} guarantees 
the $\cU$-minimality of $\cU_i$. 

Let us confirm the $\cU$-expansiveness of 
$\cU_i$. Recall that 
$\cU_{(m, n)}(K_i) = \cU_{(m, n)}(K'_{i})$ 
is $\cU$-expansive for both $f$ and $g$ 
(the $\cU$-expansiveness for $g$ is a consequence
of the support, see Lemma~\ref{l.exlo}).
By Lemma~\ref{l.expu-refi}, we know that
$\cU_{(m, n+\nu)}(K'_{i})$ is $\cU$-expansive
for $g$. Finally, Lemma~\ref{l.expu-macha} 
concludes that $\cU_i$, which is a matching
Markov partition of $\cU_{(m, n+\nu; g)}(K'_{i})$,
is $\cU$-expansive for $g$. Also, note that the 
matching property implies the generating property of 
$\cU_i$.

{\bf Step 3: Recovery.}

The filtrating Markov partitions $\cU_i$ satisfies
all the desired condition except property $(\ell)$. 
Let us perform the final perturbation to 
recover it. Recall that $\cU_i$ is 
$10\varepsilon - 2\varepsilon = 8\varepsilon$-robust 
(this is a consequence of Theorem~\ref{t.isola}
and the fact 
$\cU_{(m,n)}$ is $10\varepsilon$-robust for $f$)
and contains 
a periodic point $p_i$ which is still 
$\varepsilon$-flexible and having a large
stable manifold in $\cU_i$ (this is a consequence
of Theorem~\ref{t.isola}). 
Furthermore, the homoclinic class of $p_i$
in $\cU_i$ is not trivial, because it contains the 
circuit $S_{i, \nu}$.
Now we apply 
Theorem~\ref{t.relative} to $g$: Then 
we take a diffeomorphism $h$ which is 
$4\varepsilon$-close to $g$ such that
the relative
homoclinic class of $p_{i, h}$ in ${\cU_i}$
satisfies property 
$(\ell_{\cU_i})$.
In particular, the chain recurrence class of $p_{i,h}$
in $\cU_i$
satisfies property $(\ell)$. 

Since the supports of the perturbations for $\cU_1$
and $\cU_2$ are disjoint, we can perform the 
perturbation without any interference.
Note that 
the structure of $\cU_i$ as a filtrating Markov partition 
is preserved, so we can keep the transitivity, 
$\cU$-expansiveness, smallness of the diameters and the 
largeness of the minimum period.
Also, $\cU$-minimality is preserved for the perturbation 
is local, see Lemma~\ref{l.minimal-robust}.
Thus we constructed 
the desired filtrating Markov partitions $\cU_i$ up
to $6\varepsilon$-perturbations whose support 
is contained in $\cU$.
\end{proof}

\begin{rema}\label{r.local-why}
Let us see the usefulness of the concept
of the stability under local perturbations. 
Consider the perturbation 
we made by Theorem~\ref{t.isola} in 
Step 2 of above proof. 
We want to produce a new filtrating Markov partition
keeping the minimality. 
We know that 
$\cU_{(m, n)}(K_i)$ is $\cU$-minimal and 
this property is $C^1$-robust. 
However, the size of the perturbation 
by Theorem~\ref{t.isola} is $2\varepsilon$ and
we are not sure if the $C^0$-robustness for the 
minimality is greater than that. We can shrink 
$\varepsilon$ as small as we want, but it would 
change the $C^1$-robustness of the minimality.

This brings a non-trivial problem, but if we know 
that the $\cU$-minimality is stable under the 
local perturbations, we can circumvent it.
\end{rema}

Note that if $(\cU, p)$ satisfies property $(\ell)$,
then by taking refinements we obtain 
the assumption of Proposition~\ref{p.C1-expansive}.
Thus
Proposition~\ref{p.C1-expansive}, together with 
Proposition~\ref{p.principle-expansive} implies item~\ref{i.2} of Theorem~\ref{t.minimal-non}: 
\begin{theo}\label{t.miniex}
 Let $(\cO,\cR, p_f)$ be a $C^1$-open set of diffeomorphisms on a closed $3$-manifold $M$ admitting a transitive filtrating Markov partition 
 $\cR$ and a periodic point  $p_f$ 
 varying continuously with $f$ such that $(f,\cR,[p_f])$ satisfies the property $(\ell)$. 

Then there is a residual subset $\cG\subset \cO$ 
such that every $f\in\cG$ has
an uncountable set of chain recurrence classes which are all minimal and expansive.
\end{theo}

\subsection{Constructing an example}
Theorem~\ref{t.miniex} gives an
example of $C^1$-locally generic existence of $C^1$-diffeomorphisms having minimal, expansive 
aperiodic classes. This is a consequence of Baire's 
category theorem and because of that 
there is a lack of concrete information 
about diffeomorphisms we obtained.

On the other hand, it is possible to construct a concrete 
example of such diffeomorphism by applying 
Proposition~\ref{p.C1-expansive} successively.
Let us see this. First, by applying 
Proposition~\ref{p.C1-expansive} repeatedly 
one can prove the following (we omit the proof 
for it is straightforward):
\begin{prop}\label{p.exa}
Let $M$ be a smooth closed connected 
three dimensional manifold. 
Let us consider $f \in \mathrm{Diff}^1(M)$
which satisfies the assumption of 
Proposition~\ref{p.C1-expansive}. 
Then, one can construct a Cauchy sequence (in the 
$C^1$-distance)
$(f_n)_{n \geq 1}$ and a nested sequence of 
transitive, generating 
filtrating Markov partitions
 $(\cU_n)$ satisfying the following:
\begin{itemize}
\item $f_1 =f$.
\item $\{(f_n, \cU_n)\}_{n\geq 1}$ is a nested sequence
for an aperiodic class.
\item $\cU_{n+1}$ is $\cU_{n}$-expansive for chains
with respect to $f_{n+1}$.
\item $\cU_{n+1}$ is $\cU_{n}$-minimal for $f_{n+1}$.
\end{itemize}
\end{prop}
Now, take a $C^1$-diffeomorphism 
$f_{\infty} = \lim_{n\to\infty}f_n$. By  Proposition~\ref{p.minimale} and Proposition~\ref{l.expansive},
we know that $\Lambda:=\cap_{n \geq 1} \cU_n$
is an expansive, minimal aperiodic class. 
This construction is more concrete than the result of 
Theorem~\ref{t.miniex} since we have some information
on $f_{\infty}$, for instance it coincides with $f$ 
outside $\cU_1$. This method would be 
useful, if one wishes to construct some peculiar 
example keeping control for some of the part of 
the diffeomorphism.
Note that we can perform such a concrete construction
since we have Proposition~\ref{p.minimale} and Proposition~\ref{l.expansive},
which are stated in terms of sequence of maps $(f_n)$.

We can obtain the similar results in the later 
examples, but for the sake of simplicity 
we will not state them explicitly.

%
%
%
%

\subsection{Ergodic diameters of sub Markov partitions}
\label{ss.C1infi}

In this subsection, we discuss the size of the  space of 
invariant measures supported on a special kind of 
sub Markov partition. We begin with a definition.

\begin{defi}
Given a (possibly non-filtrating) 
Markov partition $\cV$, its \emph{ergodic diameter}
$\delta_{\mathrm{erg}}(\cV)$ is the diameter of the set of probability measures obtained as accumulation points of the $\cV$-pseudo orbits of points
$(x_i) \subset \cV$, that is,
the diameter of the closed set 
$\mathcal{P}_{\infty, \mathrm{pseudo}}(\cV)$,
see Section~\ref{ss.unique}. 
\end{defi}

\begin{lemm}\label{l.ergodiam} 
Consider a circuit of points $K$ of a diffeomorphism $f$ consisting of a unique periodic orbit $\gamma$ and a finite set of homoclinic orbits of $\gamma$. Let $\cU$ be a filtrating Markov partition containing $K$
such that $\gamma$ has a large stable manifold in $\cU$. 

Then for every $n, m \geq 0$, $\cU_{(m, n)}(K)$ is a
transitive Markov partition. Moreover, for any $\eta>0$  there is $m_0, n_0>0$ such 
that for any $m\geq m_0, n \geq n_0$,
the ergodic diameter of $\cU_{(m, n)}(K)$ 
is smaller than $\eta$: 
\[
\delta_{\mathrm{erg}}(\cU_{(m, n)}(K))<\eta.
\]

\end{lemm}
\begin{proof} $\cU_{(m, n)}(K)$ is a 
sub Markov partition of $\cU_{(m, n)}$. 
Let us see the transitivity of it. 
Every rectangle contains 
at least a point of $K$ which is chain recurrent. Therefore, for any $\varepsilon>0$ 
one can find an $\varepsilon$-pseudo orbit in $K$ 
connecting any pair of connected components of 
$\cU_{(m, n)}(K)$.
This defines a $\cU_{(m, n)}(K)$-chain of 
connected components between any pair 
of connected components.

Now, let us discuss the ergodic diameter.
Recall that $\cU_{(m, n)}(K)$ converges to $K$
as $m, n \to \infty$, see Lemma~\ref{l.mrk-sml}.
Now, for $m, n$ large enough, 
$\cU_{(m, n)}(K)$-pseudo orbits are 
$\varepsilon$-close to some 
$\varepsilon$-pseudo orbits of points in $K$, 
for $\varepsilon$ arbitrarily small. 
Note that a pseudo-orbit in $K$ has two types 
of orbit segments: One which follows $\gamma$
or the one which follows the homoclinic orbits. 
When $m, n$ are large and 
$\varepsilon$ is small, 
the former one tends to $\gamma$. Also, 
for the latter one the parts of the orbit
which follows $\gamma$ will have 
an arbitrarily large portion. 
Thus, when $m, n$ are large then 
for every sufficiently large $k$ 
the measures 
$\cP_{k, \mathrm{pseudo}}(\cU_{(m,n)}(K))$
are all very close to the Dirac measures supported on 
$\gamma$.

Therefore the ergodic diameter of $\cU_{(m,n)}(K)$ is arbitrarily small for every sufficiently large 
$m$ and $n$.  
\end{proof}

\begin{rema}
If $\cU$, $\cV$ are possibly non-filtrating Markov partitions such that 
$\cV$ is matching to $\cU$,
then for every $n$ we have 
\[
\cP_{k, \mathrm{pseudo}}(\cV) \subset \cP_{k, \mathrm{pseudo}}(\cU).
\]
If $\cV$ contains a periodic orbit, then we know 
that $\cP_{k, \mathrm{pseudo}}(\cV)\neq \emptyset$.
Thus we have
\[
\delta_{\mathrm{erg}}(\cV) \leq \delta_{\mathrm{erg}}(\cU).
\]
\end{rema}

\subsection{Local genericity of 
minimal aperiodic classes supporting 
infinitely many ergodic measures}

Recall that we have fixed 
$\lambda_0>0$ satisfying $\prod_1^{+\infty} (1+\lambda_0^i) \cdot \sum_1^{+\infty}\lambda_0^i\leq \frac12,$ see Section~\ref{ss.infierg}.
Also, recall that $(\ell_\cV)$ is a version of property $(\ell)$ 
for a relative homoclinic class in $\cV$, see Section~\ref{s.preliminar}.

The aim of this section is to prove the following.

\begin{prop}\label{p.C1minnergo}

Let us consider the following family of $C^1$-robust
properties $(\cP^n)_{n \geq 1}$ for 
a $C^1$-diffeomorphism
$f$ having a transitive filtrating Markov partition 
$\cU$ containing periodic orbits 
$\Gamma_n :=\{\gamma_i\}_{i=1,\ldots,n}$:
There is $(m_0, n_0)$ such that $\cU_{(m_0, n_0)}$
has $n$ mutually distinct sub Markov 
partitions $(\cW_i)_{i=1,\ldots,n}$ such that 
\begin{itemize}
 \item[(A)] 
 $\gamma_i \subset \cW_i$ and $(f, \cW_i , p_i)$ satisfies property $(\ell_{\cW_i})$, where $p_i$
 is a point of $\gamma_i$.
 \item[(B)] $\delta_{\mathrm{erg}}(\cW_i)<\lambda_0^{n+1} \rho(\cM_n)$ where $\cM_n=\{\mu_1,\dots,\mu_n\}$ is the set of Dirac measures supported on the orbits of $\Gamma_n$ 
and $\rho(\cdot)$ denotes the independence radius,
see Section~\ref{ss.infierg}. 
 \item[(C)] For any connected component $C$ of $\cU$ one has $C\cap\cW_i\neq\emptyset$. 
\end{itemize}
This family of properties $(\cP^n)$
is type $\fP_{M, \infty}$
(see Defintion~\ref{d.mini-infi-prop}), that is,
given any $C^1$-neighborhood $\rO$ of $f$ and
for any $\delta>0$ and $N>0$ there is 
$g\in \rO$ such that the following holds:
\begin{enumerate}
\item The support $\supp(g,f)$ is strictly 
contained in $\cU$.
 \item There are disjoint, transitive filtrating Markov 
 partitions $\cU_1,\cU_2 \subset \mathring{\cU}$ 
 which are both $\cU$-minimal for chains (with respect to $g$).
 \item For $g$, the continuations 
 $\{\gamma^g_i\}$ are defined and
 $\{\gamma^g_i\}$ are outside $\cU_1 \cup \cU_2$.
 \item $\cdiam(\cU_j)$, $\cdiam(g(\cU_j)) < \delta$ and the minimum period
of $\cU_j$ is larger than $N$ for $j=1, 2$.
 \item There are periodic orbits $\gamma_i^j\subset \cU_j$,$j\in\{1, 2\}$,  $i\in\{1,\dots, n+1\}$.
 \item For any $j\in\{1,2\}$ and $i\in\{1,\dots, n\}$ one has  
\[\mathfrak{d}(\mu^g_i,\mu^{j}_i)< \lambda_0^{n+1} \rho (\{\mu_i^n \mid i=1,\ldots, n\})\]
where $\mu^g_i,\mu_i^j$ are the Dirac probabilities associated to $\gamma^g_i$ and $\gamma_i^j$, respectively. 
 \item For each $j=1,2$, there is a pair of 
 integers $(m_j, n_j)$ such that 
 in the refinements $\cU_{j, (m_j, n_j)}$  
 there are disjoint transitive sub Markov partitions 
 $\cW^j_i$ of $\cU_{j, (m_j, n_j)}$,  
 ($i=1,\dots,n+1$) which 
 satisfies the following: 
\begin{itemize}
 \item There are periodic orbits
 $O(p_i^{j})\subset \cW^j_i$ ($i=1,\ldots, n+1$) such that
 $(g, \cW^j_i, p_i^{j})$ satisfies 
 the property $(\ell_{\cW_i^j})$.
 \item $\delta_{\mathrm{erg}}(\cW^j_i)<\lambda_0^{n+2} \rho(\cM^j_{n+1})$ where $\cM^j_{n+1}=\{\mu^j_1,\dots,\mu^j_{n+1}\}$ is the set of Dirac measures supported on the orbits of $p_i^{j}$.
 \item For any connected  
 component $C$ of $\cU_j$ one has 
 $C\cap \cW_{i}^j \neq\emptyset$
 for every $j=1,2$ and $i=1,\ldots,n+1$.
\end{itemize}
 \end{enumerate}
Note that condition (7) guarantees that 
$(g, \cU_j, \{\gamma^j_{i}\}_{1\leq i \leq n+1})$ 
satisfies property $\cP^{n+1}$. 
\end{prop}

\begin{proof} We choose $\varepsilon >0$ such 
that any $10\varepsilon$-perturbation of $f$ belongs to $\rO$
and the filtrating Markov partition $\cU$ is $10\varepsilon$-robust. 

In the course of the proof, the confirmations of
conditions (1, 2, 4) are easy or similar in the proof of 
Proposition~\ref{p.C1-expansive}.
Thus we avoid the detailed explanation of them and
concentrate on how to obtain 
conditions (3, 5, 6, 7).

{\bf Step 1: Preparation.}

As $f$ satisfies $(\ell_{\cW_i})$ in $\cW_i$, according to Theorem~\ref{t.flex} there is an arbitrarily $C^1$-small perturbation $f_0$
of $f$ supported in $\bigcup \cW_i$ such 
that there are $\varepsilon$-flexible points whose
orbits are $\varepsilon$-dense in the relative homoclinic class of $p^{f_0}_i$ in $\cW_i$, $i=1,\dots, n$. 
We choose $f_0$ such that 
assumptions (A), (B) and (C) still hold
(note that a local perturbation on a filtrating Markov partition 
does not change the ergodic diameter, see Remark~\ref{r.robust}).

This allows us to choose a family $p^j_i$,
$j=1,2$, $i=1,\dots,n+1$ with the following properties. 
 \begin{itemize}
 \item All the points $p^j_i$, $j=1,2$, $i=1,\dots, n$ 
 are $\varepsilon$-flexible and 
 have large stable manifolds in $\cW_i$ (thus in $\cU$). 
\item $p^j_{n+1}$ is $\varepsilon$-flexible and 
has a large stable manifold in $\cU$ (for that 
we only need to choose one from one of $\cW_i$). 
 \item As the consequence of assumption 
 (C), $p^j_{i}$ are pairwise 
 homoclinically related in $\cU$. 
Also, each $p^j_i$ has homoclinic point to itself.
  \item For $i= 1,\dots, n$, the orbits of $p^j_i$, 
  are contained in $\cW_i$ and meet every connected component of $\cU$. 
 \end{itemize}
Thus we chose two circuits $K_j$ $(j=1,2)$ 
consisting of the following objects:
\begin{itemize}
\item Periodic orbits $\{\cO(p^j_i)\}$, $i=1,\ldots, n+1$.
\item Heteroclinic orbits between 
$p^j_{i}$ and $p^j_{i'}$ for $i \neq i'$,
\item Homoclinic orbits for each $p^j_{i}$.
\end{itemize}

We apply Theorem~\ref{t.aff} to $K_j$ such that 
for sufficiently fine refinement $\cU'$,
the sub Markov partition $\cU'(K_j)$ are affine for
$j=1, 2$.
We denote the perturbed
map by $f_1$.

Then we apply Lemma~\ref{l.min-sml-nbd} for $K_j$: 
We choose a small neighborhood $\cO_j$ of $K_j$
satisfying the conclusion for $j=1, 2$.
Then for every 
sufficiently large $m$ and $n$ we have 
$\cU'_{(m, n;f_1)}(K_j)$ is a subset of $\cO_j$.
Note that we may assume every connected component 
of them has a small diameter by letting $m, n$ large. 
In particular, we may assume 
that $\cU'_{(m, n)}(K_j)$ are affine, 
does not 
contain $\{p^{f_0}_i\}$ (points in $\Gamma^{f_0}_n$) and 
$\cU'_{(m, n)}(K_1)$, $\cU'_{(m, n)}(K_2)$ 
are disjoint.

{\bf Step 2. Expulsion of circuits.}

According to Theorem~\ref{t.isola}, there is $g$
which is $2\varepsilon$-$C^1$-close to $f_1$ 
such that 
the support of $g$ is strictly contained in
$\cU'_{(m, n)}(K_1) \cup \, \cU'_{(m, n)}(K_2)$,  
there are transitive filtrating Markov partitions 
$\cU_j \subset \cU'_{(m, n+\nu)}(K_j)$ 
(where $\nu$ is some positive integer)
and circuits $K'_{j}\subset \cU_j$ ($j=1,2$) with the following properties:
\begin{itemize}
 \item The periodic $g$-orbits of the circuit $K'_{j}$ are the $f_0$-periodic 
 orbits $p^j_i$, $j=1,2$, $i=1,\dots, n+1$.
 \item $p^j_i$ ($j=1,2$, $i=1,\dots, n+1$) 
are still $\varepsilon$-flexible and have large stable manifolds in $\cU_j$. 
\item $\cU_j \subset \cU_{(m, n+\nu)}(K'_j)$ are matching 
and they are $8\varepsilon$-robust. 
\end{itemize}
Notice that by the choice of $\cU'_{(m, n)}(K'_j)$, 
the filtrating Markov partitions $\cU_1$ and $\cU_2$ are both $\cU$-minimal.

As the periodic orbits $p^j_i$ for $f_1$ and $g$ 
coincide and for $i=1\dots,n$ these orbits are contained in $\cW_{i}$, we have 
\[
\mathfrak{d}(\mu^{g}_i,\mu^j_i) < 
\lambda_0^{n+1}\rho(\{\mu^{g}_i \mid i=1,\ldots, n \}) \mbox{ for } i=1,\dots,n \mbox{ and } j=1,2
\]
where $\mu^j_i$ is the Dirac measure 
supported on the orbit of $p^j_i$.

{\bf Step 3: Recovery.}  

It remains to recover properties in (7).
For each $p^j_i$, since $\cU_j$ is transitive, we can
find a homoclinic orbit of it such that it visits every 
rectangle of $\cU_j$, see Lemma~\ref{l.homoclinic}.
Consider a circuit $S_{i, j}$
which consists of the periodic orbit 
$\cO(p^j_i)$ and such a homoclinic orbit.
Then consider $\cU_{j,(m, n)}(S_{i, j})$. Since 
$p^j_i$ has a large stable manifold in $\cU_j$, by
Lemma~\ref{l.ergodiam} for every sufficiently 
large $m, n$ we have
$$\delta_{\mathrm{erg}}(\cU_{j,(m, n)}(S_{i, j}))<\lambda_0^{n+2} \rho(\cM^j_{n+1}).$$ 

Furthermore, by definition of $S_{i, j}$ for any connected 
component $C$ of $\cU_j$ one has 
$$\cU_{j,(m, n)}(S_{i, j})\cap C\neq \emptyset.$$

We fix sufficiently large $m, n$ and put
$\cW^j_i:= \cU_{j, (m,n)}(S_{i, j})$.
Now Theorem~\ref{t.relative} asserts that 
there is 
a $4\varepsilon$-perturbation $h$ of $g$ 
with support strictly contained in 
$\cW^j_i$ having property $(\ell_{\cW^j_i})$.
Note that they are disjoint if $m, n$ are sufficiently large,
considering the fact that $p^j_i$ has large stable manifold. 
Also, this perturbation keeps the structure of the Markov partitions $\cU_j$
and $\{\cW^j_i\}$, and does not destroy
the properties we obtained.
Thus we have the conclusion.
\end{proof}

To obtain 
item~\ref{i.1} of Theorem~\ref{t.minimal-non},
we are left to show that property $(\ell)$
implies 
$\cP^1$. This can be observed as follows:
\begin{itemize}
\item If $(\cU, p)$ satisfy property $(\ell)$, then 
by adding arbitrarily small perturbation we may 
assume that there is a hyperbolic periodic point  
$p_0$ homoclinically related to $p$, having large
stable manifold, $\varepsilon$-flexible for 
very small $\varepsilon$ 
and visits every rectangle of $\cU$.
\item Consider a circuit $K_0$ consisting of $p_0$ and 
a homoclinic point. 
By Theorem~\ref{t.relative} we may 
assume that $(p_0, \cU_{(m, n)}(K_0))$ satisfy
property $\ell_{\cU_{(m, n)}(K_0)}$ up to some 
small perturbation.
\end{itemize}
Thus, by
Proposition~\ref{p.C1minnergo} and 
Proposition~\ref{p.principle-genurgo} 
we have the following:
\begin{theo} Let $\rO$ 
be a $C^1$-open set of $\mathrm{Diff}^1(M)$
(where $M$ is  
a closed $3$-manifold) such that 
every $f \in \rO$ admits 
a transitive filtrating Markov partition having a hyperbolic periodic point $p_f$ varying continuously with $f$ over $\rO$ and 
$(f, \cU, p_f)$ satisfies the property $(\ell)$. 

Then there is a residual subset $\rR\subset \rO$ such that every $f\in\rR$ has
an uncountable set of aperiodic chain recurrence classes 
which are all minimal but are supporting infinitely many ergodic measures.
\end{theo}

\subsection{$C^1$-genericity of transitive non-minimal aperiodic classes}\label{s.C1nonmin}

In this subsection, we give an example 
of a property of dynamical systems 
which is type $\fP_{\mathrm{TnM}, A, B}$.
We begin with a result about filtrating Markov partitions.

\begin{lemm}\label{l.transi2} 
Let $\cU$ be a filtrating Markov partition 
and $K\subset \mathring\cU$ be a circuit such that every
periodic point has a large stable manifold. 
Assume that there is a point $p\in K$ such that 
for every connected component $C$ of $\cU$ there 
are $i<0<j$ such that $f^i(p), f^j(p) \in \mathring C$.  
Then, for every sufficiently large $m$ and $n$, 
$(\cU_{(m, n)}, \cU_{(m, n)}(p))$ 
is stably $\cU$-transitive for local perturbations.
\end{lemm}

\begin{proof} 
Let us confirm the 
second condition of Definition~\ref{d.trachai}
and the assumption of Lemma~\ref{l.stab-tran},
see Remark~\ref{r.statran}.

Since $\cU$ is transitive, we deduce that 
every rectangle of $\cU$ contains at least one rectangle 
of $\cU_{(m, n)}$ and it is in its interior 
if $m, n \geq 1$. This implies the second condition of 
Definition~\ref{d.trachai}.
Now, let us confirm 
the assumption of Lemma~\ref{l.stab-tran}.
Consider a connected component $A$ of $\cU$.
By definition, there is $i^+(A)$ such that 
$f^{i^+(A)}(p) \in \mathring{A}$. 
We denote the distance between
$\{f^{i^+(A)}(p)\}$ and 
$M \setminus \mathring{A}$ by $\delta$,
where $M$ is the ambient manifold.
For each $f^j(p)$, $j=0,\ldots, i^+(A)$,
we choose a sufficiently small compact neighborhood 
$O_{j}$ such that the distance between
$f^{i^+(A) -j}(O_{j})$
and $M \setminus \mathring{A}$ is larger than 
$\delta /2$.


The assumption that every periodic orbit of  $K$ has
a large stable manifold implies that as $m, n \to \infty$
the Markov partition $\cU_{(m, n)}(K)$ 
satisfies $\cdiam\, \cU_{(m, n)}(K) \to 0$,
see Lemma~\ref{l.mrk-sml}. 
Thus, when $m, n$ are large then 
for every $\cU_{(m, n)}(K)$-chain of 
connected components $(U_i)_{i=0,\ldots, i_0}$
where $i_0 \leq i^+(A)$, 
we have $U_{i_0} \subset O_{i_0}$,
since when each connected component is small 
then the chain of connected components is almost equal to the 
true orbit.
Accordingly we have $f^{i^+(A) - i_0}(U_{i_0}) \subset \mathring{A}$.
Thus we can obtain the sufficient condition
for the stable $\cU$-transitivity in Lemma~\ref{l.stab-tran} for 
a fixed connected component $A$ by taking sufficiently 
fine refinement.

%

The same argument holds for backward iterations.
Since the number of the 
connected components of $\cU$ is finite, by choosing 
$m, n$ for each connected component and 
letting $m, n$ larger than 
all of them, we obtain the desired refinement.
\end{proof}

\begin{lemm}\label{l.trans-macha}
Let $\cU$, $\cV$ be a filtrating Markov partition
of $f \in \mathrm{Diff}^1(M)$ 
such that $\cV \subset \mathring{\cU}$ and 
$(\cV, V_0)$ is stably $\cU$-transitive. 
If $\cW$ is a matching Markov partition of $\cV$,
then $(\cW, W_0)$ is stably $\cU$-transitive, too,
where $W_0$ is the unique rectangle of $\cW$
contained in $V_0$.
\end{lemm}
\begin{proof}
Let $A$ be any connected component of $\cU$. 
Then there is a connected component $C$ of 
$\cV$ such that $C \subset \mathring{A}$. 
Since $\cW$ is matching to $\cV$, there is 
a connected component $C'$ of $\cW$ contained 
in $C$. Thus, $C' \subset \mathring{A}$ and
this shows the second condition of Definition~\ref{d.trachai}.

Let us check the first condition for $g$ whose 
support is contained in $\cW$. 
Since $\cV$ is stably $\cU$-transitive, 
there is $i^+(A)$ such that 
$g^{i^+(A)}(V_{0}) \subset \mathring{A}$, 
which implies $g^{i^+(A)}(W_{0}) \subset \mathring{A}$.  The confirmation 
for the backward iteration is similar.
\end{proof}

Let us prove the main result.

\begin{prop}\label{p.C1-non-minimal} 
Let $f$ be a diffeomorphism of a closed $3$ manifold $M$.  Let $A,B\subset M$ be two disjoint compact subsets. 

Consider the following $C^1$-robust 
property $\cP_3$ for a
filtrating set $\cU$:
\begin{itemize}
\item $\cU$ is a transitive filtrating Markov partition whose rectangles are either disjoint from $A$ (resp. $B$) or included in $\mathring{A}$ (resp. $\mathring{B}$).
Thus we can define the restriction of 
Markov partition $\cU|_{A}$ (resp. $\cU|_B$),
which consists of rectangles contained in $A$ (resp. $B$), see Definition~\ref{d.tra-nonmi}.  
\item There is a periodic point $p\in\cU$ 
such that $(f, \cU, p)$ satisfies property $(\ell)$. 
\item The restriction of $\cU$ to $A$
(resp. $B$) contains a transitive sub Markov partition
and it satisfies property $(\ell_{\cU|_A})$ 
$(resp. (\ell_{\cU|_B}))$. 
\end{itemize}

This property is
 type $\fP_{\mathrm{TnM}A,B}$, 
 that is, if $f$ satisfies these conditions for 
$(\cU, p)$, 
then, for any $\delta>0$, $N>0$ and any $C^1$-neighborhood $\rO$ of $f$ there is $g\in \rO$ 
such that:
\begin{itemize}
 \item $\supp(g, f)$ is strictly contained in $\cU$.
 \item There are disjoint transitive filtrating Markov 
 partitions $\cU_1$, $\cU_2$ whose minimum periods 
 are larger than $N$.
 \item $\cU_1,\cU_2 \subset \mathring{\cU}$
 are both stably $\cU$-transitive for 
 local perturbations with respect to $g$.
 \item $\cdiam(\cU_i), \cdiam g(\cU_i)<\delta$ for $i=1,2$.
 \item $\cU_1$, $\cU_2$ satisfy the property $\cP_3$,
 that is:
 \begin{itemize}
\item The restrictions of $\cU_i$ ($i=1,2$) 
to $A$ (resp. $B$) is well defined.
\item $\cU_i$ satisfies property $(\ell)$.
\item  $\ell_{\cU_{i}|_A}$ 
(resp. $\ell_{\cU_{i}|_B}$) 
contains a transitive sub Markov partition
 satisfying property 
 $(\ell_{\cU_{i}|_A})$ (resp. $(\ell_{\cU_{i}|_B})$).  
 \end{itemize}
\end{itemize}
\end{prop}
\begin{proof}
Similar to the proofs of 
Proposition~\ref{p.C1-expansive} and
Proposition~\ref{p.C1minnergo}, we divide 
the proof into three steps. Again the arguments for 
obtaining the smallness of the diameter of new
Markov partitions and the largeness of the 
minimum periods are almost the same, so 
we keep the explanation of them short.

{\bf Step 1: Preparation.}

We fix $\varepsilon >0$ such that every diffeomorphism
which is $10\varepsilon$-$C^1$-close to $f$ is 
contained in $\rO$. Also, we assume that
the filtrating Markov partition $\cU$ is
$10\varepsilon$-robust.

First, Theorem~\ref{t.flex} allows us to perform an arbitrarily small perturbation such that there are six
$\varepsilon$-flexible points $p^i_0$, $p^i_A$ and $p^i_B$ $(i=1, 2)$ with large stable manifolds in $\cU$
satisfying the following:
\begin{itemize}
 \item The orbit of $p^i_0$ meets every rectangles of
 $\cU$. 
 \item The orbit of $p^i_A$  
 is contained in $A$. There is a homoclinic orbit 
 $x^i_A$ of $p^i_A$ whose orbit is contained in $A$.
 \item The same condition holds for $p^i_B$,
 $B$ and a homoclinic orbit $x^i_B$.
\end{itemize}
Let us see how to take such points. To take $p^i_A$,
we apply Theorem~\ref{t.flex} to the transitive Markov
partition satisfying property $(\ell_{\cU|A})$ in the 
hypothesis. Then using the largeness of the stable 
manifold of $p^i_A$ we can choose a homoclinic 
point $x^i_A$. 
The construction of $p^i_B$ and $x^i_B$ can be 
done similarly. For the choice of $p^i_0$ apply
Theorem~\ref{t.flex} to the
transitive Markov partition $\cU$.
By abuse of notation we denote the perturbed map by $f$. Note that we may assume that these six periodic 
points have large periods.

Note that the largeness of the stable manifolds of
$p^i_0$ and $p^i_A$, $p^i_B$ implies that
there are heteroclinic orbits connecting among them.
We choose such heteroclinic orbits and denote them by 
$y^i_{0A}$, $y^i_{A0}$,$y^i_{0B}$ and $y^i_{B0}$ (where $y^i_{0A}$ is 
a heteroclinic orbit from $p_{0}^i$ to $p_{A}^i$, etc...).
Then we take two
disjoint circuits $K_i$ ($i =1, 2$) 
such that:
\begin{itemize} 
\item Periodic points of $K_i$ are $p_0^{i}$, $p^i_A$ and $p^i_B$.
\item Homoclinic orbits are $x^i_A$ and $x^i_B$.
\item Heteroclinic orbits are $y^i_{0A}$, $y^i_{A0}$
and $y^i_{0B}$ and $y^i_{B0}$.
\end{itemize}
Note that we have sub circuits $K^i_A$ 
(resp. $K^i_B$) of $K_i$ which consists of a periodic point 
$p^i_A$ (resp. $p^i_B$) and a homoclinic orbit
$x^i_A$ (resp. $x^i_B$) which are contained in 
$A$ (resp. $B$).

For each $K_i$, we apply Theorem~\ref{t.aff}.
Then by an arbitrarily small perturbation we may assume 
that for every large $m$ and $n$ 
the Markov partition $\cU_{(m, n)}(K_i)$ is affine.

Note that 
$p^i_0 \in \cU_{(m, n)}(K_i)$ 
for every $m$ and $n$. Thus 
for every sufficiently large $m$ and $n$ we can apply Lemma~\ref{l.transi2}. We fix such large $m$ and $n$.
We also assume that $\cdiam(\cU_{(m, n)}(K_i))$
is very small.

{\bf Step 2: Expulsion.} 

Now Theorem~\ref{t.isola} allows us to take a $4\varepsilon$-perturbation $g$ of the diffeomorphism supported in $\cU_{(m, n)}(K_1) \cup \,\cU_{(m, n)}(K_2)$ such that $g$ admits transitive filtrating Markov partitions $\cU_i\subset \cU_{(m, n+\nu)}(K_i)$ 
(where $\nu$ is some positive integer) 
with the following properties: 
\begin{itemize}
\item $\cU_i$ is matching 
to $\cU_{(m, n+\nu)}(K_i)$.
\item There are circuits $K'_{i}\subset \cU_i$ which
are similar to $K_i$ for $i=1, 2$.
\item The periodic orbits of $K'_{i}$
 coincides with the those of $K_i$.
 \item Every periodic point of $K'_{i}$ 
 is $\varepsilon$-flexible 
 and has a large stable manifold in $\cU_i$.
\item $\cU_i$ is stably $\cU$-transitive 
for local perturbations.
\end{itemize}

{\bf Step 3: Recovery.}

Note that for $\cU_i$ we can define restrictions
to $A$ and $B$, for it is a subset of $\cU_{(m, n)}$.
We apply Theorem~\ref{t.relative} to 
$\cU_i|_A$ and $\cU_i|_B$ (note that these
perturbations do not interfere, for the supports are disjoint):
As $\cU_i\cap A$ and $\cU_i\cap B$ contain sub circuit containing a periodic point and a homoclinic point, 
one deduces that $\cU_i\cap A$ and $\cU_i\cap B$ contain a non-trivial transitive sub Markov partition
satisfying property $(\ell_{\cU_i|_A})$
and $(\ell_{\cU_i|_B})$ 
up to a $4\varepsilon$-perturbation. 
Note that it implies that 
$\cU_i$ satisfies property $(\ell)$.

By the choice of $m$ and $n$, Lemma~\ref{l.transi2} 
and Lemma~\ref{l.trans-macha}
imply that $\cU_i$ is stably $\cU$-transitive. 
Thus, $\cU_i$ satisfy all the announced properties,
and it completes the proof. 
\end{proof}

Now to conclude the item \ref{i.4} of Theorem~\ref{t.minimal-non} we only need to 
prove the following:
\begin{lemm}
Suppose that a chain recurrence class $[p]$,
where $p$ is a hyperbolic periodic point 
of a $C^1$-diffeomorphism $f$ satisfies 
property $(\ell)$. Then, there is a $C^1$-open neighborhood $\rO$ of $f$ such that there is 
a transitive filtrating Markov partition $\cU$ and 
disjoint compact set $A$, $B$ such that 
$(f, \cU)$ satisfies the assumption of 
Proposition~\ref{p.C1-non-minimal}.
\end{lemm}
\begin{proof}
For the proof, we only need to repeat the proof of 
Proposotion~\ref{p.C1-non-minimal}. 
Let us explain this. 

First, by assumption we know that $[p]$ contains
a non-trivial homoclinic class in a filtrating 
Markov partition which we denote by $\cU$. Thus by 
Theorem~\ref{t.relative} up to an arbitrarily 
small perturbation there are three different
$\varepsilon$-flexible periodic points $p_0, p_A, p_B$
with large 
stable manifolds where $\varepsilon$ can be chosen 
arbitrarily small, in particular, we may assume that
$\cU$ is $10\varepsilon$-robust. 
Since they are homoclinically related and their 
homoclinic class is non-trivial, 
we may assume that 
\begin{itemize}
\item $p_A$, $p_B$ has a homoclinic point $x_A, x_B$
respectively.
\item $p_A, p_0$ are homoclinically related with 
the heteroclinic orbits $y_{A0}, y_{0A}$.
\item $p_B, p_0$ are homoclinically related with 
the heteroclinic orbits $y_{B0}, y_{0B}$.
\end{itemize}
Then consider the circuit with periodic points 
$p_A, p_B, p_0$ and homo/heteroclinic orbits
$x_A, x_B, y_{A0}, y_{0A}, y_{B0}$ and $y_{0B}$
and denote it by $K$.
We denote the sub circuit consisting of 
$p_A$ and $x_A$ by $K_A$, and 
$p_B$ and $x_B$ by $K_B$.
Since every periodic point of $K$ has a large stable 
manifold, by Lemma~\ref{l.mrk-sml} 
we have $\cU_{(m, n)}(K)$ is very close to $K$.
Thus choosing $m$ and $n$ very large we may assume that 
$\cU_{(m, n)}(K_A)$ and $\cU_{(m, n)}(K_B)$ are disjoint. 

Now we apply Theorem~\ref{t.isola} to $K$ and
$\cU_{(m, n)}(K)$: Then by $2\varepsilon$-perturbation 
we obtain a diffeomorphism $g$ such that 
 there is a transitive filtrating Markov partition 
 $\cU'$ containing $K_g$ (the continuation of $K$).
 Note that $\cU'(K_{Ag}), \cU'(K_{Bg})$ 
 are disjoint sub Markov partitions. Thus we may choose 
 two disjoint compact sets $A, B$ such that 
 $\cU'(K_{Ag}) \subset \mathring{A}$ and 
 $\cU'(K_{Bg}) \subset \mathring{B}$. 
 
 Now we apply Theorem~\ref{t.relative} to 
 $\cU'(K_{Ag})$ and $\cU'(K_{Bg})$:
 Up to a $4\varepsilon$-perturbation we have
 that $\cU'(K_{Ag})$ and $\cU'(K_{Bg})$ satisfy
 property $(\ell_{\cU'(K_{Ag})})$ and
 $(\ell_{\cU'(K_{Bg})})$ respectively, having almost 
 the same circuit $K_g$. 
 Note that at this moment $(\cU', p_0)$ satisfies 
 property $(\ell)$. Thus letting $A$ and $B$ as above
 we obtain the conclusion.
 \end{proof}
\subsection{$C^1$-genericity of non-transitive aperiodic classes}

In this subsection, we complete the 
proof of the locally generic existence of non-transitive 
uniquely ergodic aperiodic classes. 

\begin{prop}\label{p.C1-non-transitive} 
Let $f$ be a diffeomorphism of 
a closed three manifold $M$ and
$A$, $B$ be disjoint compact subsets of $M$. 
Consider the following $C^1$-robust 
property $\cP_4$ for a filtrating set $\cU$:
\begin{itemize}
\item $\cU$ is a transitive filtrating Markov partition 
having unique rectangles contained in $\mathring{A}$
and $\mathring{B}$ which we denote 
by $A_0$ and $B_0$ respectively. 
\item There is a circuit of points $L$ which consists of 
one periodic orbit $\cO(p)$ with 
a large stable manifold in $\cU$ 
and three homoclinic orbits of $p$
which we denote by $\gamma_j$, $j\in \{0,1,2\}$.  
 \item  $L\cap A_0$,  $L\cap B_0$ are 
 singletons. We denote them by $x_1$ and 
 $x_2$ respectively. 
 We have $x_i \in \gamma_i$ for $i=1, 2$.
 \item  Let $L_0$ be 
 the sub circuit of $L$ consisting of 
 $\cO(p)$ and $\gamma_0$. Then 
  $(f, \cU(L_0))$ satisfies property 
  $(\ell_{\cU(L_0)})$. 
\end{itemize}
This property $\cP_4$ 
is type $\fP_{\mathrm{nT}, A, B}$ 
(see Proposition~\ref{p.principle-nontrurgo}).
Namely, for any $\delta_1>0, \delta_2>0$, $N>0$, $M>0$ and any $C^1$-neighborhood $\rO$ of $f$ there is $g\in \rO$ such that:
\begin{itemize}
 \item $\supp (g, f)$ is strictly contained in $\cU$.
 \item There are disjoint, 
 transitive filtrating Markov partitions 
 $\cU_i \subset \mathring{\cU}$ ($i=1, 2$) 
 whose minimum periods are larger than $N$.
 \item $\max \{\cdiam(\cU_i), \cdiam (g(\cU_i))\}<\delta_1$ for $i=1,2$.
 \item The ergodic diameter of $\cU_i$ 
 (see Section~\ref{ss.C1infi}) 
 is less than $\delta_2$ for $i=1,2$.
 \item For $i=1, 2$, the rectangles $A_0$, $B_0$ of $\cU$ contain exactly one rectangle of $\cU_i$, which 
 we denote by $A^i_1$ and $B^i_1$ respectively. 
 \item Any $\cU_i$-chain of components 
 of length $M$ starting from or ending at $A^i_1$
 does not contain $B^i_1$ for $i=1, 2$.
 \item $\cU_i$ satisfies the property $\cP_4$ for $i=1, 2$. Namely, 
\begin{itemize}
 \item There is a circuit $L^i$ consisting of one periodic
 orbit $\cO(p^i)$ and three 
 homoclinic orbits $\gamma^i_j$ ($j=0, 1, 2$) 
 of $\cO(p^i)$.
 \item  $L^i\cap A_0$,  $L^i\cap B_0$ are 
 singletons. We denote them by $x^i_1$ and $x^i_2$. We have $x^i_j \in \gamma^i_j$ for $j=1, 2$.
 \item Let $L^i_0$ be the sub circuit of $L^i$
 consisting of $\cO(p^i)$ and $\gamma^i_0$. 
 Then $(g, \cU_i(L^i_0))$ satisfies property
 $(\ell_{\cU_i(L^i_0)})$.  
\end{itemize}
\end{itemize}
\end{prop}

\begin{proof}
We choose $\varepsilon>0$ such that any $10\varepsilon$-perturbation of $f$ belongs to $\rO$
and the filtrating Markov partition 
$\cU$ is $10\varepsilon$-robust. 

{\bf Step 1: Preparation.}

As $f$ satisfies property $(\ell_{\cU(L_0)})$, 
by Theorem~\ref{t.flex} we know 
there is an arbitrarily small perturbation of $f$ (still denoted by $f$) having two 
distinct $\varepsilon$-flexible 
periodic points $p^i$ ($i=1,2$) whose orbits are 
contained in $\cU(L_0)$, 
with large stable manifolds in $\cU$. 
Note that by assumption, $\cU(L_0)$ does not 
contain $A_0$ and $B_0$

Let us show that there is a circuit $K_i$ consisting of the orbit of $p^i$ and $3$ homoclinic orbits of $p^i$ as follows:
\begin{itemize}
\item $\gamma^i_0$: 
whose orbit is disjoint from $A_0$ and $B_0$,
\item $\gamma^i_1$: whose orbit is disjoint from $B_0$ and passes $A_0$,
 \item $\gamma^i_2$: whose orbit is disjoint from $A_0$ and passes $B_0$.
\end{itemize}
The existence of $\gamma^i_0$ is 
just a consequence of the non-triviality
of the relative homoclinic class of $p^i$ 
in $\cU(L_0)$.

Let us see how to find $\gamma^i_1$. 
The existence of the 
homoclinic orbit $\gamma_1$ in the assumption
tells us there is a $\cU$-chain
of points starting from $\cU(p^i)$, passing $A_0$
only once, ending at $\cU(p^i)$,
and do not pass $B_0$. Thus by Lemma~\ref{l.homoclinic} and Remark~\ref{r.nontri}
(note that the orbit of $p^i$ does not pass $A_0$), 
there is a homoclinic 
orbit $\gamma^i_1$ of $p^i$ which has such an
itinerary. This gives us the desired homoclinic orbit.
Similarly, we can find $\gamma^i_2$.
We put $A_0 \cap \gamma^i_1 = x^i_1$ and 
$B_0 \cap \gamma^i_2 = x^i_2$ for $i= 1, 2$.

Now we apply Theorem~\ref{t.aff} to $K_i$.
Up to an arbitrarily small perturbation, 
for every sufficiently large $m$ and $n$ 
we have that $\cU_{(m, n)}(K_i)$ are affine.
By abuse of notation, we denote the perturbed 
diffeomorphism by $f$, too.
Since $\cdiam (\cU_{(m, n)}(K_i)) \to 0$ 
as $m, n \to \infty$
by Lemma~\ref{l.mrk-sml}, we know that given $N>0,
M>0$, $\delta_1 >0$ and $\delta_2 >0$ by choosing $m$ and $n$ large we have the following conditions:
\begin{itemize}
\item $\cdiam(\cU_{(m, n)}(K_i)), 
\cdiam (g(\cU_{(m, n)}(K_i))) < \delta_1$ for every $g$ which is 
$10\varepsilon$-$C^1$-close to $f$.
\item The minimum period of $\cU_{(m, n)}(K_i))$
is larger than $N$ for $i=1, 2$.
\item Each $\cU_{(m, n)}(K_i)$-chain of components
 of length $M$
starting from or ending at $\cU_{(m, n)}(x^i_1)$ 
does not contain the rectangle 
$\cU_{(m, n)}(x^i_2)$.
\item $\delta_{\mathrm{erg}}(\cdiam(\cU_{(m, n)}(K_i))) <\delta_2$ (see Lemma~\ref{l.ergodiam}).
 \end{itemize}
Let us explain how to obtain the third 
condition. If $m, n$ are large then 
$\cU_{(m, n)}(K_i)$ is very close to $K_i$. 
Thus every $\cU_{(m, n)}(K_i)$-chain of connected 
components starting from $\cU_{(m, n)}(x^i_1)$
needs to be long to reach 
$\cU_{(m, n)}(\cO(p^i))$. Especially, it should be 
long to reach $\cU_{(m, n)}(x^i_2)$ 
as well.

{\bf Step 2: Expulsion.}

 Now Theorem~\ref{t.isola} allows us to perform 
a $2\varepsilon$-perturbation $g$ of $f$ 
supported in 
$\cU_{(m, n)}(K_1)\cup \cU_{(m, n)}(K_2)$ 
such that $g$ admits filtrating transitive 
Markov partitions 
$\cU_i\subset \cU_{(m, n+\nu)}(K_i)$, where 
$i=1,2$ and $\nu$ is some positive integer, 
with the following properties: 
\begin{itemize}
\item $\cU_i$ is matching to
 $\cU_{(m, n+\nu)}(K_i)$ and
$8\varepsilon$-robust.
 \item There are circuits $K'_{i}\subset \cU_i$ 
 which is similar to $K_i$ 
 and whose periodic orbits coincides with 
 $\cO(p^i)$ for $i=1, 2$.
 \item $p^i$ is $\varepsilon$-flexible and has a large stable manifold in $\cU_i$ for $i=1, 2$.
\end{itemize}
Note that the three conditions for $\cU_{(m, n)}(K_i)$ 
(smallness of the connected components, 
non-existence of chains and the smallness of the ergodic diameters)
inherits to $\cU^i$.

Since $\cU_i$ is matching to 
$\cU_{(m, n+\nu)}(K_i)$
and $\cU_{(m, n+\nu)}(K_i)$ 
contains unique rectangle
in $A_0$ (resp. $B_0$),
there is a unique rectangle of $\cU_i$ 
contained in $A_0$ (resp. $B_0$)
which we denote by $A^i_1$ (resp. $B^i_1$).
Note that $K'_{i}\cap A^i_1$ 
(resp. $K'_{i}\cap B^i_1$) is a singleton 
and it is contained in the homoclinic orbits 
which conjugates to $\gamma^i_1$ 
(resp. $\gamma^i_2$).

{\bf Step 3: Recovery.}

Recall that $p^i$ has the third homoclinic orbit 
in $K'_{i}$ which is disjoint from $A^i_1\cup B^i_1$. 
Let $K'_{i, o}$ denote the circuit consisting of 
$\cO(p^i)$ and the homoclinic orbit 
which conjugates to 
$\gamma^i_0$.
Applying Theorem~\ref{t.relative} to  
$\cU^i(K'_{i, o})$, we know that up to 
a $4\varepsilon$-perturbation 
$\cU_i(K'_{i, o})$
satisfies the condition $(\ell_{\cU_i(K'_{i, o})})$.

Note that this last perturbation's support is contained
in $\cU_i$, thus it does not change the set of 
$\cU_i$-chains, and as a result keeps the ergodic 
diameter small.
\end{proof}

Let us complete the proof of the existence 
of non-transitive classes.

\begin{theo}
Let $\rO$ be a $C^1$-open set of diffeomorphisms on a closed $3$-manifold $M$ admitting a transitive filtrating Markov partition $\cU$ containing a periodic point 
$p_f$  varying continuously with $f$ such that
$(f,\cU, p_f)$ satisfies the property $(\ell)$. 

Then there is a $C^1$-residual subset $\rR\subset \rO$ 
such that every $f\in\rG$ has
an uncountable set of chain recurrence classes which are 
not transitive and uniquely ergodic.
\end{theo}

\begin{proof}
We only need to confirm that property $(\ell)$ implies 
the existence of transitive Markov partition satisfying 
the condition required for the ``first-level'' of the 
induction process. Then 
Proposition~\ref{p.C1-non-transitive}
guarantees that $C^1$-generically we have 
uncountably many 
nested sequences of filtrating Markov partitions 
satisfying the assumption of 
Lemma~\ref{l.nontr}, and the result holds. 

Suppose that $(f, \cU)$ satisfies condition $(\ell)$.
First, by Theorem~\ref{t.flex} we know that up 
to an arbitrarily small perturbation and for any 
arbitrarily small $\varepsilon >0$
there is a periodic point $p$ which is 
$\varepsilon$-flexible 
with a large stable 
manifold, having a non-trivial homoclinic class.
Since the homoclinic class is non-trivial, $p$ has
infinitely many homoclinic points whose orbits 
are mutually distinct. 
We choose three of them and call them $x_i$ ($i=0,1,2$).

Now, using the largeness of the stable manifold of $p$,
we choose a sufficiently fine refinement of $\cU$ such 
that there is a unique rectangles
$A$ and $B$ containing $x_1$ and 
$x_2$ respectively. Then consider the circuit $L$ consisting 
of $\cO(p)$ and $\cO(x_i)$ (i=0, 1, 2). We denote the sub-circuit
$L_0$ consisting of $\cO(p)$ and $\cO(x_0)$. 
Then, apply Theorem~\ref{t.relative}
to $\cU(L^0)$ to 
obtain the property $(\ell_{\cU(L^0)})$.
\end{proof}


\vspace{.5cm}

{\bf Christian Bonatti} bonatti@u-bourgogne.fr\\  
Institut de Math\'{e}matiques de Bourgogne, UMR 5584 du CNRS,  Universit\'{e} de
Bourgogne Franche-Comt\'{e},\\ F-21000, Dijon, France.

\vspace{.5cm}
{\bf Katsutoshi Shinohara} ka.shinohara@r.hit-u.ac.jp\\
Graduate School of Business Administration, Hitotsubashi University,
2-1 Naka, Kunitachi, Tokyo 186-8601, Japan.


\begin{thebibliography}{99}

%

\bibitem[AS]{AS} R. Abraham and S. Smale, {\em Nongenericity of $\Om$-stability,\/} Global Analysis I, Proc. Symp. Pure Math A.M.S., {\bf 14}, 5--8 (1968).

\bibitem[Ar]{Ar}
A. Artigue,
\emph{Minimal expansive systems and spiral points,}
Topology and its Applications,
{\bf 194}, 
166--170 (2015).

\bibitem[Bo]{Bo} C. Bonatti, 
\emph{Towards a global view of dynamical systems, 
for the $C^1$-topology,} Ergod. Th. \& Dynam. Sys., 
{\bf 31}(4) (2011), 959--993.

\bibitem[BC]{BC}
C. Bonatti, S. Crovisier,
{\it R\'{e}currence et g\'{e}n\'{e}ricit\'{e}},
Invent. Math., 158 (2004), no. 1, 33--104. 
%
%
%
\bibitem[BD$_1$]{BD1} 
C. Bonatti, L. D\'{i}az,
{\it On maximal transitive sets of generic diffeomorphisms},
Publ. Math. Inst. Hautes \'{E}tudes Sci., no. 96, (2002), 171--197.
%
%
%
\bibitem[BD$_2$]{BD2} Ch. Bonatti and L.J. D\'\i az, {\em Connexions h\'et\'eroclines et g\'en\'ericit\'e d'une infinit\'e de puits ou de sources,\/} Ann. Sci. \'Ecole Norm. Sup., {\bf 32}, 135--150, (1999).
    
\bibitem[BDP]{BDP} 
C. Bonatti, L. D\'{i}az, E. Pujals,
{\it A $C^1$-generic dichotomy for diffeomorphisms: Weak forms of hyperbolicity or infinitely many sinks or sources}, 
Ann. of Math., Volume 158, Number 2 (2003), 355--418.

\bibitem[BDV]{BDV}
C. Bonatti, L. D\'{i}az, M. Viana, 
{\it Dynamics beyond uniform hyperbolicity. A global geometric and probabilistic perspective},
Encyclopaedia of Mathematical Sciences, 102. 
Mathematical Physics, III. Springer-Verlag, Berlin, 2005. xviii+384 pp. 
%

\bibitem[BS$_1$]{BS1} C. Bonatti and K. Shinohara, 
\emph{Flexible periodic points,} 
Ergod. Th. \& Dynam. Sys.,
{\bf 35}(5) (2015), 1394--1422.

\bibitem[BS$_2$]{BS2} C. Bonatti and K. Shinohara, 
\emph{Volume hyperbolicity and Wildness}, 
Ergod. Th. \& Dynam. Sys.,
{\bf 38}(3) (2018), 886--920.

\bibitem[BS$_3$]{BS3} C. Bonatti, K. Shinohara \emph{A mechanism for ejecting a horseshoe from a partially hyperbolic chain recurrence class}, preprint. 

\bibitem[By]{By}
P. Boyland, \emph{Weak disks of Denjoy minimal sets}, 
Ergodic Theory and Dynamical Systems, {\bf 13} (4), 597--614 (1993).

\bibitem[Co]{Co} C. Conley, {\em Isolated invariant sets and Morse index,\/} CBMS Regional Conference Series in Mathematics, {\bf 38}, AMS Providence, R.I., (1978). 

%

\bibitem[GM]{GM}
JM. Gambaudo, M. Martens, 
\emph{Algebraic Topology for Minimal Cantor Sets},
 Ann. Henri Poincar\'{e} {\bf 7}, 
 423--446 (2006).

\bibitem[Gr]{Gr} C. Grillenberger, 
{\em Construction of strictly ergodic systems I. Given entropy,}
Z. Warhrscheinlichkeitstheorie, {\bf 25}, 323--334 (1972).

\bibitem[Ha]{Ha} S. Hayashi, {\em Connecting invariant manifolds and the  solution of the $C^1$-stability and $\Om$-stability conjectures for flows,\/} Ann. of Math., {\bf 145}, 81--137, (1997) and Ann. of Math., {\bf 150}, 353--356, (1999).

\bibitem[Ku]{Ku} Y. Kupka \emph{Contribution à la th\'eorie des champs g\'en\'eriques.}  Contributions to Differential Equations {\bf 2}, 457--484 (1963). 

\bibitem[Ma$_1$]{Ma1} R. Ma\~n\'e, {\em  Contributions to the stability conjecture,\/}  Topology,  {\bf 17}, 386-396, (1978).

\bibitem[Ma$_2$]{Ma2} R. Ma\~n\'e, {\em  An ergodic closing lemma,\/} Ann. of Math.,  {\bf 116}, 503-540, (1982).

\bibitem[Ma$_3$]{Ma-can} R. Ma\~n\'e, 
{\em Expansive homeomorphisms and topological dimension,\/} Trans. of AMS., 
{\bf 252} (1979), 313--319.


\bibitem[Ne$_1$]{Ne} S. Newhouse, \emph{ Nondensity of axiom ${\rm A}({\rm a})$ on $S^{2}$. } (1970)  Global Analysis (Proc. Sympos. Pure Math., Vol. XIV, Berkeley, Calif., 1968)  pp. 191--202 Amer. Math. Soc., Providence, R.I.


\bibitem[Ne$_2$]{N} S.Newhouse, {\em Diffeomorphisms with infinitely many sinks,\/} Topology, {\bf 13}, 9--18, (1974).


\bibitem[Pa]{Pa} J. Palis {\em On the $C^1$ $\Omega$-stability conjecture} Publ. Math. I.H.E.S. {\bf 66}, 211-215, (1988).

\bibitem[Po]{Po} H. Poincar\'{e},
\emph{Les m\'{e}thodes Nouvelles de la M\'{e}canique 
c\'{e}leste}, 1899, new edition by Les grands
classiques Gauthier-Villars, librairie Blanchard, Paris, 1987.

\bibitem[Pu]{Pu} C. Pugh, {\em The closing lemma,\/} Amer. J. Math., {\bf 89}, 956--1009, (1967). 

\bibitem[Ro]{Ro} J. Robbin, {\em A structural stability theorem} Ann. of Math., {\bf 94}, 447-493, (1971).

\bibitem[Ro$_2$]{Ro2} C. Robinson, {\em Structural stability for $C^1$-diffeomorphisms} J. Diff. Equ.{\bf 22}, 28-73, (1976). 

\bibitem[Sa]{Sa}
A.Sannami, 
\emph{The stability theorems for discrete dynamical systems on two-dimensional manifolds}, 
Nagoya Math. J. {\bf 90} (1983), 1--55.

\bibitem[Sh]{Sh} M. Shub,{\em Topological transitive diffeomorphism on $T^4$\/}, Lect.  Notes in Math.,  {\bf 206},   39  (1971).

\bibitem[Sm]{Sm} S. Smale, \emph{Stable manifolds for differential equations and diffeomorphisms.} Ann. Scuola Norm. Sup. Pisa (3) 17 1963 97--116. 

\end{thebibliography}
\end{document}